\documentclass{article}

\usepackage{etex}
\usepackage{graphicx}
\usepackage{latexsym}
\usepackage{amsfonts}
\usepackage{amssymb}
\usepackage{amsmath}
\usepackage{verbatim}
\usepackage{url}
\usepackage[standard, thmmarks]{ntheorem}
\usepackage[margin=3cm]{geometry}
\usepackage{fancyhdr}

\newtheorem{ax}{Axiom}
\newtheorem*{ax0}{Axiom}
\newtheorem{thm}{Theorem}[section]
\newtheorem{cor}[thm]{Corollary}
\newtheorem{lem}[thm]{Lemma}
\newtheorem{prop}[thm]{Proposition}

\newtheorem{defn}[thm]{Definition}

\newcommand{\Z}{\mathbb Z}

\newcommand{\Q}{\mathbb Q}

\newcommand{\F}{\mathcal F}

\newcommand{\M}{\mathcal M}

\newcommand{\Ob}{\text{Ob\,}}

\newcommand{\To}{\longrightarrow}

\newcommand{\W}{\mathcal{W}}

\DeclareMathOperator{\Ab}{Ab} \DeclareMathOperator{\Set}{Set}
\DeclareMathOperator{\Mor}{Mor}

\input xy
\xyoption{all}

\begin{document}

\title{Sutured Floer homology, sutured TQFT and non-commutative QFT}

\author{Daniel V. Mathews}%

\date{}

\maketitle

\begin{abstract}
We define a ``sutured topological quantum field theory'', motivated by the study of sutured Floer homology of product $3$-manifolds, and contact elements. We study a rich algebraic structure of suture elements in sutured TQFT, showing that it corresponds to contact elements in sutured Floer homology. We use this approach to make computations of contact elements in sutured Floer homology over $\Z$ of sutured manifolds $(D^2 \times S^1, F \times S^1)$ where $F$ is finite. This generalises previous results of the author over $\Z_2$ coefficients. Our approach elaborates upon the quantum field theoretic aspects of sutured Floer homology, building a non-commutative Fock space, together with a bilinear form deriving from a certain combinatorial partial order; we show that the sutured TQFT of discs is isomorphic to this Fock space.
\end{abstract}

\tableofcontents

\section{Introduction}

\subsection{Chord diagrams and signs}

\begin{figure}[h]
\centering
\includegraphics[scale=0.5]{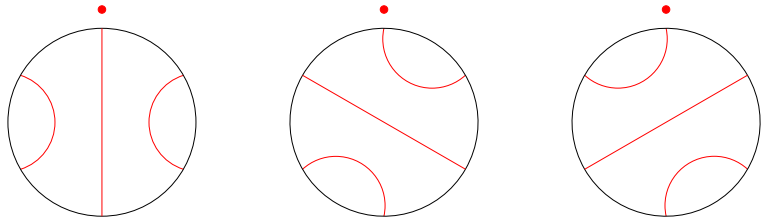}
\caption{A bypass triple.} \label{fig:21}
\end{figure}

This paper, like its prequel \cite{Me09Paper}, is about fun with chord diagrams. A chord diagram $\Gamma$ is a finite collection of non-intersecting properly embedded arcs in a 2-dimensional disc $D^2$. Two chord diagrams are considered equivalent if they are homotopic relative to endpoints. Fixing $2n$ points on $\partial D^2$, there are finitely many chord diagrams of $n$ arcs (or chords) connecting them --- in fact, $C_n$ of them, the $n$'th Catalan number. In \cite{Me09Paper} we considered the $\Z_2$ vector space $SFH_{comb}(T,n)$ generated by chord diagrams of $n$ chords, subject to a relation called the \emph{bypass relation}. The bypass relation says that whenever $3$ chord diagrams coincide, except within a disc on which the chords appear as shown in figure \ref{fig:21}, they sum to zero.

In \cite{Me09Paper}, we showed that this vector space describes precisely contact elements in the sutured Floer homology of manifolds of the form $(D^2 \times S^1, F \times S^1)$, where $F$ is finite, $|F| = 2n$, with $\Z_2$ coefficients. However sutured Floer homology can be defined over $\Z$, and the structure of contact elements should be similar. Therefore, our vector space $SFH_{comb}$ should generalise to $\Z$ coefficients.

Over $\Z$, contact elements have a sign ambiguity \cite{HKM08}. So a chord diagram should represent an element in a $\Z$-module (i.e. abelian group) $V$, \emph{up to sign}, as in the \emph{lax vectors} of \cite{Conway}. The sum of two lax vectors $\pm v, \pm w$ is not well-defined: choosing representatives (lifts) $v,w$, we could have $(\pm v) + (\pm w)$ equal to $\pm (v+w)$ or $\pm (v-w)$. A three-term relation such as the bypass relation says that three lax vectors somehow sum to zero; equivalently, it says that 3 lax vectors $\pm u, \pm v, \pm w$ are related so that $\pm (u+v)$ or $\pm (u-v)$ equals $\pm w$. 

We would like to resolve these ambiguities. We shall do so in this paper, and along the way we shall unearth much extra structure.

\subsection{Stackability resolves signs}

Central to this resolution will be a bilinear form, which we shall denote $\langle \cdot | \cdot \rangle$, on our abelian group $V$. The corresponding bilinear form over $\Z_2$ was called $m(\cdot,\cdot)$ in \cite{Me09Paper}. It can be defined by \emph{stacking} chord diagrams, an operation defined in \cite{Me09Paper} and developed in section \ref{sec:sutured_surfaces} below. We place two chord diagrams as the lids of a cylinder with vertical chords running along its sides, then round corners and chords in a specific way to obtain a sphere with curves on it. The chord diagrams are said to be \emph{stackable} if we obtain a single connected curve on the sphere. We will have a map
\[
\langle \cdot | \cdot \rangle \; : \; V \otimes V \To \Z
\]
which takes stackable chord diagrams to $\pm 1$, and non-stackable chord diagrams to $0$; and which reduces to $m(\cdot, \cdot)$ modulo $2$. The key property we shall use to obtain coherent signs for chord diagrams is that while there is ambiguity in $\pm 1$, there is no ambiguity in $\pm 0 = 0$. The bypass relation tells us that a chord diagram should be given by $\pm (u+v)$ or $\pm (u-v)$, for some lax vectors $\pm u, \pm v$. However if we know, say, that $\langle u | w \rangle = \langle v | w \rangle = 1$ for some $w$, then $\langle \pm (u-v) | w \rangle = 0$ while $\langle \pm (u+v) | w \rangle = \pm 2$; we can then usefully distinguish between the two.

We shall consider specific \emph{basis} chord diagrams, as described in \cite{Me09Paper}; these correspond to basis elements of $V$. We shall find a coherent way to choose signs for all these basis elements, using the above idea, so as to describe contact elements in general.

\subsection{A non-commutative Fock space and partial order}

In quantum field theory a Fock space $\F$ is an algebraic object whose elements can represent states of several particles. We might say that $x \in \F$ represents the presence of one particle, $y \in \F$ the presence of a second particle, and $xy$ the presence of both. A Fock space has creation and annihilation operators which add and remove particles, and which are adjoint with respect to an inner product. In a commutative (bosonic) Fock space, $xy=yx$. In an anticommutative (fermionic) Fock space, $xy=-yx$. We shall consider a noncommutative $\F$, however, in which neither $xy=yx$ nor $xy=-yx$ holds. Roughly speaking, our $\F$ will contain linear combinations of (noncommutative) products of $x$ and $y$.

There are natural creation and annihilation operators on $\F$, which insert or delete a specified symbol in a word in a specified place. There are many more creation and annihilation operators in the noncommutative case than in the commutative or anticommutative cases, and they obey relations similar to a simplicial set. In an appropriately noncommutative way, creation and annihilation operators are adjoint with respect to the bilinear form $\langle \cdot | \cdot \rangle$.

In \cite{Me09Paper} we considered a certain partial order. Consider words on $x$ and $y$ (in \cite{Me09Paper} we used $-$ and $+$). Say that $w_1 \leq w_2$ if $w_2$ can be obtained from $w_1$ by moving some (possibly none) of the $x$'s to the right (equivalently, by moving some of the $y$'s to the left). Thus $xxyy \leq yxyx$ but $xyyx \nleq yxxy$.

Now define a bilinear form $\langle \cdot | \cdot \rangle \; : \; \F \otimes \F \To \Z$ as a boolean version of the partial order $\leq$, as follows. For two words $w_1, w_2$ in $x$ and $y$, let $\langle w_1 | w_2 \rangle = 1$ if $w_1 \leq w_2$ and $0$ otherwise; then extend to $\F$ linearly. We imagine $\langle \cdot | \cdot \rangle$ as a noncommutative version of an inner product. We can think of $\F$ as a ``Fock space of two non-commuting particles''.

We shall prove the following. (A precise version is \ref{thm:mini_main_isomorphism}; a more detailed statement is theorem \ref{thm:main_isomorphism}.)
\begin{thm}
\label{thm:main_thm_rough}
The Fock space and bilinear form defined from $\leq$ are isomorphic to the abelian group $V$ generated by chord diagrams, and bilinear form defined by stacking:
\[
(\F, \langle \cdot | \cdot \rangle) \cong (V, \langle \cdot | \cdot \rangle).
\]
\end{thm}

That is, chord diagrams give elements of $\F$ up to sign; basis chord diagrams give basis elements of $\F$ (i.e. words), up to sign; any three bypass-related chord diagrams give lax elements of $\F$ summing to zero in an appropriate sense; and stackability is described by the boolean version of $\leq$. This was shown in \cite{Me09Paper} mod $2$; we will show it still holds over $\Z$.

An element of $V \cong \F$ corresponding to a chord diagram $\Gamma$ can be expressed in terms of basis elements of $\F$, i.e. words in $x$ and $y$. As considered at length in \cite{Me09Paper}, among the words appearing in such an element, there are well-defined first and last words $w_- \leq w_+$ with respect to $\leq$. In fact, the correspondence $\Gamma \mapsto (w_-, w_+)$ gives a bijection between chord diagrams and pairs of words comparable with respect to $\leq$. We will show the same properties as shown in \cite{Me09Paper} over $\Z_2$ also hold over $\Z$; in fact every chord diagram gives an element in $\F$ which is a linear combination of words with all coefficients $\pm 1$.

Furthermore, it is not difficult to show that $\langle \cdot | \cdot \rangle$ is nondegenerate. Thus there is a \emph{duality} operator $H$ on $\F$ such that for all $u,v$, $\langle u | v \rangle = \langle v | H u \rangle$. Note that for commutative $\langle \cdot | \cdot \rangle$, i.e. $\langle u | v \rangle = \langle v | u \rangle$, we have $H=1$; and in the anticommutative case, $H=-1$. We will show that $H$ is \emph{periodic}, in the sense that some power of $H$ is the identity; this is, in a sense, a generalisation of commutativity and anticommutativity. This is a purely algebraic statement, but a direct algebraic proof is not at all clear; we prove it by showing that $H$ is equivalent to the operation of \emph{rotating} a chord diagram.

We note that the key feature distinguishing our Fock space from the usual sort is the bilinear form $\langle \cdot | \cdot \rangle$, based on the partial order $\leq$. Partial orders appear prominently in the theory of \emph{causal sets} (see e.g. \cite{Bombelli_et_al_87,Sorkin05}), one approach to quantum gravity. But we have not been able to find anything like the above structure in the theoretical physics literature.

\subsection{Sutured Topological Quantum Field Theory}

We have been somewhat vague about what $V$ is above: it is an object built out of chord diagrams and related objects, but related to sutured Floer homology. We now explain what $V$ is.

All of the above is motivated by sutured Floer homology, contact structures and contact elements; but it has become quite independent of its starting point. The idea is that we can build an algebraic structure, with many of the same properties, purely by reference to surfaces (like discs), arcs drawn on them (like chord diagrams), and various topological operations. To a surface $\Sigma$ with some markings $F$ on the boundary, we shall associate an abelian group $V(\Sigma, F)$. When arcs $\Gamma$ are drawn appropriately on the surface --- a set of \emph{sutures} --- we shall associate an element (up to sign) $c(\Gamma)$ in the abelian group $V(\Sigma,F)$, called a \emph{suture element}. These associations are natural with respect to certain topological operations, such as gluing and the stacking operation. Thus we have something similar to a topological quantum field theory, which we shall call \emph{sutured topological quantum field theory}. It will be constructed without reference to any Heegaard decompositions or holomorphic curves. 

We shall define sutured TQFT axiomatically. We shall impose more axioms than merely the associations described above, in order that, at least for discs, sutured TQFT is unique. We regard these axioms as fairly natural, and explain their rationale as we introduce them. A version of the bypass relation, resolving the sign ambiguities mentioned above, will appear in sutured TQFT, but we do not impose it as an axiom. We shall instead define the bilinear form $\langle \cdot | \cdot \rangle$ from stackability (which is nothing but a certain inclusion, or gluing, of surfaces) and impose a certain nondegeneracy condition on it. From this, the bypass relation shall follow.

The abelian group $V$ discussed above is in fact the sutured TQFT of discs $V(D^2)$; the isomorphism of theorem \ref{thm:main_thm_rough} is in fact $V(D^2) \cong \F$. All the structure in sutured TQFT of discs --- including suture elements, various maps arising from various gluings and inclusions, and the stacking map --- corresponds to structure in $\F$. And conversely, all the structure in $\F$ --- including creation and annihilation operators, simplicial structure, and duality operator --- corresponds to structure appearing on discs in sutured TQFT.

As discussed by Honda--Kazez--Mati\'{c} in \cite{HKM08} and at length by the author in \cite{Me09Paper}, the $SFH$ of such product manifolds $(\Sigma \times S^1, F \times S^1)$, and their contact elements, have properties similar to a $(1+1)$-dimensional topological quantum field theory. Our sutured TQFT is a ``pure'' version of this TQFT, abstracted from $SFH$ itself. The point of this paper is that, in describing this TQFT axiomatically, we find surprising properties familiar from (non-topological!) quantum field theory, namely all the structure of the Fock space $\F$. The sutured TQFT of discs is the QFT of two non-commuting particles, where the two usages of ``QFT'' are quite distinct.

We remark that sutured TQFT, as we define it, is very similar in some respects to a \emph{planar algebra} \cite{Jones99}. To a surface with boundary is associated an algebraic object; to curves on the surface, dividing it into positive and negative regions, are associated distinguished elements of the algebraic object; surfaces may be glued together, giving maps of algebraic objects which are natural with respect to the distinguished elements; surfaces with many boundary components give operators. However there are several distinctions: sutured TQFT is not restricted to planar surfaces; sutured TQFT has no canonical form of multiplication (except on discs, where it does not coincide with planar algebra multiplication); the element associated to a set of curves in sutured TQFT has a sign ambiguity; and in the author's limited knowledge, nothing like a bypass relation or inner product has been studied in the context of planar algebras. Clearly the two subjects have wildly disparate motivations: planar algebras being motivated by the study of subfactors and von Neumann algebras; and sutured TQFT being motivated by contact geometry and sutured Floer homology. We shall leave these connections for future investigation, but they are striking, and we wonder how deep they are.

\subsection{$SFH$ gives a sutured TQFT}

Although $V$ is defined from the axioms of sutured TQFT, it is motivated by sutured Floer homology; indeed it is designed to be isomorphic to the $SFH$ of certain manifolds. We will show that the sutured Floer homology of product manifolds $SFH(\Sigma \times S^1, F \times S^1)$, where $\Sigma$ is a surface with boundary and $F \subset \partial S$ is finite, forms a sutured TQFT. Since sutured TQFT on discs is unique, any sutured TQFT on discs is isomorphic to $SFH(D^2 \times S^1, F \times S^1)$. However sutured TQFT on more complicated surfaces is not unique, and can be zero at higher genus. But we also consider an additional axiom, also satisfied by $SFH(\Sigma \times S^1, F \times S^1)$, which ensures nontriviality at higher genus.

The upshot of this paper, then, is that three structures are equivalent, in the case of discs: sutured Floer homology $SFH(D^2 \times S^1, F \times S^1)$, sutured TQFT $V(D^2)$, and the non-commutative Fock space $\F$.

Moreover, all the structure of $SFH$ of product manifolds $(\Sigma \times S^1, F \times S^1)$ (a very restricted class of manifolds, to be sure) can be described without considering Heegaard decompositions, holomorphic curves, or contact structures --- all of which are involved subjects and essential to the definition of contact elements in sutured Floer homology. We obtain proofs about contact elements which are both ``holomorphic curve free'' and ``contact geometry free''. In a subsequent paper we shall use sutured TQFT to give a proof that the contact element of a torsion contact structure is zero \cite{GHV, Massot09}.

\subsection{What this paper does}

The above remarks give an introduction to the main ideas in this paper, but are not a complete description. This paper performs several tasks, as follows.

First, it makes computations in sutured Floer homology. Namely, it classifies contact elements in the sutured Floer homology, with $\Z$ coefficients, of sutured manifolds $(D^2 \times S^1, F \times S^1)$, where $F \subset \partial D^2$ is finite. This generalises the results over $\Z_2$ in \cite{Me09Paper}. Thus we are able to extend our results from that paper regarding the structure of contact elements to $\Z$ coefficients.

Second, it elaborates greatly upon the quantum-field-theoretic aspects of sutured Floer homology and contact elements for sutured manifolds of the form $(\Sigma \times S^1, F \times S^1)$; in particular, it proves that when $\Sigma$ is a disc, $SFH$ is isomorphic to the Fock space $\F$, which is a formal algebraic model of a certain (extremely simple) non-commutative QFT, with creation and annihilation operators, a non-commutative ``inner product'', and more. In fact, because the formal algebraic structure of $\F$ is easiest to define, we shall describe $\F$ precisely and deduce various properties of it --- and then show that $SFH$ is isomorphic to it. 

Third, it abstracts from $SFH$ and describes it axiomatically as sutured TQFT. These axioms constrain sutured TQFT on discs to be isomorphic to the Fock space $\F$ of two non-commuting particles, although (at least with the axioms as we state them) sutured TQFT is not unique on more complicated surfaces. We show that the sutured Floer homology of manifolds $(\Sigma \times S^1, F \times S^1)$, together with contact elements, forms a sutured TQFT.

Moreover, we may eschew sutured Floer homology until the last minute; it is not necessary to the discussion of sutured TQFT or non-commutative QFT, although we do use it to prove that a sutured TQFT exists! In a sense, our SFH results are more general than SFH: they are about suture elements in sutured TQFT, of which contact elements in SFH of product manifolds form an example.

\subsection{Structure of this paper}

The above considerations determine the structure of this paper. First  we first establish the formal algebraic structure of our non-commutative quantum field theory (section \ref{sec:algebra}). Then we axiomatically introduce sutured topological quantum field theory (section \ref{sec:sutured_TQFT}), and deduce many properties of it. We show that $SFH(\Sigma \times S^1, F \times S^1)$ forms a sutured TQFT (section \ref{sec:SFH}). And then we demonstrate the isomorphism between sutured TQFT of discs, and non-commutative QFT (section \ref{sec:sutured_TQFT_discs}), and prove some more properties of this structure.

Thus, we begin this paper by describing, in great detail, the non-commutative algebra of $\F$ (section \ref{sec:algebra}). Much of the detail in this section (especially the bi-simplicial structure \ref{sec:bi-simplicial}, differentials, commutation relations, normal form \ref{sec:differentials_commutation_normal_form} and Temperley--Lieb representation \ref{sec:Temperley-Lieb}) can be skipped on a first reading: the important details are the Fock space $\F$, the creation and annihilation operations, the bilinear form $\langle \cdot | \cdot \rangle$, and the duality operator $H$. Similarly, other sections on sutured TQFT (variations on axioms, the isomorphism between duality and rotation, periodicity) are also quite technical.

\subsection{Acknowledgments}

This paper was written during the author's visit to the Mathematical Sciences Research Institute in March 2010, and during the author's postdoctoral fellowship at the Universit\'{e} de Nantes, supported by the ANR grant ``Floer power''.

\section{Algebraic non-commutative QFT}
\label{sec:algebra}

\subsection{Fock space}

Let $S$ be a set with two elements, $S = \{x,y\}$. Let $\M$ be the free monoid on $S$, $\M = S^* = \{x,y\}^*$, i.e. the set of all finite words (including the empty word, which is the identity) on $\{x,y\}$, under the operation of concatenation.

Let $\F$ be the monoid ring of $\M$ over $\Z$. That is, $\F$ consists of finite $\Z$-linear combinations of finite words on $\{x,y\}$; multiplication is concatenation, now extended linearly. Alternatively, $\F$ is the polynomial ring $\F = \Z[x,y]$ generated over $\Z$ by two non-commuting indeterminates $x,y$. The empty word, denoted $1$, is a multiplicative identity. 

Clearly $\M$ and $\F$ have several gradings: degree $n_x$ in $x$, degree $n_y$ in $y$, and linear combinations of these. Clearly multiplication adds these gradings, making $\F$ into a bi-graded ring. Let $\M_n$ denote the subset of $\M$ consisting of words of length $n$, i.e. with total degree $n$, and $\F_n$ the additive subgroup generated by $\M_n$, i.e. linear combinations of words of length $n$. Thus, as graded abelian group, $\F = \oplus_{n \geq 0} \F_n$ and $\F_n = \left( \Z x \oplus \Z y \right)^{\otimes n}$. Denote by $\M_{n_x, n_y}$ the subset of $\M$ consisting of words of degree $n_x, n_y$ in $x,y$ respectively. Let $\F_{n_x, n_y}$ be the additive subgroup of $\F$ generated by $\M_{n_x, n_y}$, i.e. linear combinations of words of degree $n_x, n_y$ in $x,y$ respectively. 

As an alternative notation, thinking of $x$ as having degree $-1$ and $y$ as having degree $1$, let $\M_n^e$ denote the subset of $\M$ consisting of words of length $n$ and degree $e$, i.e. such that $n_y - n_x = e$. Then $\F_n^e$ is the additive subgroup of $\F$ generated by $\M_n^e$. So $n,e,n_x,n_y$ are related by
\[
n = n_x + n_y, \quad e = -n_x + n_y, \quad n_x = \frac{n-e}{2}, \quad n_y = \frac{n+e}{2}.
\]
And we have, obviously,
\[
\M_{n_x, n_y} = \M_{n_x + n_y}^{n_y - n_x}, \quad \M_n^e = \M_{\frac{(n-e)}{2}, \frac{n+e}{2}}, \quad
\F_{n_x, n_y} = \F_{n_x + n_y}^{n_y - n_x}, \quad \F_n^e = \F_{\frac{n-e}{2}, \frac{n+e}{2}}.
\]
Further, clearly,
\[
\M = \bigsqcup_n \M_n = \bigsqcup_{n,e} \M_n^e = \bigsqcup_{n_x,n_y} \M_{n_x,n_y}, \quad 
\F = \bigoplus_n \F_n = \bigoplus_{n,e} \F_n^e = \bigoplus_{n_x, n_y} \F_{n_x, n_y},
\]
where $n$ varies over all non-negative integers, and $e \in \Z$ satisfies $|e| \leq n$, and $e \equiv n$ mod $2$. As abelian groups
\[
\F_{n} \cong \Z^{2^n} \quad \text{and} \quad \F_n^e \cong \Z^{\binom{n}{n_x}} = \Z^{\binom{n}{n_y}} = \Z^{\binom{n}{(n-e)/2}} = \Z^{\binom{n}{(n+e)/2}}.
\]

\subsection{Creation and annihilation operators}
\label{sec:algebraic_creation_annihilation}

On $\F$ there we can define many simple operations $\F \To \F$, which are abelian group ($\Z$-module) homomorphisms. In fact any operation on words (i.e. a function from the free monoid $\M$ to itself) induces a linear operator (abelian group endomorphism) on $\F$. Better, any function $\M \To \M \cup \{0\}$ induces a linear operator on $\F$, regarding $0$ as the zero in $\F$. Let the following operations $\F \To \F$ be induced from functions on words $w \in \M$ as follows. Here $s$ denotes a letter in $\{x,y\}$.
\begin{enumerate}
\item
\begin{enumerate}
\item \emph{Initial annihilation}
$a_{s,0} w$: if $w$ begins with an $s$, delete it; else return $0$.
\item \emph{Internal annihilation}
$a_{s,i} w$, for $1 \leq i \leq n_s w$: delete the $i$'th $s$ in $w$.
\item \emph{Final annihilation}
$a_{s,n_s w + 1} w$: if $w$ ends with an $s$, delete it; else return $0$.
\end{enumerate}
\item
\begin{enumerate}
\item \emph{Initial creation}
$a^*_{s,0} w$: prepend an $s$ to the beginning of $w$.
\item \emph{Internal creation}
$a^*_{s,i} w$, for $1 \leq i \leq n_s w$: replace the $i$'th $s$ in $w$ with $ss$.
\item \emph{Final creation}
$a^*_{s, n_s + 1} w$: append an $s$ to the end of $w$.
\end{enumerate}
\end{enumerate}

The first three of these are collectively called \emph{annihilation} operators. The names initial, internal, and final should be clear; initial and final annihilation are collectively called \emph{terminal} annihilation. The last three are collectively \emph{creation} operators; we have final, internal, initial, and terminal creations. On each $\F_n^e$ there are $n+4$ creation operators and $n+4$ annihilation operators, each of which map to some $\F_{n \pm 1}^{e \pm 1}$:
\[
\begin{array}{cc}
a_{x,i} \; : \; \F_n^e \To \F_{n-1}^{e+1}, & a_{y,i} \; : \; \F_n^e \To \F_{n-1}^{e-1}, \\
a^*_{x,i} \; : \; \F_n^e \To \F_{n+1}^{e-1}, & a^*_{y,i} \; : \; \F_n^e \To \F_{n+1}^{e+1}.
\end{array}
\]
Here we take $\F_{-1}^{e} = \{0\}$: on the empty word (identity) there are four annihilation operators, all of which give $0$. The empty word, or identity, we can also call the \emph{vacuum}.

These operators satisfy various relations which are easily checked. Thinking of $x$ and $y$ as two different species of particle, we consider separately the relations between annihilation and creation which are inter- and intra-species:
\begin{enumerate}
\item
\emph{Between $x$ and $y$; inter-species.}
\begin{enumerate}
\item 
In almost every case, $x$-annihilation/creation and $y$-creation/annihilation commute. That is, for $0 \leq i \leq n_x + 1$ and $0 \leq j \leq n_y + 1$, except for $(i,j) = (0,0)$ or $(n_x + 1, n_y + 1)$:
\begin{align*}
a_{x,i} \circ a_{y,j} = a_{y,j} \circ a_{x,i}, &\quad
a^*_{x,i} \circ a_{y,j} = a_{y,j} \circ a^*_{x,i}, \\
a_{x,i} \circ a^*_{y,j} = a^*_{y,j} \circ a_{x,i}, &\quad
a^*_{x,i} \circ a^*_{y,j} = a^*_{y,j} \circ a^*_{y,j}.
\end{align*}
\item
Initial $x$-annihilation/creation and initial $y$-annihilation/creation do not commute:
\begin{align*}
a_{x,0} \circ a_{y,0} \neq a_{y,0} \circ a_{x,0}, &\quad
a^*_{x,0} \circ a_{y,0} \neq a_{y,0} \circ a^*_{x,0}, \\
a_{x,0} \circ a^*_{y,0} \neq a^*_{y,0} \circ a_{x,0}, &\quad
a^*_{x,0} \circ a^*_{y,0} \neq a^*_{y,0} \circ a^*_{x,0}.
\end{align*}
\item
Final $x$-annihilation/creation and final $y$-annihilation/creation do not commute:
\begin{align*}
a_{x, n_x + 1} \circ a_{y, n_y + 1} \neq a_{y, n_y + 1} \circ a_{x, n_x + 1}, &\quad
a^*_{x, n_x + 1} \circ a_{y, n_y + 1} \neq a_{y, n_y + 1} \circ a^*_{x, n_x + 1}, \\
a_{x, n_x + 1} \circ a^*_{y, n_y + 1} \neq a^*_{y, n_y + 1} \circ a_{x, n_x + 1}, &\quad
a^*_{x, n_x + 1} \circ a^*_{y, n_y + 1} \neq a^*_{y, n_y + 1} \circ a^*_{x, n_x + 1}.
\end{align*}
\end{enumerate}
\item
\emph{Among $x$ or among $y$; intra-species.} Set $s$ to be $x$ or $y$; we now only consider annihilations and creations of the same type.
\begin{enumerate}
\item
Annihilations commute after a shift. For $0 \leq i < j \leq n_s + 1$ (not $i=j$):
\[
a_{s,i} \circ a_{s,j} = a_{s,j-1} \circ a_{s,i}.
\]
\item
Annihilations and creations usually commute, with a possible shifting, and are sometimes inverses of each other. For $0 \leq i,j \leq n_s + 1$:
\[
a_{s,i} \circ a^*_{s,j} = \left\{ \begin{array}{cc} 
				      a^*_{s,j-1} \circ a_{s,i} & i < j \\
				      1				& i=j,j+1 \\
				      a^*_{s,j} \circ a_{s,i-1} & i > j+1.
				\end{array} \right.
\]
\item
Creations commute, after a shift. For $0 \leq i \leq j \leq n_s + 1$:
\[
a^*_{s,i} \circ a^*_{s,j} = a^*_{s,j+1} \circ a^*_{s,i}.
\]
\end{enumerate}
\end{enumerate}

\subsection{Bi-simplicial structure}
\label{sec:bi-simplicial}

As noted in \cite{Me09Paper}, the intra-species relations are identical to those of a \emph{simplicial set}. Indeed we have \emph{two} simplicial structures on $\F$, one for $x$ and one for $y$. Under the simplicial structure for $x$, words with degree $n_x$ in $x$ have $n_x + 2$ annihilation operators ($n_x$ internal and $2$ terminal), regarded as face maps, and $n_x + 2$ creation operators, regarded as degeneracy maps; thus words with degree $n_x$ in $x$ should be regarded as $(n_x + 1)$-dimensional degenerate simplices. Similarly, under the $y$-simplicial structure, words with degree $n_y$ in $y$ have $n_y + 2$ annihilation and creation operators, regarded respectively as face and degeneracy maps, and should be regarded as $(n_y+1)$-dimensional degenerate simplices.

Thus $\F$ has the structure of a (bi-)simplicial object in the category of abelian groups, and there are natural contravariant functors
\[
\mathbb{F}_x, \mathbb{F}_y \; : \; \Delta \To \Ab
\]
with ``image $\F$'', $\F = \oplus_n \mathbb{F}_x ( {\bf n} ) = \oplus_n \mathbb{F}_y ( {\bf n} )$, which we now describe. Recall the simplicial category $\Delta$ has as objects the (set-theoretic) non-negative integers ${\bf n} = \{0, 1, \ldots, n-1\}$ and as morphisms order-preserving functions ${\bf m} \to {\bf n}$ between them. $\Ab$ is the category of abelian groups. 

In fact not just $\F$, the monoid ring on the free monoid $\M$, but the free monoid $\M$ itself has two simplicial structures, arising from functors
\[
\mathbb{F}_x, \mathbb{F}_y \; : \; \Delta \To \Set,
\]
with ``image $\M$'', $\sqcup_n \mathbb{F}_x ( {\bf n} ) = \sqcup_n \mathbb{F}_y ( {\bf n} ) = \M \cup \{0\}$. (We must adjoin $0$ since terminal annihilation operators may return $0$. The zero is different from the empty word/vacuum $1$.) Since the elements of the monoid $\M$ form a free basis for $\F$ as an abelian group, and functions between them extend to abelian group endomorphisms of the monoid group $\F$, the functors $\mathbb{F}_x, \mathbb{F}_y$ to $\Set$ extend linearly to the functors to $\Ab$.

The functor $\mathbb{F}_x$ can be described by regarding a word $w$ on $\{x,y\}$ as a function. In fact it will be useful to describe such words by functions in several different ways, and we now pause to describe these functions.
\begin{defn}
Let $w \in \M_{n_x, n_y}$, and number the $x$'s (resp. $y$'s) in $w$ left to right from $1$ to $n_x$ (resp. $n_y$). Define six functions $f_w^x, f_w^y, g_w^x, g_w^y, h_w^x, h_w^y$ as follows.
\begin{enumerate}
\item
$f_w^x: \{1, 2, \ldots, n_x\} \To \{0, 1, \ldots, n_y\}$, $f_w^x(i)$ is the number of $y$'s (strictly) to the left of the $i$'th $x$ in $w$. 
\item
$f_w^y: \{1, 2, \ldots, n_y\} \To \{0, 1, \ldots, n_x\}$, $f_w^y(i)$ is the number of $x$'s (strictly) to the left of the $i$'th $y$ in $w$.
\item
$g_w^x: \{1, 2, \ldots, n_x + n_y\} \To \{0, 1, \ldots, n_x\}$, $g_w^x(i)$ is the number of $x$'s in $w$ up to (and including) the $i$'th symbol. 
\item
$g_w^y: \{1, 2, \ldots, n_x + n_y\} \To \{0, 1, \ldots, n_y\}$, $g_w^y(i)$ is the number of $y$'s in $w$ up to (and including) the $i$'th symbol.
\item
$h_w^x: \{1, 2, \ldots, n_x\} \To \{1, 2, \ldots, n_x + n_y\}$, $h_w^x(i)$ is the position of the $i$'th $x$ in $w$.
\item
$h_w^y: \{1, 2, \ldots, n_y\} \To \{1, 2, \ldots, n_x + n_y\}$, $h_w^y(i)$ is the position of the $i$'th $y$ in $w$.
\end{enumerate}
\end{defn}
For $s \in \{x,y\}$, clearly $f_w^x, f_w^y$ are increasing. The $g_w^s$ are slowly increasing, $g_w^s (i+1) - g_w^s (i) \in \{0,1\}$, and $g_w^s(1) \in \{0,1\}$, $g_w^s (n_x + n_y) = n_s$. (Note $g_w^s(i)$ is like a baseball team's score after $i$ innings, as described in \cite{Me09Paper}.) The $h_w^s$ are strictly increasing. It's clear that any of these functions corresponds to a unique word $w$; moreover, there is a bijection between $\M_{n_x, n_y}$ and functions with these respective properties. These functions are clearly related in many ways; this is clear since any one determines a word, and hence determines all the others. For instance, $h_w^s (i) = f_w^s (i) + i$; the images of $h_w^x$ and $h_w^y$ form a partition of $\{1, \ldots, n_x + n_y\}$; and $h_w^s (i) = \min \{ (g_w^s)^{-1}(i) \}$; also $g_w^x (i) + g_w^y (i) = i$.

For the moment we only need the $f_w^s$. As the $f_w^s$ are increasing, i.e. order-preserving, they can be regarded as morphisms in $\Delta$; after shifting elements of sets to obtain set-theoretic integers ${\bf n} = \{0, \ldots, n-1\}$, we have $f_w^x \in \Mor_\Delta( {\bf n_x}, {\bf n_y + 1} )$; moreover as noted $\M_{n_x, n_y}$ and $\Mor_\Delta( {\bf n_x}, {\bf n_y + 1} )$ are bijective.

The idea of $\mathbb{F}_x$ is to take each object ${\bf n}$ in $\Delta$ to the set of all words in $\M$ with degree $n$ in $x$. From the foregoing such words are in bijection with $\sqcup_{n_y} \Mor_\Delta( {\bf n}, {\bf n_y + 1})$, which we can denote $\Mor_\Delta( {\bf n}, \cdot )$. And then $\mathbb{F}_x$ takes the morphism $g: {\bf m} \to {\bf n}$ of $\Delta$ to the function
\[
\Mor_\Delta \left( {\bf n}, \cdot \right) \stackrel{\mathbb{F}_x g}{\To} \Mor_\Delta \left( {\bf m}, \cdot \right)
\]
given by pre-composition with $g$. Note $\mathbb{F}_x (g)$ takes each $\Mor_\Delta ( {\bf n}, {\bf n_y + 1} )$ to $\Mor_\Delta ( {\bf m}, {\bf n_y + 1} )$, hence preserves $n_y$ (which makes sense as this is the structure of adding and deleting $x$'s).

However this structure does not allow for terminal creation and annihilation operators $a^*_{x,0}, a^*_{x,n_x+1}$. Hence we modify the basic idea of the above description a little. The idea is to append and prepend an $x$ to the beginning and end of each word, and then the above gives us the structure we want. Define the subset $\Mor_\Delta^T ({\bf m}, {\bf n})$ of $\Mor_\Delta ({\bf m}, {\bf n})$ to be those order-preserving maps ${\bf m} \to {\bf n}$ which take $0 \mapsto 0$ and $m-1 \mapsto n-1$, i.e. ``terminals to terminals''. Call them \emph{terminal-preserving morphisms}. It's clear that there is a natural bijection $\Mor_\Delta ({\bf m}, {\bf n}) \cong \Mor_\Delta^T ( {\bf m + 2}, {\bf n})$ given by shifting a map by one and setting its values on $0$ and $m+1$. It's also clear that the composition of a morphism in $\Mor_\Delta^T ({\bf m}, {\bf n})$ and a morphism in $\Mor_\Delta^T({\bf n}, {\bf k})$ is a morphism in $\Mor_\Delta^T({\bf m}, {\bf k})$: terminal-preserving morphisms are closed under composition.

Words of degree $n_x, n_y$ in $x,y$ are in bijective correspondence with words of degree $n_x + 2, n_y$ in $x,y$ which begin and end with $x$. Similarly, non-decreasing functions ${\bf n_x} \to {\bf n_y + 1}$ are in bijective correspondence with $\Mor_\Delta^T ({\bf n_x + 2}, {\bf n_y + 1})$. We define $\mathbb{F}_x$ to take ${\bf n} \in \Ob \Delta$ to the set of all words of degree $n-2$ in $x$, along with $0$. That is,
\begin{align*}
\mathbb{F}_x \left( {\bf n} \right) &= \{0\} \cup \left\{ \text{words of $x$-degree $n-2$ in $\M$} \right\}
\cong \{0\} \cup \left\{ \begin{array}{c} \text{words of $x$-degree $n$ in $\M$} \\ \text{which begin and end with $x$} \end{array} \right\} \\
&\cong \{0\} \cup \bigsqcup_{n_y} \Mor_\Delta ({\bf n-2}, {\bf n_y + 1})
\cong \{0\} \cup \bigsqcup_{n_y} \Mor_\Delta^T \left( {\bf n}, {\bf n_y + 1} \right) \\
&\cong \{0\} \cup \Mor_\Delta^T \left( {\bf n}, \cdot \right)
\end{align*}
For $n=0$ or $n=1$ then, $\mathbb{F}_x ({\bf n}) = \{0\}$. We then define $\mathbb{F}_x$ to take a morphism $g \in \Mor_\Delta({\bf m}, {\bf n})$ to the function
\[
\{0\} \cup \Mor_\Delta^T \left( {\bf n }, \cdot \right) \stackrel{\mathbb{F}_x g}{\To} \{0\} \cup \Mor_\Delta^T \left( {\bf m}, \cdot \right)
\]
which takes $0 \mapsto 0$ and which pre-composes functions by $g$, if such pre-composition gives a terminal-preserving morphism; else gives $0$. (Note that if $g$ is terminal-preserving then such a composition is certainly terminal-preserving; if $g$ is not terminal-preserving then the composition may or may not be terminal-preserving.) This $\mathbb{F}_x g$ preserves the degree in $y$, unless it maps to $0$.

This gives the contravariant functor $\mathbb{F}_x \; : \; \Delta \To \Set$ with image $\M$ as described; this functor then extends to one $\Delta \To \Ab$ with image $\F$ as described. 

For each $n_x \geq 0$ and $0 \leq i \leq n_x + 1$, we define $a_{x,i}$ to be the unique morphism in $\Mor_\Delta ({\bf n_x + 1}, {\bf n_x + 2})$ which has image $({\bf n_x + 2}) \backslash \{i\} = \{0, 1, \ldots, i-1, i+1, \ldots, n_x + 1\}$, i.e. which takes $0 \mapsto 0, \; 1 \mapsto 1, \ldots, i-1 \mapsto i-1, \; i \mapsto i+1, \ldots, n_x \mapsto n_x + 1$. Such an $a_{x,i}$, under $\mathbb{F}_x$, gives a map $\mathbb{F}_x ({\bf n_x + 2}) \to \mathbb{F}_x ({\bf n_x + 1})$, i.e.
\[
\{0\} \cup \left\{ \text{words of $x$-degree $n_x$ in $\M$} \right\} \To \{0\} \cup \left\{ \text{words of $x$-degree $n_x-1$ in $\M$} \right\},
\]
which one easily sees that, extended linearly to $\F$, is our $a_{x,i}$ as originally defined. Note that the terminal annihilation operators $a_{x,0}, a_{x,n_x+1}$ are not terminal-preserving (a terminal is annihilated rather than preserved!) but all other $a_{x,i}$ are terminal-preserving; hence $a_{x,0}, a_{x,n_x+1}$ may sometimes return zero, but the other $a_{x,i}$ do not.

Similarly, for each $n_x \geq 0$ and $0 \leq i \leq n_x + 1$, let $a_{x,i}^*$ be the unique morphism in $\Mor_\Delta ( {\bf n_x + 3}, {\bf n_x + 2} )$ which is surjective and takes the value $i$ twice, i.e. $0 \mapsto 0, \; 1 \mapsto 1, \ldots, i \mapsto i, \; i+1 \mapsto i, \ldots, n_x + 2 \mapsto n_x + 1$. Then $\mathbb{F}_x (a_{x,i}^*) \; : \mathbb{F}_x ({\bf n_x + 2}) \to \mathbb{F}_x ({\bf n_x + 3})$ is a map
\[
\{0\} \cup \left\{ \text{words of $x$-degree $n_x$ in $\M$} \right\} \To \{0\} \cup \left\{ \text{words of $x$-degree $n_x + 1$ in $\M$} \right\}
\]
and again, extended to $\F$, is $a_{x,i}^*$ as originally defined. As each $a_{x,i}^*$ is terminal-preserving, $0 \mapsto 0$ but nothing else maps to $0$.

The functor $\mathbb{F}_y$ is defined similarly, but reversing the roles of $x$ and $y$, and using the increasing functions $f_w^y$. The two functions $f_w^x$ and $f_w^y$ give two ``dual'' ways of looking as a word as an order-preserving map.

\subsection{Partial order}
\label{sec:partial_order}

Suppose we have two words $w_0, w_1$, and consider the functions $f_{w_i}^s, g_{w_i}^s, h_{w_i}^s$, for $i \in \{0,1\}$ and $s \in \{x,y\}$. Inequalities on all these functions are essentially equivalent.
\begin{lem}
Let $w_0, w_1 \in \M$ have the same $x$- and $y$-degree. The following inequalities are all equivalent:
\[
f_{w_0}^x \leq f_{w_1}^x, \quad f_{w_0}^y \geq f_{w_1}^y, \quad 
g_{w_0}^x \geq g_{w_1}^x, \quad g_{w_0}^y \leq f_{w_1}^y, \quad
h_{w_0}^x \leq h_{w_1}^x, \quad h_{w_0}^y \geq h_{w_1}^y.
\]
\end{lem}

\begin{Proof}
The inequalities on $f_{w_i}^s$ and $h_{w_i}^s$ are clearly equivalent, since $h_w^s (i) = f_w^s (i) + i$. The inequalities on $g_{w_i}^s$ and $h_{w_i}^s$ are equivalent since the $g_{w_i}^s$ are slowly increasing and $h_w^s (i) = \min \{ (g_w^s)^{-1}(i) \}$. The inequalities on $h_{w_i}^x$ and $h_{w_i}^y$ are equivalent since their images are complementary.
\end{Proof}

If any (hence all) of these inequalities hold, we say that $w_0 \leq w_1$. It's clear this gives a partial order on $\M$; $\leq$ only relates words in the same $\M_{n_x, n_y}$. It's clear that $\leq$ is a sub-order of the lexicographic (total) ordering on each $\M_{n_x,n_y}$ ($x$ comes before $y$).

It's obvious that if $f \leq g$, then for any order-preserving function $h$, $h \circ f \leq h \circ g$. Recall that the creation operators $a_{s,i}^*$ are terminal-preserving morphisms ${\bf n_s + 3} \to {\bf n_s + 2}$, and under $\mathbb{F}_s$ map $\Mor_\Delta^T({\bf n_s + 2}, \cdot) \to \Mor_\Delta^T({\bf n_s + 3}, \cdot)$ by pre-composition. Hence
\[
  w_0 \leq w_1 \quad \text{implies} \quad a_{s,i}^* w_0 \leq a_{s,i}^* w_1
\]
and in fact the converse is also true. 

The same is not true for annihilation operators. For one thing, terminal annihilations may map to zero, so that for instance $a_{s,i} w_0$ might be a word but $a_{s,i} w_1$ zero, and no comparison possible. Even in the nonzero case we may have $w_0 \nleq w_1$ but $a_{s,i} w_0 \leq a_{s,i} w_1$; e.g. $yxxy \nleq xyyx$ but, applying $a_{1,x}$ to both sides, $yxy \leq yyx$. It is however true that if $w_0 \leq w_1$ and $a_{s,i} w_0 \neq 0$, $a_{s,i} w_1 \neq 0$ then $a_{s,i} w_0 \leq a_{s,i} w_1$.

Later on (section \ref{sec:variations_of_nondegeneracy}) we shall to examine this partial order in more detail; we will need a notion of \emph{difference} between words, and a notion of \emph{minimum} and \emph{maximum} of two words.

\begin{defn}
\label{defn:difference}
Given two words $w_0, w_1 \in \M_{n_x,n_y}$ (comparable or not), the \emph{difference} between $w_0$ and $w_1$ is
\[
d(w_0, w_1) = \sum_{i=1}^{n_x} h_{w_1}^x (i) - h_{w_0}^x (i).
\]
\end{defn}
Since $h_w^x (i)$ gives the position of the $i$'th $x$ in $w$, $h_{w_1}^x (i) - h_{w_0}^x (i)$ gives the difference in position between the $i$'th $x$'s in $w_0$ and $w_1$. If we regard $x$'s as pawns and $y$ as empty squares on a $1 \times n$ chessboard, pawns moving left to right, then $d(w_0, w_1)$ is the number of signed pawn moves required to go from $w_0$ to $w_1$. If $w_0 \leq w_1$, then only left-to-right moves are required (all terms in the sum are positive) and $d(w_0, w_1)$ is the number of such moves.

\begin{lem}
\label{lem:word_decomposition}
For any $w_0, w_1 \in \M_{n_x,n_y}$ (comparable or not), there are decompositions
\[
w_0 = w_0^0 w_0^1 \cdots w_0^{2k-1}, \quad w_1 = w_1^0 w_1^1 \cdots w_1^{2k-1}
\]
where each $w_0^i, w_1^i$ have the same numbers of each symbol, i.e. $w_0^i, w_1^i \in \M_{n_x^i, n_y^i}$ for some $n_x^i, n_y^i$, and such that $w_0^i \leq w_1^i$ if $i$ is even and $w_1^i \leq w_0^i$ if $i$ is odd. The words $w_0^0, w_1^0, w_0^{2k-1}, w_1^{2k-1}$ might be empty, but the other $w_j^i$ are not.
\end{lem}

\begin{Proof}
Recall $g_w^s(i)$ is the number of instances of the symbol $s$ in $w$ up to and including the $i$'th position; and $g_w^x (i) + g_w^y (i) = i$. In particular, $g_{w_0}^x (i) = g_{w_1}^x (i)$ iff $g_{w_0}^y (i) = g_{w_1}^y (i)$, iff up to the $i$'th position, both $w_0$ and $w_1$ have the same number of $x$'s and $y$'s. We split the words $w_0, w_1$ at such locations.

The difference $g_{w_0}^y - g_{w_1}^y$ increments slowly: $g_{w_0}^y (i+1) - g_{w_1}^y (i+1) = g_{w_0}^y (i) - g_{w_1}^y (i) + \delta$ where $\delta \in \{-1,0,1\}$. In the baseball interpretation of \cite{Me09Paper}, the difference in score between two teams changes by at most $1$ each innings. When scores are level, we split the words $w_0, w_1$; once on each interval of the game in which scores remain level.

Thus we split $w_0, w_1$ into sub-words such that the difference has a constant sign on each sub-word. Any baseball game splits into sub-games on which one team has the lead. At the end of each sub-game, scores are level; so each $w_0^i, w_1^i$ lie in the same $\M_{n_x^i, n_y^i}$.
\end{Proof}

Having done this, let
\[
w_- = w_0^0 w_1^1 w_0^2 w_1^3 \cdots w_1^{2k-1}, \quad w_+ = w_1^0 w_0^1 w_1^2 w_0^3 \cdots w_0^{2k-1}.
\]
One can easily verify
\[
\begin{array}{ccc}
f_{w_-}^x = \min(f_{w_0}^x, f_{w_1}^x),& g_{w_-}^y = \min( g_{w_0}^y, g_{w_1}^y ),& h_{w_-}^x = \min( h_{w_0}^x, h_{w_1}^x ), \\
f_{w_+}^x = \max(f_{w_0}^x, f_{w_1}^x),& g_{w_+}^y = \max( g_{w_0}^y, g_{w_1}^y ),& h_{w_+}^x = \max( h_{w_0}^x, h_{w_1}^x ),
\end{array}
\]
and similar relations. In particular, $w_- \leq w_+$ and, although there might be many ways to split $w_0, w_1$ (baseball can remain tied for several innings), the resulting $w_-, w_+$ are unique. Thus we define $w_-$ to be the \emph{minimum} and $w_+$ the \emph{maximum} of the pair $(w_0, w_1)$. Note that if $w_0 \leq w_1$ then the minimum of $(w_0, w_1)$ is $w_0$ and the maximum is $w_1$.

\subsection{Adjoints, bilinear form}
\label{sec:adjoints_bilinear_form}

The provocative notation $a_{s,i}$ and $a_{s,i}^*$ for annihilation and creation operators suggests that they should be adjoint. They are indeed adjoint with respect to the partial order $\leq$, but only in one direction each. We have
\[
a_{y,i}^* \; w_0 \leq w_1 \quad \text{iff} \quad w_0 \leq a_{y,i} \; w_1
\]
for all $i$ from $0$ to $n_y (w_1)$, the $y$-degree of $w_1$ (note $w_1$ has $y$-degree one larger than that of $w_0$). We also have
\[
a_{x,i} \; w_0 \leq w_1 \quad \text{iff} \quad w_0 \leq a_{x,i}^* \; w_1
\]
for $0 \leq i \leq n_x (w_0)$ (here $w_0$ has $x$-degree one larger than that of $w_1$).

These inequalities hold even for terminal creations and annihilations: if a terminal annihilation gives $0$ then we count the inequality as false. To see this, note that $a_{y,0}^* w_0$ starts with $y$; if $w_1$ begins with $x$ then $a_{y,0}^* w_0 \nleq w_1$ and $w_0 \nleq a_{y,0} w_1 = 0$. Similar considerations apply for the operators with $(s,i) = (y,n_y + 1), (x,0), (x,n_x+1)$.

(Note these ``adjoint inequalities'' are only true in one direction. It is not true that $a_{y,i} w_0 \leq w_1$ iff $w_0 \leq a_{y,i}^* w_1$: for instance $xxy = a_{y,1} yxxy \leq xyx$ but $yxxy \nleq a_{y,1}^* xyx = xyyx$. Similarly, it is not true that $a_{x,i}^* w_0 \leq w_1$ iff $w_0 \leq a_{x,i} w_1$: for instance $a_{x,1}^* yxy = yxxy \nleq xyyx$ but $yxy \leq a_{x,1} xyyx = yyx$.)

We therefore introduce the notation $\langle \cdot | \cdot \rangle$, which is a boolean version of $\leq$. We can define $\langle \cdot | \cdot \rangle \; : \; \M \times \M \To \{0,1\}$, so that for two words $w_0, w_1 \in \M$, $\langle w_0 | w_1 \rangle = 1$ iff $w_0 \leq w_1$; otherwise $\langle w_0 | w_1 \rangle = 0$. We then extend linearly to a map $\F \otimes \F \To \Z$. 

We can easily verify various properties of $\langle \cdot | \cdot \rangle$.
\begin{itemize}
\item
It is bilinear over $\Z$ (by definition).
\item
It is not symmetric, indeed far from it: if $\langle w_0 | w_1 \rangle = \langle w_1 | w_0 \rangle = 1$ then $w_0 = w_1$.
\item
``Words have norm one'': $\langle w | w \rangle = 1$.
\item
The decomposition $\F = \oplus_{n_x, n_y} \F_{n_x,n_y}$ is orthogonal with respect to $\langle \cdot | \cdot \rangle$: if $w_0, w_1$ have different $x$- or $y$-degree then they are not related by $\leq$, hence $\langle w_0 | w_1 \rangle = 0$. 
\item
Creation and annihilation operators $a_{s,i}, a_{s,i}^*$ are partially adjoint as above,
\[
\langle a_{x,i} w_0 | w_1 \rangle = \langle w_0 | a_{x,i}^* w_1 \rangle \quad \text{and} \quad \langle w_0 | a_{y,i} w_1 \rangle = \langle a_{y,i}^* w_0 | w_1 \rangle,
\]
for $0 \leq i \leq n_x(w_0)$, $0 \leq i \leq n_y(w_1)$ respectively.
\item
``Creation operators are isometries'' (but annihilation operators are not), as discussed above,
\[
\langle w_0 | w_1 \rangle = \langle a_{s,i}^* w_0 | a_{s,i}^* w_1 \rangle.
\]
\item
It is ``multiplicative'': if $(a,c) \in \M_{n_x, n_y}$ and $(b,d) \in \M_{n'_x, n'_y}$, then $\langle ab | cd \rangle = \langle a | c \rangle \; \langle b | d \rangle$.
\item
It is nondegenerate: see below.
\end{itemize}
To see that the bilinear form $\langle \cdot | \cdot \rangle$ is nondegenerate, suppose $\langle v | \cdot \rangle = 0$, for some $v \neq 0$, $v = \sum_i a_i w_i$, a sum of words $w_i \in \M$, with coefficients $0 \neq a_i \in \Z$. By orthogonality of each $\F_{n_x,n_y}$ we may assume all $w_i$ have the same $x$- and $y$-degree. Let $w_-$ be the lexicographically first among the $w_i$; then we have $\langle w_- | w_- \rangle = 1$ but for every other $w_i$, $\langle w_i | w_- \rangle = 0$. Thus $0 = \langle v | w_- \rangle = \langle a_- w_- | w_- \rangle = a_- \neq 0$, a contradiction. Similarly if $\langle \cdot | v \rangle = 0$ then by taking $w_+$, the lexicographically last among the $w_i$, we obtain $\langle w_+ | v \rangle \neq 0$, a contradiction. This gives nondegeneracy.

In fact, our bilinear form $\langle \; | \; \rangle$ is essentially unique, in the following sense.
\begin{prop}
Suppose $B: \F \otimes \F \To \Z$ is a bilinear form such that
\begin{enumerate}
\item 
Distinct summands $\F_n^e$ are orthogonal: If $w_0 \in \M_{n_0}^{e_0}$, $w_1 \in \M_{n_1}^{e_1}$, with $(n_0, e_0) \neq (n_1,e_1)$, then $B(w_0, w_1) = 0$.
\item
Creations and annihilations are partially adjoint: for $0 \leq i \leq n_x(w_0)$, $B( a_{x,i} w_0,  w_1 ) = B( w_0,  a_{x,i}^* w_1 )$; and for $0 \leq i \leq n_y(w_1)$, $B( w_0 , a_{y,i} w_1 ) = B( a_{y,i}^* w_0 , w_1 )$. (When a terminal annihilation gives $0$, we have $B(0, \cdot) = 0$ or $B(\cdot, 0) = 0$.)
\item
$B(1,1) = 1$.
\end{enumerate}
Then $B = \langle \cdot | \cdot \rangle$.
\end{prop}

\begin{Proof}
First note that condition (ii) implies that creation operators are isometries: $B( a_{x,i}^* w_0, a_{x,i}^* w_1 ) = B( a_{x,i} a_{x,i}^* w_0, w_1 ) = B(w_0, w_1)$ and $B( a_{y,i}^* w_0, a_{y,i}^* w_1 ) = B( w_0, a_{y,i} a_{y,i}^* w_1 ) = B(w_0, w_1)$. Note that this works for all creation operators, including both terminal creations.

Given any two words $w_0, w_1 \in \M_n^e$, we note that $B(w_0, w_1)$ can be simplified if $w_0$ has a repeated $y$ or if $w_1$ has a repeated $x$, since then $w_0 = a_{y,i}^* w'_0$ or $w_1 = a_{x,i}^* w'_1$; so suppose there are no such repeated symbols. If $w_0, w_1$ begin with the same symbol then we may simplify, since then $w_0 = a_{y,0}^* w'_0$ or $w_1 = a_{x,0}^* w'_1$; so suppose they do not. If $w_0$ begins with $y$ and $w_1$ begins with $x$, then let $w_1 = a_{x,0}^* w'_1$ so that $B(w_0, w_1) = B(w_0, a_{x,0}^* w'_1) = B(a_{x,0} w_0, w'_1) = 0$. Thus we may simplify to the case where $w_0$ begins with $x$ and has no repeated $y$; and $w_1$ begins with $y$ and has no repeated $x$; but $w_0,w_1$ lie in the same $\M_n^e$.

It follows from these conditions that $w_0$ has $e \leq 0$; and that $w_1$ has $e \geq 0$. As $w_0, w_1$ lie in the same $\M_n^e$, we must have $e=0$, so that $w_0$ is of the form $(xy)^m$ and $w_1$ is of the form $(yx)^m$. Moreover, $B$ is determined from its values on such words. We now compute, for $m \geq 1$:
\begin{align*}
B \left( (xy)^m, (yx)^m \right) &= B \left( (xy)^m y, (yx)^m y \right) \\
&= B \left( a_{y,m}^* (xy)^m, (yx)^m y \right)
= B \left( (xy)^m, a_{y,m} (yx)^m y \right)
= B \left( (xy)^m, (yx)^{m-1} xy \right) \\
&= B \left( (xy)^{m-1}, (yx)^{m-1} \right)
\end{align*}
(In the first line we apply a terminal creation.  In the second line we use adjoint relations to reorder some symbols. In the third line each we remove terminal creations.) Applying this repeatedly we have $B( (xy)^m, (yx)^m ) = B(1,1) = 1$, and hence $B = \langle \cdot | \cdot \rangle$.
\end{Proof}

We remark that, since the ``stackability map'' of \cite{Me09Paper}, mod $2$, easily satisfies the hypotheses of this proposition, this gives another proof of the result, proved directly in \cite{Me09Paper}, that the stackability map mod $2$ is the boolean version of $\leq$.

\subsection{Duality}
\label{sec:duality}

By nondegeneracy, the bilinear form $\langle \cdot | \cdot \rangle$ gives a duality map on each $\F_{n_x,n_y}$, which is an isomorphism, at least over the rationals. In fact, as $\langle \cdot | \cdot \rangle$ is asymmetric, there are two such maps, which are isomorphisms over the rationals. (We will shortly see that these are also isomorphisms over $\Z$.)
\[
\begin{array}{cc}
\iota_- \; : \; \F_{n_x, n_y} \stackrel{\cong}{\To} \left( \F_{n_x, n_y} \right)^*, & v \mapsto \langle v | \cdot \rangle \\
\iota_+ \; : \; \F_{n_x, n_y} \stackrel{\cong}{\To} \left( \F_{n_x, n_y} \right)^*, & v \mapsto \langle \cdot | v \rangle 
\end{array}
\]
Composing these two maps in the two possible directions gives two inverse automorphisms of each $\F_n^e$.
\[
\begin{array}{c}
H = \iota_+^{-1} \circ \iota_- \; : \; \F_n^e \To \F_n^e, \\
H^{-1} = \iota_-^{-1} \circ \iota_+ \; : \; \F_n^e \To \F_n^e.
\end{array}
\]
It's clear from the definition that $\langle u | v \rangle = \langle v | H u \rangle = \langle H^{-1} v | u \rangle$. Thus $H$ is ``unitary'' with respect to $\langle \cdot | \cdot \rangle$, $H^* = H^{-1}$; $H$ is an isometry, $\langle u | v \rangle = \langle Hu | Hv \rangle$.

As mentioned in the introduction, if $\langle \cdot | \cdot \rangle$ were symmetric (resp. antisymmetric), then $H=1$ (resp. $-1$). We will show that $H$ is periodic, the period depending only on the degrees $n_x, n_y$. We will prove this theorem in section \ref{sec:periodicity} using sutured TQFT; we know no direct algebraic proof of this result.
\begin{thm}
\label{thm:H_periodic}
On $\F_n^e = \F_{n_x,n_y}$, $H^{n+1} = (-1)^{n_x n_y}$. In particular $H^{2n+2} = 1$. The period of $H$ is $2n+2$, if both $n_x, n_y$ are odd; else the period is $n+1$.
\end{thm}

As a free abelian group, $\F$ has basis $\M$, i.e. the monomials. There is another bilinear form on $\F$, which we denote by a dot $\cdot$, with respect to which this basis is orthonormal. Let $w_0 \cdot w_1 = 1$ if $w_0 = w_1$ and $w_0 \cdot w_1 = 0$ otherwise. This $\cdot$ is clearly symmetric, $\Z$-bilinear, and nondegenerate. Creation and annihilation are no longer adjoint operators. Creation operators are still isometries: $w_0 \cdot w_1 = a_{s,i}^* w_0 \cdot a_{s,i}^* w_1$ for $s \in \{x,y\}$ and $0 \leq i \leq n_s + 1$. Annihilation operators are still not isometries, e.g. $xy \cdot yx = 0$ but $(a_{x,1} xy ) \cdot (a_{x,1} yx) = y \cdot y = 1$. It is still the case that if $w_0 \cdot w_1 = 1$ and $a_{s,i} w_0, a_{s,i} w_1 \neq 0$ then $w_0 \cdot w_1 = a_{s,i} w_0 \cdot a_{s,i} w_1$; but this does not say much since in this case $w_0 = w_1$.

By nondegeneracy of both bilinear forms $\langle \cdot | \cdot \rangle$ and $\cdot$, it follows that there are operators $Q_+, Q_-: \F_n^e \To \F_n^e$, over $\Q$, which are isomorphisms over $\Q$, intertwining the two forms, i.e. such that for all $u,v \in \F_n^e$
\[
u \cdot v = \langle u | Q_+ v \rangle = \langle Q_- u | v \rangle
\]
But since $\cdot$ is symmetric, we have also
\[
u \cdot v = v \cdot u = \langle v | Q_+ u \rangle = \langle Q_- v | u \rangle = \langle H^{-1} Q_+ u | v \rangle = \langle u | H Q_- v \rangle.
\]
It follows that
\[
H = Q_+ Q_-^{-1},
\]
and hence $(Q_- Q_+^{-1})^{n+1} = (-1)^{n_x n_y}$ on $\F_{n_x, n_y}$.

Given a word $w$, it consists of alternating blocks of $x$ and $y$'s. Read $x$ left to right. An $x$ which is the last in its block, other than possibly an $x$ at the end of $w$, is called a \emph{exceptional} $x$. So the exceptional $x$'s are precisely those followed by a $y$. Let $E_w^x$ denote the set of exceptional $x$'s in $w$. For every subset $T \subseteq E_w^x$, let $\psi_T^x : \M_{n_x, n_y} \To \M_{n_x, n_y}$ denote the operation of taking each $x$ in $T$, and its immediately following $y$, and replacing this $xy$ with $yx$, i.e. ``moving the $x$ forwards one position'' or ``moving the $y$ backwards one position''. For example, if $w = x_1 x_2 y_1 y_2 x_3 y_3 y_4 x_4 y_5 x_5$ (we use subscripts to distinguish distinct letters $x$ and $y$), then $E_w^x = \{x_2, x_3, x_4\}$, and if $T = \{x_2, x_4\}$ then $\psi_T^x w = x_1 y_1 x_2 y_2 x_3 y_3 y_4 y_5 x_4 x_5$. Similarly we denote $E_w^y$ the set of exceptional $y$'s in $w$, i.e. those which are followed by an $x$. For every subset $T \subseteq E_w^y$ denote by $\psi_T^y$ the operation of taking each $y$ in $T$, and its following $x$, and replacing this $yx$ with $xy$, ``moving the $y$ forwards one position'' or ``moving the $x$ backwards one position. If our word is a one-dimensional chessboard, $x$'s are pawns, and $y$'s are empty squares, then $\psi^x$ advances a pawn to the right, and $\psi^y$ advances a pawn to the left. 

We then have the following explicit formulas for $Q_\pm$ and their inverses, describing them on the basis $\M_n^e$ for $\F_n^e$ as free abelian group.
\begin{prop}
\label{prop:Q_formulas}
For any word $w \in \M_n^e$,
\begin{align*}
Q_+ w = \sum_{T \subseteq E_w^y} (-1)^{|T|} \psi_T^y w, \quad Q_+^{-1} w = \sum_{w_i \leq w} w_i \\
Q_- w = \sum_{T \subseteq E_w^x} (-1)^{|T|} \psi_T^x w, \quad Q_-^{-1} w = \sum_{w_i \geq w} w_i
\end{align*}
\end{prop}
(So, for example, if $w = xyxy$ then $Q_- w = xyxy - xyyx - yxxy + yxyx$ and $Q_+ w = xyxy - xxyy$.)

Note it follows that both $Q_+, Q_-$ and their inverses have integer coefficients; hence $Q_+, Q_-$ are automorphisms of each $\F_{n_x, n_y}$ with integer coefficients; hence so is $H$. Clearly the inner product $\cdot$ induces isomorphisms between $\F_{n_x, n_y}$ and the dual $(\F_{n_x, n_y})^*$; as the two bilinear forms are related by the automorphisms $Q_\pm$, it follows that $\iota_\pm$ are isomorphisms also. Hence all of the above holds over $\Z$, not just over $\Q$.

\begin{Proof}
First, consider $Q_+^{-1}$. For any words $w_0,w_1 \in \M_n^e \subset \F_n^e$ we have
\[
w_0 \cdot Q_+^{-1} w_1 = \langle w_0 | w_1 \rangle = \left\{ \begin{array}{cl} 1 & w_0 \leq w_1, \\ 0 & \text{otherwise.} \end{array} \right.
\]
As words are all orthogonal with respect to $\cdot$, $Q_+^{-1} w_1$ contains $w_0$ with coefficient $1$ whenever $w_0 \leq w_1$, otherwise with coefficient $0$. Thus $Q_+^{-1} w$ is as claimed. The proof for $Q_-^{-1}$ is similar.

Now consider $Q_-$. We will show that, for all $w, w_0 \in \M_n^e$,
\[
\sum_{T \subseteq E_w^x} (-1)^{|T|} \left\langle \psi_T^x w \; | \; w_0 \right\rangle = 
\left\langle \sum_{T \subseteq E_w^x} (-1)^{|T|} \psi_T^x w \; \Bigg| \; w_0 \right\rangle = w \cdot w_0 = \left\{ \begin{array}{cl} 1 & w_0 = w, \\ 0 & \text{otherwise,} \end{array} \right.
\]
which gives the desired expression for $Q_-$. Thus we consider when $\langle \psi_T^x w | w_0 \rangle$ is $0$ or $1$.

Consider the function $f_w^x: \{1, \ldots, n_x\} \To \{0, \ldots, n_y\}$ defining $w$ i.e. $f_w^x(i)$ is the number of $y$'s to the left of the $i$'th $x$ in $w$. We note that $\psi_T^x$ has the following effect on the function $f_\cdot^x$:
\[
f_{\psi_T^x w}^x (i) = \left\{ \begin{array}{cl} f_w^x (i) & \text{if the $x$ numbered $i$ is not in $T$,} \\
f_w^x (i) + 1 & \text{if the $x$ numbered $i$ is in $T$.} \end{array} \right.
\]
Thus, $\langle \psi_T^x w | w_0 \rangle = 1$ iff $f_w^x (i) \leq f_{w_0}^x (i)$ for all $i$ numbering $x$'s of $w$ not in $T$, and $f_w^x(i) + 1 \leq f_{w_0}^x (i)$ for all $i$ numbering $x$'s of $w$ in $T$.

So, let $\{1, \ldots, n_x \}$, numbering the $x$'s in $w$, be partitioned into three sets $S_0, S_1, S_2$ as follows:
\begin{align*}
S_0 &= \left\{ i \; : \; f_w^x (i) \leq f_{w_0}^x -1 \right\} \\
S_1 &= \left\{ i \; : \; f_w^x (i) = f_{w_0}^x (i) \right\} \\
S_2 &= \left\{ i \; : \; f_w^x (i) \geq f_{w_0}^x + 1 \right\}
\end{align*}
If $S_2$ is nonempty, then some $x$ in $w$ lies to the right of the corresponding $x$ in $w_0$. So $\langle w | w_0 \rangle = 0$, and any $\langle \psi_T^x w | w_0 \rangle = 0$; also $w_0 \cdot w = 0$. Thus, we may assume $S_2 = \emptyset$. 

If $T$ contains an element $i$ of $S_1$, then the $i$'th $x$ in $\psi_T^x w$ lies to the right of the corresponding $x$ in $w_0$, so $\langle \psi_T^x w | w_0 \rangle = 0$. We have $\langle \psi_T^x w | w_0 \rangle = 1$ iff $T \subseteq S_0$ and $S_2 = \emptyset$. Let $T_0 = S_0 \cap E_w^x$. Then
\[
\sum_{T \subseteq E_w^x} (-1)^{|T|} \langle \psi_T^x w | w_0 \rangle = \sum_{T \subseteq T_0} (-1)^{|T|} \langle \psi_T^x w | w_0 \rangle = \sum_{T \subseteq T_0} (-1)^{|T|}.
\]
The first equality follows since for every $T$ containing an element outside $S_0$, $\langle \psi_T^x w | w_0 \rangle = 0$. The second equality follows since whenever $T \subseteq T_0$, we have $T \subseteq S_0$ and hence $\langle \psi_T^x w | w_0 \rangle = 1$. But the final sum is $0$, unless $T_0$ is the empty set. So the expression above is $1$ only when $S_2 = \emptyset$, and $T_0 = E_w^x \cap S_0 = \emptyset$, i.e. $E_w^x \subseteq S_1$. That is, the expression is $1$ iff $w \leq w_0$ and every block-ending $x$ in $w$ lies in the same position as the corresponding $x$ in $w_0$. This implies $w = w_0$, so that the expression for $Q_-$ is as claimed.

The case of $Q_+$ is similar.
\end{Proof}

As a corollary we have an explicit description for $H$. 
\begin{cor}
\label{cor:H_formula}
\[
H w = Q_+ Q_-^{-1} w = \sum_{w_i \geq w} \sum_{T \subseteq E_{w_i}^y} (-1)^{|T|} \psi_T^y w_i
\]
\qed
\end{cor}

We shall discuss $H$ in more detail in sections \ref{sec:duality_explicity}, \ref{sec:duality_recursively} and \ref{sec:periodicity}.

\subsection{Differentials, commutation relations, normal form}
\label{sec:differentials_commutation_normal_form}

In quantum field theory, for a creation operator $a^*$ and annihilation $a$, one usually has $[a,a^*] = 1$ or $\{a,a^*\} = 1$. A term $aa^*$ creates an particle, and then annihilates every one among that new set of particles; while $a^* a$ annihilates each among the original set of particles, then a new particle is created. Thus $aa^*$ and $a^* a$ count an identical set of situations, except that $aa^*$ counts the original state once more. For bosons then $[a,a^*] = aa^* - a^* a = 1$, and for fermions $\{a,a^*\} = aa^* + a^* a = 1$.

Our non-commutative setting is a little more symmetric than the commutative case: our creations create once, and annihilations annihilate once. The standard physical creation and annihilation are asymmetric: creation creates once, but annihilation annihilates everything in its turn. But in our situation, we can take a sum (or alternating sum) over $i$ of annihilation operators $a_{s,i}$, which should behave similarly to the usual annihilation operator. In any case, define operators
\[
  a_s = \sum_{i=1}^{n_s} a_{s,i}, \quad \quad d_s = \sum_{i=1}^{n_s} (-1)^i a_{s,i} \quad \text{on $\F_{n_x, n_y}$}.
\]
In the first case we obtain $[a_s ,a_{s,0}^*] = 1$; in the second case we have $\{d_s, a_{s,0}^*\} = 1$. 

Note that $a_s$ is nothing more than partial differentiation by the symbol $s \in \{x,y\}$, if we regard each word as a monomial in the two non-commuting variables $x,y$, $a_s = \frac{\partial}{\partial s}$. It follows that the two corresponding operators in $x,y$ commute, $a_x a_y = a_y a_x$. Obviously $a_s$ is compatible with the ring structure on $\F$, satisfying a Leibniz rule: $a_s(w_0 w_1) = (a_s w_0) w_1 + w_0 ( a_s w_1 )$.

On the other hand, $d_s$ is similar to an exterior differential; $d_s^2 = 0$; $d_s(w_0 w_1) = (d_s w_0) w_1 + (-1)^k w_0 (d_s w_1)$ where $k$ is the degree of $w_0$ in $s$; and $d_x d_y = d_y d_x$. We may regard $d_x, d_y$ as boundary operators arising from the two simplicial structures on $\F$. (For this we could also take $d_s = \sum_{i=0}^{n_s + 1} (-1)^i a_{s,i}$, the sum including terminal annihilations; this is also a differential, $d_s^2 = 0$ and $[d_x, d_y] = 0$. Our choice fits better with the physical analogy.)

The differential $d_x$ is a linear operator $\F_{n_x,n_y} \To \F_{n_x - 1, n_y}$, and $d_y$ is $\F_{n_x,n_y} \To \F_{n_x, n_y - 1}$. Thus the abelian groups $\F_{n_x, n_y}$ become a double chain complex, and the relation $\{d_s, a_{s,0}^*\} = 1$ says that $a_{s,0}^*$ is a chain homotopy from $1$ to $0$; thus the homology of this complex is trivial. In fact, the relation $d_s a_{s,0}^* + a_{s,0}^* d_s = 1$ explicitly shows that any ``closed'' element $z$, $d_s z = 0$, is a ``boundary'' since $d_s a_{s,0}^* z + a_{s,0}^* d_s z = d_s \left( a_{s,0}^* z \right) = z$. (See also the argument of proposition 7.5 of \cite{Me09Paper}, following \cite{FraSimplicial}.)

The simplicial set relations show us how to commute creation and annihilation operators of the same species $s \in \{x,y\}$. It follows immediately from them that any sequence of $l$ $s$-creation and $k$ $s$-annihilation operators can be written uniquely in the form (annihilating then creating)
\[
a_{s,i_1}^* a_{s,i_2}^* \cdots a_{s,i_k}^* a_{s,j_1} a_{s,j_2} \cdots a_{s,j_l}
\]
where $i_1 > i_2 > \cdots > i_k$ (note strict inequality), and $j_1 \geq j_2 \geq \cdots \geq j_l$. Alternatively, such a sequence of creation and annihilation operators can be written in the form (creating then annihilating)
\[
a_{s,j_1} a_{s,j_2} \cdots a_{s,j_l} a_{s,i_1}^* a_{s,i_2}^* \cdots a_{s,i_k}^*
\]
satisfying the same inequalities.

\subsection{Temperley--Lieb representation}
\label{sec:Temperley-Lieb}

We can now introduce more operators on $\M$, extended linearly to $\F$:
\[
T_{s,i} = a_{s,i} - a_{s,i+1}, \quad T_{s,i}^* = a_{s,i}^* - a_{s,i+1}^*,
\]
for $s \in \{x,y\}$ and $0 \leq i \leq n_s$.

These operators $T_{s,i}$ have various obvious properties. For $1 \leq i \leq n_s - 1$, $T_{s,i} w = 0$ iff $w$ has the $i$'th and $(i+1)$'th $s$ symbols adjacent. And $T_{s,0} w = 0$ iff $w$ begins with $s$; $T_{s,n_s} w = 0$ iff $w$ ends with $s$.
We have
\[
a_{s,i+1} T_{s,i}^* = a_{s,i+1} \left( a_{s,i}^* - a_{s,i+1}^* \right) = 1 - 1 = 0,
\]
and so if we define
\[
U_{s,i} = T_{s,i}^* a_{s,i+1} = \left( a_{s,i}^* - a_{s,i+1}^* \right) a_{s,i+1}
\]
then we have $\left( U_{s,i}^* \right)^2 = 0$. For $1 \leq i \leq n_s - 1$, the effect of $U_{s,i}$ on a word $w$ is to return $0$, if the $i$'th and $(i+2)$'th $s$ symbols are in the same block. Otherwise it gives difference between the two words, one obtained by moving the $(i+1)$'th $s$ back to be adjacent to the $i$'th $s$, the other obtained by moving the $(i+1)$'th $s$ forward to the $(i+2)$'th $s$. This makes combinatorially clear why $\left( U_{s,i} \right)^2 = 0$. It's then also clear that $U_{s,i}$ and $U_{s,j}$ commute when $|i-j| \geq 2$.

We can then consider $U_{s,i}, U_{s,j}$ when $|i-j| = 1$. We can compute (e.g. putting various products in normal form) that
\begin{align*}
U_{s,i} U_{s,i+1} U_{s,i} &= - U_{s,i} \\
U_{s,i+1} U_{s,i} U_{s,i+1} &= - U_{s,i+1}
\end{align*}

Thus the $U_{s,i}$ satisfy relations similar to the \emph{Temperley-Lieb algebra}. This algebra is defined by the relations $U_i^2 = \delta U_i$, $U_i U_{i+1} U_i = U_i$, $U_{i+1} U_i U_{i+1} = U_{i+1}$, and $U_i U_j = U_j U_i$ for $|i-j| \geq 2$ \cite{Temperley-Lieb, Baxter89, Bernstein_Frenkel_Khovanov}; we have $\delta = 0$ and some sign changes, a ``twisted'' representation.

\subsection{A distinguished subset}
\label{sec:distinguished_subset}

We will now define some distinguished elements in $\F$.

First, consider the initial creation operators. $a_{x,0}^*, a_{y,0}^*$. From the vacuum $1$, applying these operators gives precisely all of the elements of $\M$, hence a basis for $\F$. This distinguished basis subset $\M \subset \F$ is preserved under all creation and non-terminal annihilation operators, and of course is preserved under multiplication (in $\M$ and in $\F$). The set $\{0\} \cup \M$ is also preserved under all creation and annihilation operators and closed under multiplication.

Now we consider a larger set of operators, which are in a sense ``positive''. Define the set
\[
\mathcal{C}^1 = \left\{ a_{s,i}^*, a_{s,i}, T_{s,i}^* \right\} \cdot 1,
\]
i.e. the set of elements which can be obtained by applying these operators to the vacuum $1$. That is, $\mathcal{C}_1$ is the orbit of $1$ under the action of the operators $a_{s,i}^*, a_{s,i}, T_{s,i}^*$, over all $s \in \{x,y\}$ and all $0 \leq i \leq n_s + 1$ in each $\F_n^e$. Similarly, define
\[
\begin{array}{c}
\mathcal{C}^2 = \left\{ a_{s,0}^*, a_{s,0}, H \right\} \cdot 1, \\
\mathcal{C}^3 = \left\{ a_{s, n_s + 1}^*, a_{s, n_s + 1}, H \right\} \cdot 1, \\
\end{array}
\]
So $\mathcal{C}^2$ is the orbit of $1$ under the action of the initial creation and annihilation operators and $H$; $\mathcal{C}^3$ is the orbit of $1$ under the action of final creation and annihilation operators and $H$. In fact these are all the same (and hence one can find many other sets of operators which give this set as their orbit).
\begin{prop}
\label{prop:C1_equals_C2}
$\mathcal{C}^1 = \mathcal{C}^2 = \mathcal{C}^3$.
\end{prop}

Thus we may write $\mathcal{C}$ for this set. Let $\mathcal{C}_n = \mathcal{C} \cap \F_n$, and similarly $\mathcal{C}_{n_x, n_y} = \mathcal{C} \cap \F_{n_x, n_y}$, $\mathcal{C}_n^e = \mathcal{C} \cap \F_n^e$.
\begin{thm}\
\label{thm:properties_of_suture_elements}
\begin{enumerate}
\item
The set $\mathcal{C}$ is not closed under addition or multiplication by $\Z$. However $\mathcal{C}$ is closed under negation, multiplication, and creation and annihilation operators, and the operators $T_{s,i}, T_{s,i}^*, U_{s,i}$.
\item
\begin{enumerate}
\item
Each element $v \in \mathcal{C}$ lies in a particular $\F_n^e$ ($v \in \mathcal{C}_n^e$) and is of the form $v = \sum_i a_i w_i$ where $a_i = \pm 1$ and the $w_i$ are distinct words in $\M_n^e$.
\item
Among the words $w_i$ there is a lexicographically first $w_-$ and last $w_+$. Then for all $i$, $w_- \leq w_i \leq w_+$. If a word $w_0$ among the $w_i$ is comparable to all these $w_i$ with respect to $\leq$, then $w_0 \in \{w_-, w_+\}$. If $w_- \neq w_+$ then $\sum a_i = 0$.
\item
For every pair $w_- \leq w_+$ in $\M_n^e$ there are precisely two $v \in \mathcal{C}_n^e$ with $w_-, w_+$ being the lexicographically first and last words occurring in $v$; one is the negative of the other.
\item
The number of pairs $w_- \leq w_+$ in $\M_n^e$ is the Narayana number $N_n^e$, and $\sum_e N_n^e = C_n$, the $n$'th Catalan number. Thus $|\mathcal{C}_n| = 2 C_n$ and $|\mathcal{C}_n^e| = 2 N_n^e$.
\end{enumerate}
\item
The operators $Q_+, Q_-$ do not preserve $\mathcal{C}$; but $Q_+ \mathcal{C} = Q_- \mathcal{C}$, and these have the same cardinality as $\mathcal{C}$ in each grading: $|Q_\pm \mathcal{C}_n^e| = 2 N_n^e$.
\item
For any $v \in \mathcal{C}$, $\langle v | v \rangle = 1$ and $\langle v | Hv \rangle = 1$. For any $v_0, v_1 \in \mathcal{C}$, $\langle v_0 | v_1 \rangle \in \{-1,0,1\}$.
\item
If the sum or difference of two distinct nonzero elements $u,v$ of $\mathcal{C}$ is also a nonzero element of $\mathcal{C}$, then after switching signs of $u,v$ and swapping $u,v$ if necessary, the triple is $u,v,u-v$, and $\langle u | v \rangle = 1$, $\langle v | u \rangle = 0$. There exists an operator $A^*: \F_2^0 \To \F_n^e$, which is a composition of initial creation operators and applications of $H$, such that $A(xy) = u$ and $A(yx) = v$.
\item
If $u,v \in \mathcal{C}$ and $\langle u | v \rangle = 1$, then there exists a sequence $u = v_0, v_1, \ldots, v_m = v$ in $\mathcal{C}$ such that each $v_i - v_{i+1} \in \mathcal{C}$; and for each $i \leq j$, $\langle v_i | v_j \rangle = 1$.
\end{enumerate}
\end{thm}

The proof of this theorem will be given via building up a theory of sutured TQFT and proving isomorphisms between it and the above non-commutative QFT, culminating in section \ref{sec:main_isomorphism}. The Narayana numbers $N_n^e$ here are as defined in \cite{Me09Paper}.

There are some further tenuous physical analogies and speculations. Elements of $\mathcal{C}$ can be regarded as ``pure states''. The operator $H$ can perhaps be considered as a Hamiltonian generating a time evolution, its periodicity some manifestation of a term $e^{i H t}$. The operators $Q_\pm$, which taken together define $H$, can be considered as supersymmetry; $Q_+ \mathcal{C}$ can be regarded as super-partner-states. Is there some interpretation of sutured manifolds with corners having spacelike and timelike boundaries, and a Hamiltonian giving something like a partition function for contact manifolds?

\section{Sutured TQFT}
\label{sec:sutured_TQFT}

\subsection{Sutured surfaces}
\label{sec:sutured_surfaces}

For our purposes, a \emph{sutured 3-manifold} $(M, \Gamma)$ is a 3-manifold with boundary, with $\Gamma \subset \partial M$ a properly embedded oriented $1$-manifold (i.e. a collection of oriented loops on $\partial M$), such that $\partial M \backslash \Gamma = R_+ \cup R_-$, where $R_\pm$ is an oriented surface with boundary $\pm \Gamma$. (In particular, crossing $\Gamma$ along $\partial M$ takes us from $R_\pm$ to $R_\mp$.)

\begin{figure}
\centering
\includegraphics[scale=0.3]{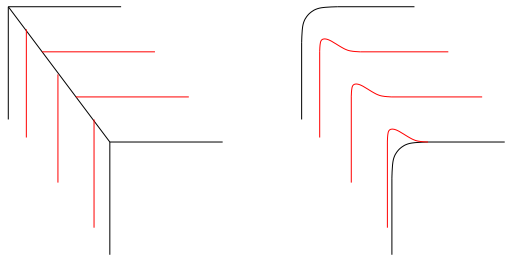}
\caption{Edge smoothing and unsmoothing.} \label{fig:1}
\end{figure}

A sutured 3-manifold may have corners on its boundary; there may be a curve $C$ on $\partial M$ along which two smooth surfaces meet. If so, sutures are required not to match along $C$, but to \emph{interleave}; see figure \ref{fig:1}. The surface may be \emph{smoothed}; in doing so, sutures are rounded and complementary regions $R_\pm$ joined as shown.

In the following, we take sutures on surfaces, but without any associated 3-manifold. Define a \emph{sutured surface} $(\Sigma, \Gamma)$ to be a compact oriented surface $\Sigma$, possibly disconnected, possibly with boundary, with $\Gamma \subset \Sigma$ a properly embedded oriented 1-submanifold; $\Gamma$ must have the property that $\Sigma \backslash \Gamma = R_+ \cup R_-$, where the $R_\pm$ are oriented as $\pm \Sigma$, and where $\partial R_\pm \cap \Gamma = \pm \Gamma$ as oriented 1-manifolds. Again $(\Sigma,\Gamma)$ may have corners, with interleaving sutures; corners may be smoothed as described above. In this paper, we will only consider sutured surfaces with nonempty boundary.

We define also a \emph{sutured background surface} (or simply \emph{sutured background}). Note that the boundary of an oriented 1-manifold may be regarded as a set of points signed $+$ or $-$. A sutured background surface $(\Sigma,F)$ is a compact oriented surface $\Sigma$ (possibly disconnected) with nonempty boundary, together with a finite set of signed points $F \subset \partial \Sigma$, such that $\partial \Sigma \backslash F = C_+ \cup C_-$, where $C_\pm$ are arcs oriented as $\pm \partial \Sigma$, and $\partial C_+ = \partial C_- = F$ as sets of signed points. (Hence each boundary component $C$ of $\Sigma$ has a positive even number of points of $F$, which cut it alternately into arcs of $C_+$ and $C_-$.) A \emph{set of sutures} $\Gamma$ on a sutured background surface $(\Sigma,F)$ is an an oriented properly embedded 1-submanifold of $\Sigma$ such that $\partial \Gamma = \partial \Sigma \cap \Gamma = F$ and such that $(\Sigma, \Gamma)$ is a sutured surface, with $\partial R_\pm = \pm \Gamma \cup C_\pm \cup F$. Again a sutured background surface may have corners.

Given a sutured background surface $(\Sigma, F)$ with $\partial \Sigma \backslash F = C_+ \cup C_-$ as above, we define \emph{gluings} of it. Consider two disjoint 1-manifolds $G_0, G_1 \subseteq \partial \Sigma$, and a homeomorphism $\tau: G_0 \stackrel{\cong}{\To} G_1$ which identifies marked points and positive/negative arcs, $G_0 \cap F \stackrel{\cong}{\to} G_1 \cap F$, $G_0 \cap C_\pm \stackrel{\cong}{\to} G_1 \cap C_\pm$. Then gluing $(\Sigma, F)$ along $\tau$ gives a surface $\#_\tau (\Sigma, F)$. If there remain marked points on each boundary component then $\#_\tau(\Sigma,F)$ is also a sutured background surface and we call $\tau$ a \emph{sutured gluing map}. If $\Gamma$ is a set of sutures on $(\Sigma,F)$ then a sutured gluing map gives a glued set of sutures $\#_\tau \Gamma$ on $\#_\tau (\Sigma,F)$.

A set of sutures $\Gamma$ on $(\Sigma,F)$ has a (relative) \emph{Euler class}, defined by $e(\Gamma) = \chi(R_+) - \chi(R_-)$. Let $|F|=2n$. On a disc, it's clear that, if $\Gamma$ contains no contractible components, then $|e(\Gamma)| \leq n-1$ and $e(\Gamma) \equiv n - 1$ mod $2$. Cutting $\Sigma$ along a properly embedded arc which intersects $\Gamma$ in one interior point has the effect that $e$ is preserved, $\chi(\Sigma)$ increases by $1$ and $n$ increases by $1$; hence $e$ and $n - \chi(\Sigma)$ are preserved. In the other direction, a gluing $\tau$ which identifies two arcs on $\partial \Sigma$, each containing one point of $F$, also preserves $e$ and $n - \chi(\Sigma)$. Hence, on a general $(\Sigma,F)$, if $\Gamma$ contains no contractible components, then $|e(\Gamma)| \leq n - \chi(\Sigma)$ and $e(\Gamma) \equiv \chi(\Sigma)$ mod $2$.

Following \cite{Me09Paper}, we label marked points $F$ on a sutured background surface $(\Sigma,F)$ as follows. Choose a basepoint of $F$ on each component of $\partial \Sigma$. Basepoints are chosen so that, with the orientation on $\partial \Sigma$ induced from the orientation on $\Sigma$, passing through the basepoint one travels from a positive to negative region. 

On a disc sutured background $(D^2, F_{n+1})$ with $|F_{n+1}| = 2(n+1)$, drawn in the plane and inheriting its orientation, we number the points of $F_{n+1}$ as in \cite{Me09Paper}. The basepoint is numbered $0$. The arc of $\partial D^2$ immediately clockwise (resp. anticlockwise) of $0$ is signed positive (resp. negative). The points of $F_{n+1}$ are numbered clockwise, modulo $2(n+1)$. Given a chord diagram $\Gamma$ on $(D^2, F_{n+1})$, from $n$ and $e$ we define $n_- = (n-e)/2$, $n_+ = (n+e)/2$, so $n_-, n_+$ are non-negative integers, $n = n_- + n_+$ and $e = n_+ - n_-$. The point of $F$ numbered $2n_+ + 1 \equiv - 2n_- - 1$ mod $2(n+1)$ is called the \emph{root point}. Cutting $\partial D^2$ at base and root points gives two arcs. The arc containing $1, \ldots, 2n_+$ is called the \emph{eastside}. The arc containing $-1, \ldots, -2n_-$ is called the \emph{westside}.

\subsection{TQFT axioms}
\label{sec:TQFT_axioms}

As in any TQFT, we would like to be associate algebraic objects to our topological objects. Here the topological objects are sutured background surfaces, sets of sutures, and gluings. Thus, we define a \emph{sutured TQFT} to be the following set of assignments with the following properties.
\begin{ax}
\label{ax:1}
To each sutured background surface $(\Sigma, F)$, assign an abelian group $V(\Sigma, F)$, depending only on the homeomorphism type of the pair $(\Sigma,F)$.
\end{ax}

\begin{ax}
\label{ax:2}
To a set of sutures $\Gamma$ on $(\Sigma,F)$, assign a subset of \emph{suture elements} $c(\Gamma) \subset V(\Sigma, F)$, depending only on the isotopy class of $\Gamma$ relative to boundary.
\end{ax}
(We would have liked a single suture element, but this will turn out not to be possible, as we discuss below in section \ref{sec:impossibility_of_coherent_signs}. Quantum states differing by a unit are physically indistinguishable.) 

\begin{ax}
\label{ax:3}
For a sutured gluing map $\tau$ of a sutured background surface $(\Sigma, F)$, assign a collection of linear maps $\Phi_\tau^i \; : \; V(\Sigma, F) \To V(\#_\tau(\Sigma, F))$.
\end{ax}
(We would have liked one canonical map, but this will turn out not to be possible.)

\begin{ax}
\label{ax:4}
For a finite disjoint union of sutured background surfaces $\sqcup_i (\Sigma_i, F_i)$, 
\[
V(\sqcup_i (\Sigma_i, F_i)) = \otimes_i V(\Sigma_i, F_i).
\]
\end{ax}

\begin{ax}
\label{ax:5}
If $\Gamma$ is a set of sutures on $(\Sigma,F)$ and $\tau$ is a gluing of $(\Sigma,F)$ then each $\Phi_\tau^i$ takes suture elements to suture elements surjectively, $c(\Gamma) \to c(\#_\tau \Gamma)$.
\end{ax}
These conditions are not quite the standard TQFT axioms but are similar to the standard axioms for a $(1+1)$-dimensional TQFT. In particular, in a (1+1)-dimensional TQFT we would assign a vector space to a $1$-manifold; instead we assign an abelian group to a decorated 2-manifold with boundary. We would assign an element of a vector space to a $2$-manifold bounded by a $1$-manifold; instead we assign it to a choice of sutures ``bounded by`` the sutured background. We would assign a homomorphism between vector spaces to a 2-dimensional cobordism between $1$-manifolds; instead we assign it to a sutured gluing of $2$-manifolds.

It follows from the above that any $V(\Sigma,F)$ can be interpreted as a space of operators. Take some components of $\partial \Sigma$ and call them \emph{incoming}; call the rest \emph{outgoing}. Write $\partial \Sigma = (\partial \Sigma)_{in} \cup (\partial \Sigma)_{out}$ and $F = F_{in} \cup F_{out}$. Suppose we have a sutured background surface $(\Sigma_{in}, F_{in})$ which has boundary $(\partial \Sigma_{in})$ identified with $(\partial \Sigma)_{in}$ by a gluing $\tau$; and $\tau$ identifies $F_{in}$ and positive/negative boundary arcs on both surfaces. The gluing gives a sutured background surface $(\Sigma_{out}, F_{out})$, where $\Sigma_{out} = \Sigma_{in} \cup \Sigma$, and a (possibly not unique) map
\[
  \Phi_\tau \; : \; V(\Sigma_{in}, F_{in}) \otimes V(\Sigma,F) \To V(\Sigma_{out}, F_{out}),
\]
which is natural with respect to suture elements: if $\Gamma, \Gamma_{in}$ are respectively sets of sutures on $(\Sigma,F)$ and $(\Sigma_{in}, F_{in})$, then $\Gamma_{out} = \Gamma_{in} \cup \Gamma$ is a set of sutures on $(\Sigma_{out}, F_{out})$, and $\Phi_\tau$ maps $c(\Gamma_{in} \cup \Gamma) \to c(\Gamma_{out})$ respectively.

Thus, a choice of sutures $\Gamma$ on $(\Sigma,F)$, together with a choice of representative $c \in c(\Gamma) \subset V(\Sigma, F)$ gives a specific $V(\Sigma_{in}, F_{in}) \To V(\Sigma_{out}, F_{out})$ which we denote $\Phi_{\tau, c}$
\[
  \Phi_{\tau, c} \; : \; V(\Sigma_{in}, F_{in}) \To V(\Sigma_{out}, F_{out}), \quad x \mapsto \Phi_\tau ( x \otimes c ).
\]
(Note the choices involved here: choice of $\Gamma$; choice of $\Phi_\tau$; choice of $c \in c(\Gamma)$.) In a tenuous sense, sutured TQFT can regarded as something like a functor from a cobordism category to an algebraic category, even though we really have no cobordisms, and we think of ``incoming'' and ``outgoing'' in the above tenuous way. With this formulation, we are viewing the situation as an inclusion of background surfaces $(\Sigma_{in}, F_{in}) \hookrightarrow (\Sigma_{out}, F_{out})$, together with a set of sutures $\Gamma$ on $(\Sigma_{out} \backslash \Sigma_{in}, F_{in} \cup F_{out})$. Note that the inclusion must be strict in the sense that $\Sigma_{in}$ lies in the interior of $\Sigma_{out}$. In this way we can regard our ``TQFT'' is a ``$2+1 = 2$ dimensional TQFT''.

In fact, we can reformulate gluing axioms \ref{ax:3} and \ref{ax:5} in terms of inclusions. This is the way maps are described in \cite{HKM08}. 

\begin{ax0}[3']
To an inclusion $(\Sigma_{in}, F_{in}) \stackrel{\iota}{\hookrightarrow} (\Sigma_{out},F_{out})$ of sutured background surfaces, with $\Sigma_{in}$ lying in the interior of $\Sigma_{out}$, together with $\Gamma$ a set of sutures on $( \Sigma_{out} \backslash \Sigma_{in}, F_{in} \cup F_{out})$, assign a collection of linear maps $\Phi_{\iota,\Gamma}^i \; : \; V(\Sigma_{in}, F_{in}) \To V(\Sigma_{out}, F_{out})$.
\end{ax0}
\begin{ax0}[5']
If $\Gamma_{in}$ is a set of sutures on $(\Sigma_{in}, F_{in})$, let $\Gamma_{out} = \Gamma_{in} \cup \Gamma$ be the corresponding set of sutures on $(\Sigma_{out}, F_{out})$. Then each $\Phi_{\iota,\Gamma}^i$ maps takes suture elements to suture elements surjectively, $c(\Gamma_{in}) \to c(\Gamma_{out})$.
\end{ax0}

\begin{lem}
Axioms 3 and 5 are equivalent to axioms 3' and 5'.
\end{lem}
\begin{Proof}
From axioms 3 and 5, for any gluing $\tau$ we have gluing maps $\Phi_\tau^i$. As discussed above, given an inclusion $(\Sigma_{in}, F_{in}) \stackrel{\iota}{\hookrightarrow} (\Sigma_{out}, F_{out})$ and sutures $\Gamma$ on $(\Sigma,F) = (\Sigma_{out} \backslash \Sigma_{in}, F_{in} \cup F_{out})$, consider the gluing $\tau$ of $(\Sigma_{in}, F_{in}) \sqcup (\Sigma, F)$ achieved by the inclusion, the gluing maps $\Phi_\tau^i : V(\Sigma_{in}, F_{in}) \otimes V(\Sigma,F) \To V(\Sigma_{out}, F_{out})$, and take a representative $c \in c(\Gamma) \subset V(\Sigma,F)$; then define
\[
\Phi_{\iota,\Gamma}^i = \Phi_{\tau,c}^i \; : \; V(\Sigma_{in}, F_{in}) \To V(\Sigma_{out}, F_{out}), \quad x \mapsto \Phi_\tau^i(x \otimes c).
\]
In the other direction, assume axioms 3' and 5', so an inclusion $\iota$ and sutures $\Gamma$ as above give inclusion maps $\Phi_{\iota,\Gamma}^i$. Now given a gluing $\tau$ of $(\Sigma,F)$, removing a neighbourhood of the boundary of $(\Sigma,F)$, after gluing, gives an inclusion $(\Sigma,F) \stackrel{\iota}{\hookrightarrow} \#_\tau(\Sigma,F)$, and lying in its interior. Moreover, there is a natural set of sutures $\Gamma$ on a neighbourhood of the boundary of $(\Sigma,F)$, taking a product neighbourhood of the boundary. After choosing $c \in c(\Gamma)$, we obtain inclusion maps $\Phi_{\iota,\Gamma}^i: V(\Sigma,F) \To V(\#_\tau(\Sigma,F))$, which we define to be $\Phi_\tau^i$. (See figure 13 of \cite{HKM08}.)
\end{Proof}

We now impose another condition, rather similar to what occurs in knot theory invariants; alternatively, saying that overtwisted contact elements are zero.
\begin{ax}
\label{ax:6}
If $\Gamma$ contains a closed contractible loop then $c(\Gamma) = \{0\}$.
\end{ax}
(Perhaps more generally we could set $c(\Gamma)$ to be some power of an indeterminate $\delta$ in this case; this would be closer to the sort of situation arising from skein relations in knot theory; perhaps this can be regarded as $\delta = 0$ in a semiclassical limit. But the analogy from contact geometry and sutured Floer homology suggests that one closed loop is equivalent to many.)

For a disc $D^2$, let $(D^2, F_n)$ denote the sutured background surface which is a disc with $2n$ points chosen on the boundary. One of these will be fixed to be a basepoint. Consider $(D^2, F_1)$. There is only one set of sutures $\Gamma$ on this $(D^2,F_1)$ with no contractible loops, namely a single arc. We will call this $\Gamma$ the \emph{vacuum} $\Gamma_\emptyset$, and $(D^2,F_1)$ the \emph{vacuum background}. We will impose the condition, usual in quantum field theory, that the vacuum is nonzero, and normalise the vacuum and its background.
\begin{ax}
\label{ax:7}
$V(D^2, F_1) = \Z$ and $c(\Gamma_\emptyset) \subseteq \{-1,1\}$.
\end{ax}
In more generality, we could set $V(D^2, F_1)$ to be a ground ring $R$ and then say that sutured TQFT is \emph{over $R$ coefficients}. In \cite{Me09Paper} we essentially considered the simpler sutured TQFT with $\Z_2$ coefficients; here we focus on $\Z$ coefficients, though we will need to make reference to $\Z_2$ coefficients and previous work at times.

Note that for any sutured background surface $(\Sigma,F)$ and sutures $\Gamma$ on it, there is an inclusion $(D^2, F_1) \hookrightarrow (\Sigma,F)$, and sutures $\Gamma'$ on $(\Sigma \backslash D^2, F \cup F_1)$, taking $\Gamma_\emptyset \mapsto \Gamma_\emptyset \cup \Gamma' = \Gamma$. Then axiom 3' gives a map $\Z \cong V(D^2, F_1) \To V(\Sigma,F)$ taking $c(\Gamma_\emptyset) \subseteq \{-1,1\} \to c(\Gamma)$ surjectively and we immediately have the following.
\begin{lem}
For any sutures  $\Gamma$ on any sutured background $(\Sigma,F)$, $c(\Gamma)$ is either a singleton or is of the form $\{c,-c\}$.
\qed
\end{lem}
Over $\Z_2$ coefficients rather than $\Z$, these two possibilities are the same, and so $c(\Gamma)$ is always a well-defined single element.

A set of sutures $\Gamma$ on $(D^2, F_n)$ with no closed curves is just a chord diagram. Any disc with a chord diagram may be included into a larger sutured disc in which the sutures simplify to the vacuum. The corresponding inclusion, after making appropriate choices, gives a map $V(D^2, F_n) \To V(D^2, F_1) = \Z$ which takes $c(\Gamma) \to c(\Gamma_\emptyset) \subseteq \{-1,1\}$ surjectively.
\begin{lem}
\label{lem:chord_diag_primitive_non-torsion}
If $\Gamma$ is a chord diagram then $0 \notin c(\Gamma)$ and every element of $c(\Gamma)$ is primitive and non-torsion.
\qed
\end{lem}

The ``interesting'' elements in each $V(\Sigma,F)$ are the suture elements. We will require our sutured TQFT to be ``minimal'' in the following sense. One can easily imagine relaxing or abandoning this axiom; we impose it largely for convenience, and to avoid technical details which, at least for now, we would rather not consider.
\begin{ax}
\label{ax:8}
Every $V(\Sigma,F)$ is spanned by suture elements.
\end{ax}

There clearly exists a non-trivial theory with the axioms so far. For example, take $V(\Sigma, F) = \oplus_\Gamma \Z c_\Gamma$, with one $\Z$ summand for each isotopy class of sutures $\Gamma$ on $(\Sigma,F)$ without contractible loops, and free basis $c_\Gamma$. Given a set of sutures $\Gamma$ without contractible loops, take $c(\Gamma) = \{c_\Gamma\}$; if $\Gamma$ has a contractible loop then $c(\Gamma) = \{0\}$. Clearly axioms 1, 2, 4, 6, 7 and 8 are satisfied. From a sutured gluing map $\tau$ on $(\Sigma,F)$ we get a natural map on sets of sutures, $\Gamma \mapsto \#_\tau \Gamma$; take $\Phi_\tau$ to map $c_\Gamma \mapsto c_{\#_\tau \Gamma}$, if $\#_\tau \Gamma$ has no contractible components, else $c_\Gamma \mapsto 0$. As the $c_\Gamma$ form a free basis of $V(\Sigma,F)$, this defines $\Phi_\tau$ uniquely. This clearly satisfies axioms 3 and 5. We will consider this ``free pre-sutured TQFT'' later in section \ref{sec:variations_of_nondegeneracy}.

We will impose a final axiom which essentially gives all the structure in the TQFT. We will define the bilinear form $\langle \cdot | \cdot \rangle \; : V(D^2, F_n) \otimes V(D^2, F_n) \To \Z$ as mentioned in the introduction, the ``non-commutative inner product'' arising from stacking sutured discs. Our axiom will be a kind of nondegeneracy up to sign.

Let $(\Sigma, \Gamma)$ be the sutured surface which is a cylinder $S^1 \times [0,1]$, with $\Gamma$ consisting of $2n$ parallel arcs of the form $\{\cdot\} \times [0,1]$; and then with a small neighbourhood of a point on one of those sutures removed (a ``leak''). Thus topologically $\Sigma$ is a pair of pants. Regard the two ends of the cylinder as incoming and the boundary of the ``leak'' as outgoing. Consider two incoming sutured background surfaces which are discs $(D^2,F_n)$. When we glue in incoming surfaces we glue them as ends of the original cylinder, so that there are corners; after gluing and rounding, we have the vacuum background $(D^2,F_1)$. Choosing a map $\Phi_\tau$ for gluing in incoming surfaces, and a $c \in c(\Gamma) \subset V(\Sigma, \Gamma)$, we obtain a map
\[
  \langle \cdot | \cdot \rangle \; = \; \Phi_{\tau,c} \; : \; V(D^2, F_n) \otimes V(D^2, F_n) \To V(D^2,F_1) = \Z.
\]
(In particular, there may be several choices for $\langle \cdot | \cdot \rangle$; for now, make an arbitrary choice; we may adjust it later.) Thus $\langle \cdot | \cdot \rangle$ describes ``inclusion into a leaky cylinder''. If we have two sets of sutures $\Gamma_0, \Gamma_1$ on $(D^2, F_n)$ and $c_0 \in c(\Gamma_0)$, $c_1 \in c(\Gamma_1)$, we glue them to $S^1 \times \{0\}$, $S^1 \times \{1\}$ respectively. We then have $\langle c_0 | c_1 \rangle \in \{-1,0,1\}$. To see why, note that after inclusion into the leaky cylinder and rounding corners, the resulting set of sutures on $(D^2, F_1)$ is either the vacuum, and $\langle c_0 | c_1 \rangle = \pm 1$, or there is a contractible suture, and $\langle c_0 | c_1 \rangle = 0$. That is, taking the ``non-leaky'' cylinder with $2n$ parallel sutures (without a small disc removed) and inserting the sutured discs, $\langle c_0 | c_1 \rangle = \pm 1$ (resp. $0$) iff after rounding corners, we have a sutured sphere with connected (resp. disconnected) sutures.

We can now state our axiom precisely.
\begin{ax}
\label{ax:nondegeneracy}
Suppose two elements $x,y \in V(D^2,F_n)$ have the following property: for any set of sutures $\Gamma$ on $(D^2,F_n)$, there exists $c \in c(\Gamma)$ such that $\langle x | c \rangle = \pm \langle y | c \rangle$. Then $x = \pm y$.
\end{ax}
Despite appearances, as mentioned in the introduction, we shall see that this axiom is essentially equivalent to the bypass relation. The formulation of nondegeneracy is somewhat unorthodox; in section \ref{sec:nondegeneracy_axioms} we consider several different alternatives and state their and state their equivalence or non-equivalence; these statements are proved in section \ref{sec:variations_of_nondegeneracy}.

Note that the statement of this axiom is limited to discs; there is no requirement at all for other surfaces. If one tries to define something similar for more complicated surfaces, there is a question of how to define it canonically; and the result of gluing two surfaces into a similar cylinder and removing a disc will no longer have genus $0$, so that there is no longer a clear answer of $\pm 1$ or $0$. In this sense (and maybe this sense alone), sutured TQFT is ``inherently planar''.

We can now give a precise version of theorem \ref{thm:main_thm_rough}. A much more general statement is given in theorem \ref{thm:main_isomorphism}.
\begin{thm}
\label{thm:mini_main_isomorphism}
In any sutured TQFT, there is an isomorphism
\[
\Bigl( V(D^2, F_{n+1}), \langle \cdot | \cdot \rangle \Bigr) \cong \Bigl( \F_n, \langle \cdot | \cdot \rangle \Bigr).
\]
\end{thm}
(Note: this theorem does not assert that any sutured TQFT actually exists!)

\subsection{Nondegeneracy axioms and bypass relations}
\label{sec:nondegeneracy_axioms}

\begin{figure}
\centering
\includegraphics[scale=0.5]{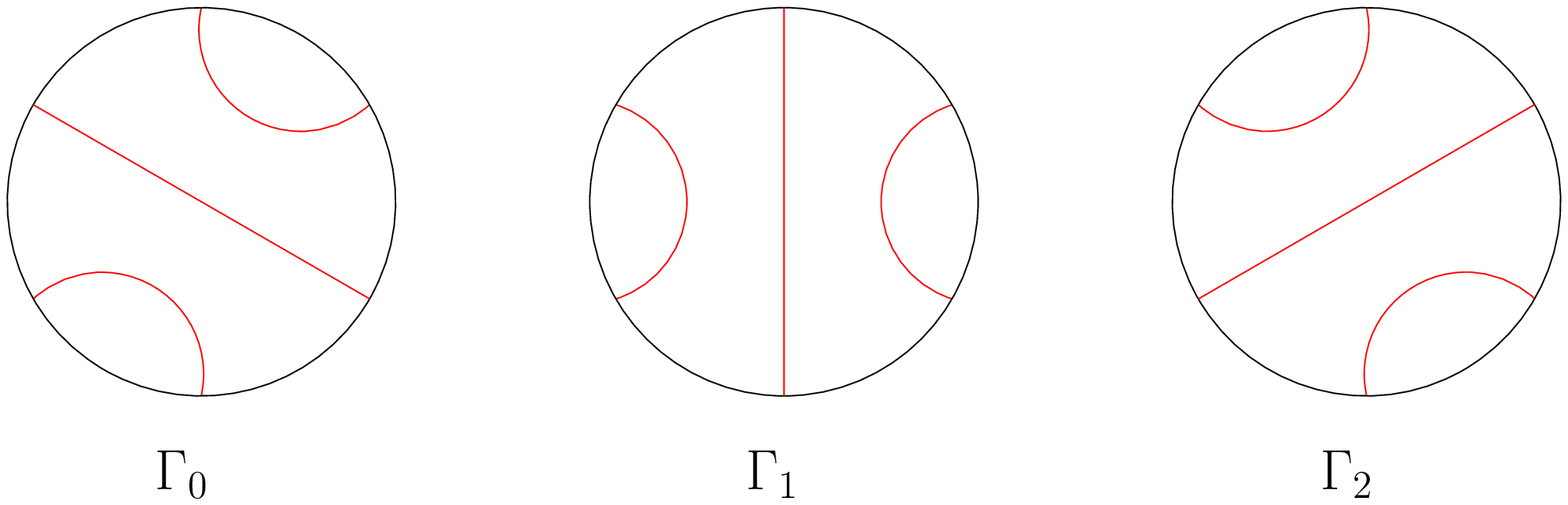}
\caption{Sutures in the bypass relation.} \label{fig:2}
\end{figure}

The nondegeneracy axiom \ref{ax:nondegeneracy} implies certain relations between suture elements. Consider $(D^2, F_3)$, and the three sets of sutures $\Gamma_0, \Gamma_1, \Gamma_2$ shown in figure \ref{fig:2}. Let $c_i \in c(\Gamma_i) \subset V(D^2, F_3)$ be suture elements.

We easily obtain:
\[
\begin{array}{rclrclrcl}
  \langle c_0 | c_0 \rangle &=& \pm 1,	&	\langle c_0 | c_1 \rangle &=& \pm 1,	&	\langle c_0 | c_2 \rangle &=& 0 \\
  \langle c_1 | c_0 \rangle &=& 0,	&	\langle	c_1 | c_1 \rangle &=& \pm 1,	&	\langle c_1 | c_2 \rangle &=& \pm 1 \\
  \langle c_2 | c_0 \rangle &=& \pm 1,	&	\langle c_2 | c_1 \rangle &=& 0,	&	\langle c_2 | c_2 \rangle &=& \pm 1
\end{array}
\]
Take $\alpha = c_0 \pm c_1$, with the plus or minus chosen so that $\langle \alpha | c_1 \rangle = 0$. Then we have
\[
\begin{array}{rclrclrcl}
  \langle \alpha | c_0 \rangle &=& \pm 1, &	\langle \alpha | c_1 \rangle &=& 0,	&	\langle \alpha | c_2 \rangle &=& \pm 1
\end{array}
\]
Thus, by axiom \ref{ax:nondegeneracy}, $\alpha = \pm c_2$. In particular, $c_0, c_1, c_2$ are linearly dependent, and $c_2 = \pm c_0 \pm c_1$ for some choice of sign. 

On a sutured surface $(\Sigma, \Gamma)$, we may consider the operation of removing an embedded disc $D$ in the interior of $\Sigma$ on which the sutures are isotopic (rel boundary of the smaller disc) to a set shown in figure \ref{fig:2} above; and then replacing the sutures on this smaller disc with a different set shown in figure \ref{fig:2}. Such an operation is called \emph{bypass surgery} and comes in two versions: \emph{up}, which takes $\Gamma_0 \mapsto \Gamma_1 \mapsto \Gamma_2 \mapsto \Gamma_0$; and \emph{down}, which takes $\Gamma_0 \mapsto \Gamma_2 \mapsto \Gamma_1 \mapsto \Gamma_0$. Bypass surgery preserves Euler class; bypass-related sutured surfaces naturally come in triples. If $\Gamma'_0, \Gamma'_1, \Gamma'_2$ are a bypass-related triple of sutures on $(\Sigma, F)$, consider including the smaller disc (considered as a sutured background surface) into a larger surface (with fixed sutures) to give $(\Sigma, F)$; choose a gluing map and a suture element for the annulus to obtain a map $V(D^2, F_3) \To V(\Sigma,F)$ taking each $c(\Gamma_i) \to c(\Gamma'_i)$; thus the linear dependency persists.
\begin{lem}
\label{lem:bypass_relation_pm}
Let $\Gamma_0, \Gamma_1, \Gamma_2$ be a bypass related triple of sets of sutures on a sutured background $(\Sigma,F)$. Choose any $c_i \in c(\Gamma_i)$. Then there exists a choice of signs such that $c_0 = \pm c_1 \pm c_2$.
\qed
\end{lem}

If we only consider the groups $V(\Sigma,F)$ mod $2$, then of course $c_0, c_1, c_2$ are single elements and $c_0 + c_1 + c_2 = 0$. In \cite{Me09Paper} we defined groups $SFH_{comb}(T,n)$, which were generated by chord diagrams of $n$ chords, subject to the relation that bypass related triples sum to zero. As such, over $\Z_2$ a sutured TQFT has each $V(D^2, F_n)$ a quotient of $SFH_{comb}(T,n)$. In fact more is true.
\begin{prop}
\label{prop:mod_2_isomorphism}
In any sutured TQFT over $\Z_2$, $V(D^2, F_{n+1}) \cong SFH_{comb}(T,n+1) \cong \Z_2^{2^{n}}$. This isomorphism takes a chord diagram $\Gamma$ in $SFH_{comb}(T,n+1)$ to the suture element $c(\Gamma) \in V(D^2, F_{n+1})$.
\end{prop}
The proof is almost immediate from the definition of sutured TQFT, given our previous work; it is given in section \ref{sec:basis_partial_order}. In sections \ref{sec:basis_partial_order} and \ref{sec:results_mod_2} we recall our previous work and consider sutured TQFT of discs mod 2.

As mentioned above, the formulation of nondegeneracy in axiom \ref{ax:nondegeneracy} is a little unorthodox. Various other formulations are equivalent; along with various formulations of the bypass relation.
\begin{prop}
\label{prop:equivalent_axioms}
In the presence of axioms 1--8 of sutured TQFT, the following axioms are equivalent.
\begin{enumerate}
\item
Suppose two elements $x,y \in V(D^2, F_n)$ have the following property: for any set of sutures $\Gamma$ on $(D^2, F_n)$ and for all $c \in c(\Gamma)$, $\langle x | c \rangle = \pm \langle y | c \rangle$. Then $x = \pm y$.
\item (Original axiom 9.)
Suppose two elements $x,y \in V(D^2, F_n)$ have the following property: for any set of sutures $\Gamma$ on $(D^2, F_n)$, there exists $c \in c(\Gamma)$ such that $\langle x | c \rangle = \pm \langle y | c \rangle$. Then $x = \pm y$.
\item
Suppose two elements $x,y \in V(D^2, F_3)$ have the following property: for any set of sutures $\Gamma$ on $(D^2, F_3)$ and for all $c \in c(\Gamma)$, $\langle x | c \rangle = \pm \langle y | c \rangle$. Then $x = \pm y$.
\item
Suppose two elements $x,y \in V(D^2, F_3)$ have the following property: for any set of sutures $\Gamma$ on $(D^2, F_3)$, there exists $c \in c(\Gamma)$ such that $\langle x | c \rangle = \pm \langle y | c \rangle$. Then $x = \pm y$.
\item
Suppose $\Gamma_0, \Gamma_1, \Gamma_2$ are a bypass triple of sutures on $(D^2, F_n)$. For any $c_0, c_1, c_2$ suture elements in $c(\Gamma_0), c(\Gamma_1), c(\Gamma_2)$ respectively, there exist $\epsilon_1, \epsilon_2 \in \{-1,1\}$ such that $c_0 = \epsilon_1 c_1 + \epsilon_2 c_2$.
\item
Suppose $\Gamma_0, \Gamma_1, \Gamma_2$ are a bypass triple of sutures on $(D^2, F_n)$. Then there exist $c_0, c_1, c_2$ suture elements in $c(\Gamma_0), c(\Gamma_1), c(\Gamma_2)$ respectively, and $\epsilon_1, \epsilon_2 \in \{-1,1\}$ such that $c_0 = \epsilon_1 c_1 + \epsilon_2 c_2$.
\item
Suppose $\Gamma_0, \Gamma_1, \Gamma_2$ are a bypass triple of sutures on $(D^2, F_3)$. For any $c_0, c_1, c_2$ suture elements in $c(\Gamma_0), c(\Gamma_1), c(\Gamma_2)$ respectively, there exist $\epsilon_1, \epsilon_2 \in \{-1,1\}$ such that $c_0 = \epsilon_1 c_1 + \epsilon_2 c_2$.
\item
Suppose $\Gamma_0, \Gamma_1, \Gamma_2$ are a bypass triple of sutures on $(D^2, F_3)$. Then there exist $c_0, c_1, c_2$ suture elements in $c(\Gamma_0), c(\Gamma_1), c(\Gamma_2)$ respectively, and $\epsilon_1, \epsilon_2 \in \{-1,1\}$ such that $c_0 = \epsilon_1 c_1 + \epsilon_2 c_2$.
\end{enumerate}
\end{prop}
Given the foregoing, it's not difficult to see that $(i) \Leftrightarrow (ii) \Rightarrow (iii) \Leftrightarrow (iv) \Rightarrow (v) \Leftrightarrow (vi) \Leftrightarrow (vii) \Leftrightarrow (viii)$. The difficult part is to show [(v)--(viii)] $\Rightarrow$ [(i)--(ii)]; we do this in section \ref{sec:variations_of_nondegeneracy}.

Still, none of the versions of nondegeneracy imposed upon $\langle \cdot | \cdot \rangle$ so far appear orthodox. A more orthodox form of nondegeneracy in fact is not equivalent.
\begin{prop}
\label{prop:inequivalent_axioms}
In the presence of sutured TQFT axioms 1--8, the following two axioms are equivalent.
\begin{enumerate}
\item 
Suppose $x \in V(D^2, F_n)$ has the following property: for any set of sutures $\Gamma$ on $(D^2, F_n)$ and for any suture element $c \in c(\Gamma)$, $\langle x | c \rangle = 0$. Then $x = 0$.
\item
Suppose $x \in V(D^2, F_n)$ has the following property: for any $y \in V(D^2, F_n)$,  $\langle x | y \rangle = 0$. Then $x = 0$.
\end{enumerate}
There exists a sutured TQFT obeying axioms 1--8 and these two alternative axioms 9, such that $V(D^2, F_n) \cong \Z^{C_n}$. Here $C_n$ is the $n$'th Catalan number; denoting the $C_n$ chord diagrams on $(D^2, F_n)$ by $\{\Gamma_i\}_{i=1}^{C_n}$, each $c(\Gamma_i) = \{\pm c_i\}$ and the the $c_i$ form a basis for each $V(D^2, F_n)$.
\end{prop}
The proof is given in section \ref{sec:variations_of_nondegeneracy}. We will construct the example explicitly; there will be \emph{no} linear relation between suture elements for distinct chord diagrams, let alone a bypass relation.

\subsection{Impossibility of coherent signs}
\label{sec:impossibility_of_coherent_signs}

We would like to have each $c(\Gamma)$ a single element. This is possible if we are working modulo $2$, as then signs do not matter. However, as long as $c(\Gamma_\emptyset)$ is non-torsion this is impossible; it remains impossible even if we relax several of the axioms. Everything in this section is a more pedantic version of material appearing in \cite{HKM08}. 
\begin{prop}
\label{prop:impossibility_of_signs}
Consider making assignments:
\begin{enumerate}
\item[1'.]
to each (homeomorphism class of) $(D^2, F_n)$, an abelian group $V(D^2, F_n)$;
\item[2'.]
to each (isotopy class of) set of sutures $\Gamma$ on $(D^2, F_n)$, a suture element $c(\Gamma) \in V(D^2, F_n)$;
\item[3'.]
to an inclusion $(D^2, F_n) \hookrightarrow (D^2, F_m)$ with a set of sutures $\Gamma$ on the intermediate $(S^1 \times I, F_n \cup F_m)$, a collection of linear maps $\Phi^i : V(D^2, F_n) \To V(D^2, F_m)$.
\end{enumerate}
Suppose such assignments satisfy:
\begin{enumerate}
\item[5'.] 
Each $\Phi^i$ is natural with respect to suture elements, i.e. for a set of sutures $\Gamma'$ on $(D^2, F_n)$, $c(\Gamma') \mapsto c(\Gamma' \cup \Gamma)$.
\item[6.]
If $\Gamma$ contains a closed contractible loop then $c(\Gamma) = 0$.
\item[7'.]
$c(\Gamma_\emptyset)$ is torsion-free. 
\item[9'.]
Let $\Gamma_0, \Gamma_1, \Gamma_2$ be the three sets of sutures on $(D^2, F_3)$ in figure \ref{fig:2}. Then the three suture elements $c(\Gamma_0), c(\Gamma_1), c(\Gamma_2) \in V(D^2, F_3)$ are linearly dependent over $\Z$.
\end{enumerate}
Then all suture elements are $0$.
\end{prop}

Here axioms 1', 2', 3' are simply the original axioms 1,2,3, restricted to discs, a particular class of gluings, and suture elements being singletons. Axiom 5' is the original axiom 5 restricted to discs, our particular class of gluings, and singleton suture elements, and axioms 7', 9' are consequences of the original axioms 7, 9 but are weaker. Axioms 4 and 8 are omitted altogether.

\begin{Proof}
Consider the three chord diagrams $\Gamma_0, \Gamma_1, \Gamma_2$, sets of sutures on $(D^2, F_3)$ of figure \ref{fig:2}. Consider the inclusion $(D^2, F_3) \hookrightarrow (D^2, F_2)$ together with the sets of sutures $\Gamma_a, \Gamma_b, \Gamma_c$ on the intermediate annulus as shown in figure \ref{fig:3}. From axiom 3' then we can choose maps $\Phi_a, \Phi_b, \Phi_c: V(D^2, F_3) \To V(D^2, F_2)$. 

\begin{figure}
\centering
\includegraphics[scale=0.3]{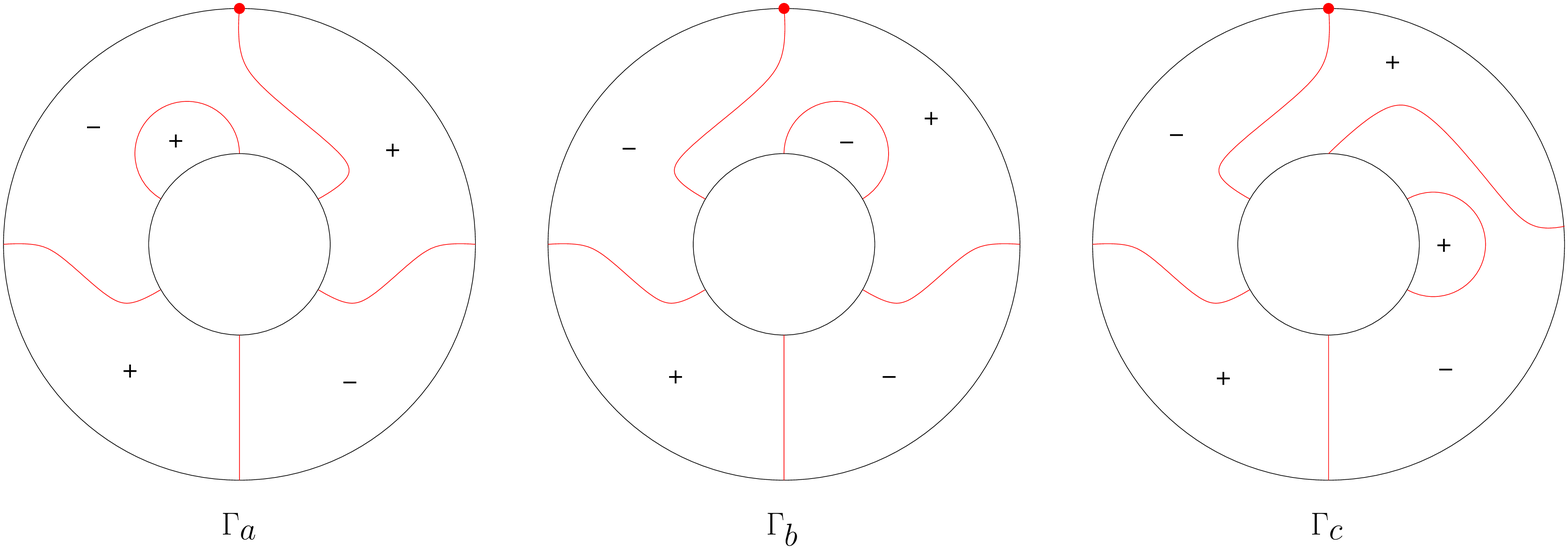}
\caption{Gluing annuli $\Phi_a, \Phi_b, \Phi_c$.} \label{fig:3}
\end{figure}

From the inclusion maps, and contractible loops giving zero, we immediately obtain
\[
\begin{array}{cccc}
\Phi_a \; : \; & c(\Gamma_0) \mapsto c(\Gamma_+)	& c(\Gamma_1) \mapsto c(\Gamma_+) 	& c(\Gamma_2) \mapsto 0\\
\Phi_b \; : \; & c(\Gamma_0) \mapsto 0			& c(\Gamma_1) \mapsto c(\Gamma_-)	& c(\Gamma_2) \mapsto c(\Gamma_-)\\
\Phi_c \; : \; & c(\Gamma_0) \mapsto c(\Gamma_+)	& c(\Gamma_1) \mapsto 0			& c(\Gamma_2) \mapsto c(\Gamma_+)
\end{array}
\]
where $\Gamma_-, \Gamma_+$ are as shown in figure \ref{fig:4}. 

\begin{figure}
\centering
\includegraphics[scale=0.3]{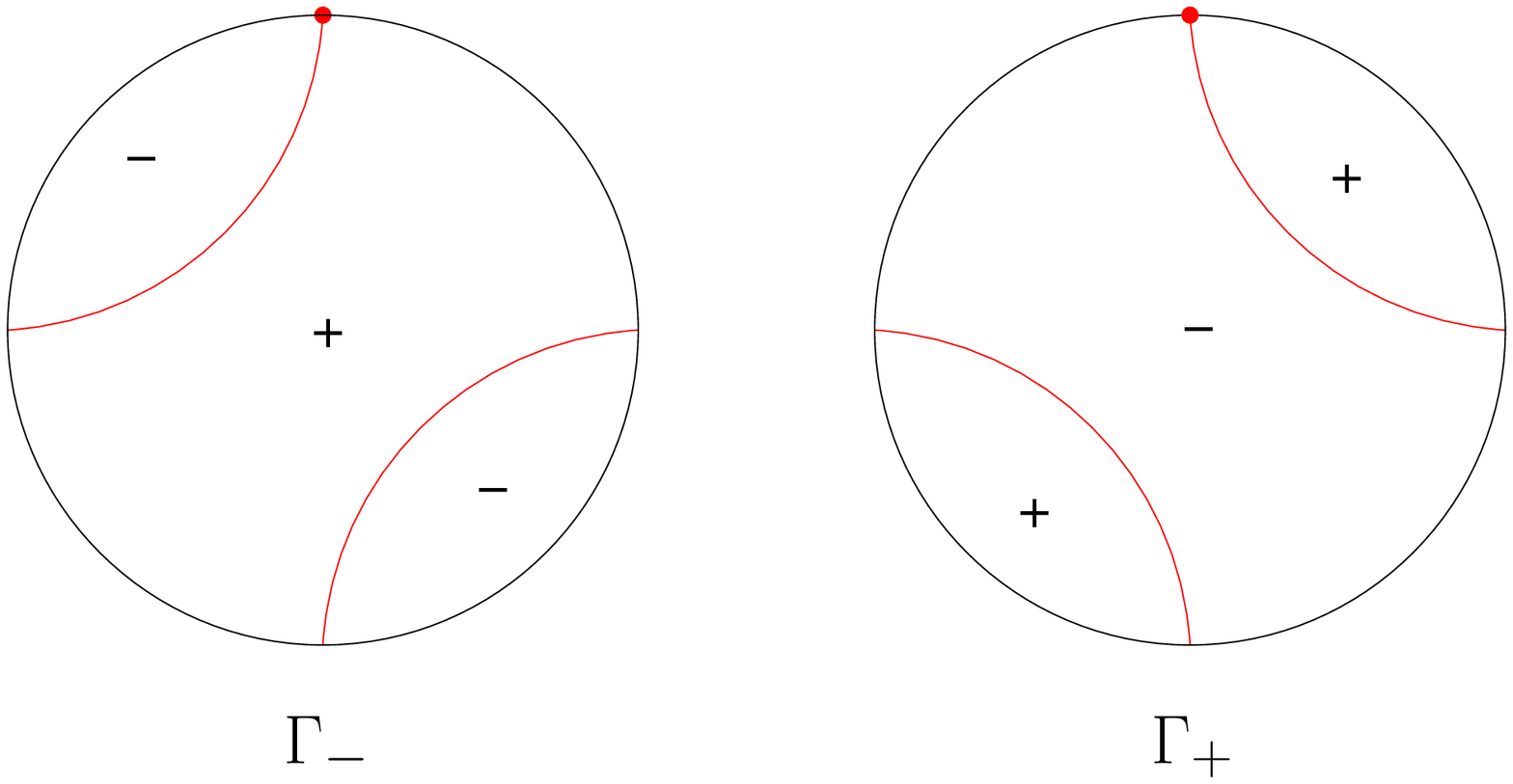}
\caption{Sutures $\Gamma_-, \Gamma_+$.} \label{fig:4}
\end{figure}

As $c(\Gamma_0), c(\Gamma_1), c(\Gamma_2)$ are linearly dependent, let $\alpha c(\Gamma_0) + \beta c(\Gamma_1) + \gamma c(\Gamma_2) = 0$ for some $\alpha, \beta, \gamma \in \Z$, not all zero. From $\Phi_a, \Phi_b, \Phi_c$ we then obtain $(\alpha + \beta) c(\Gamma_+) = 0$, $(\beta + \gamma) c(\Gamma_-) = 0$, and $(\gamma + \alpha) c(\Gamma_+) = 0$ respectively. We may then include $(D^2, F_2) \hookrightarrow (D^2, F_1)$ in various ways so that $\Gamma_+$ or $\Gamma_-$ becomes the vacuum $\Gamma_\emptyset$. Thus $(\alpha + \beta) c(\Gamma_\emptyset) = (\beta + \gamma) c(\Gamma_\emptyset) = (\gamma + \alpha) c(\Gamma_\emptyset) = 0$. As $c(\Gamma_\emptyset)$ is non-torsion, we have either $\alpha + \beta = \beta + \gamma = \gamma + \alpha = 0$ or $c(\Gamma_\emptyset) = 0$. In the first case $\alpha = \beta = \gamma = 0$, a contradiction. In the second case $c(\Gamma_\emptyset)=0$; now an inclusion into any other set of sutures gives that every suture element is $0$.
\end{Proof}

Return now to our original axioms. It follows that, at least on a disc, all chord diagrams have sign ambiguity.
\begin{prop}
\label{prop:sign_ambiguity}
Let $\Gamma$ be a set of sutures on $(D^2, F_n)$. If $\Gamma$ contains a closed loop then $c(\Gamma) = \{0\}$. Otherwise, $\Gamma$ is a chord diagram, and $c(\Gamma)$ has two distinct elements and is of the form $\{x,-x\}$.
\end{prop}

\begin{Proof}
When $\Gamma$ contains a closed loop, it contains a contractible closed loop, so $c(\Gamma) = \{0\}$. Otherwise $\Gamma$ is a chord diagram, and can be included into a larger disc, gluing an annulus to the outside, to obtain the vacuum $\Gamma_\emptyset$. Choosing a gluing map and a suture element for the annulus, we obtain a linear map $V(D^2, F_n) \To V(D^2, F_1)$ which takes $c(\Gamma) \mapsto c(\Gamma_\emptyset)$ surjectively. Conversely, $\Gamma_\emptyset$ includes into a larger disc, gluing an annulus to the outside, which gives a linear map $V(D^2, F_1) \To V(D^2, F_n)$ taking $c(\Gamma_\emptyset) \mapsto c(\Gamma)$ surjectively. The composition of these two maps is a surjection from a finite set to itself; hence a bijection. Thus for any chord diagram $\Gamma$, $|c(\Gamma)| = |c(\Gamma_\emptyset)|$. If $|c(\Gamma_\emptyset)| = 1$ then we have $|c(\Gamma)| = 1$ for all chord diagrams $\Gamma$, contradicting the previous proposition. Thus every $|c(\Gamma)| = 2$; so $c(\Gamma_\emptyset) = \{-1,1\}$, and every $c(\Gamma)$ is of the form $\{x,-x\}$.
\end{Proof}

Note that proposition \ref{prop:sign_ambiguity} does not rely on axiom 9; it relies upon the argument of proposition \ref{prop:impossibility_of_signs}, which in turn relies upon a weaker form of axiom 9. In particular, axiom 9' of proposition \ref{prop:impossibility_of_signs} is implied by any of the formulations in proposition \ref{prop:equivalent_axioms}; to see it is implied by any of axioms (i)--(iv) we use the argument of lemma \ref{lem:bypass_relation_pm}. Thus, we have the following, which we shall need later.
\begin{lem}
\label{lem:strong_sign_ambiguity}
Assume axioms 1--8 of sutured TQFT and any of the 8 alternative formulations of axiom 9 in proposition \ref{prop:equivalent_axioms}. Let $\Gamma$ be any chord diagram. Then $c(\Gamma)$ contains precisely two elements and is of the form $\{x,-x\}$.
\qed
\end{lem}

Consider now the ambiguity in gluing maps; let $\tau$ be a gluing on a sutured background surface $(\Sigma,F)$. Axiom 3 gives at least one map $\Phi_\tau: V(\Sigma,F) \To V(\#_\tau(\Sigma,F))$; axiom 7 requires that $\Phi_\tau$ takes $c(\Gamma) \to c(\#_\tau \Gamma)$ surjectively. By the above, at least on discs, and under the torsion-free assumption, $-\Phi_\tau$ will have the same properties. Thus there is no canonical such $\Phi_\tau$, and any gluing/inclusion necessarily has ambiguity.

\subsection{Creation and annihilation operators}
\label{sec:creation_and_annihilation}

We will use gluing/inclusion to define various operators in sutured TQFT. We generally follow \cite{Me09Paper}, though with some notational differences; see there for further description.

Creation operators are maps $V(D^2, F_n) \To V(D^2, F_{n+1})$; a creation operator includes a disc into a larger disc (or glues a disc to an annulus to give a larger disc) and has the effect of inserting a new outermost chord into a chord diagram in a specific place. Annihilation operators are maps $V(D^2, F_n) \To V(D^2, F_{n-1})$; an annihilation operator also includes a disc into a larger disc, but has the effect of joining two specific adjacent endpoints of a a chord diagram, closing off the region between two marked points $i,i+1$. Formally, each operator is defined by giving a specific sutured annulus. In properly defining each operator, we must choose a gluing map $\Phi_\tau^i$, and a suture element for the sutures on the annulus. This requires careful choices of signs, which we defer to section \ref{sec:choosing_coherent_basis} below; for now we simply consider the effect on sutures.

\begin{figure}
\centering
\includegraphics[scale=0.5]{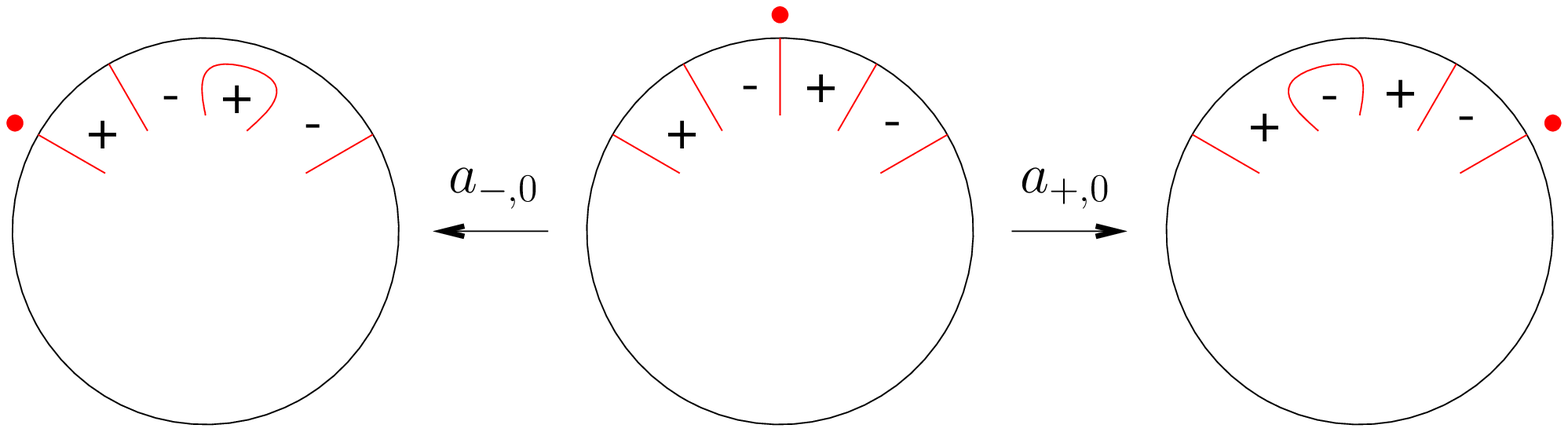}
\caption{Annihilation maps $a_{\pm,0}$.} \label{fig:6}
\end{figure}

Recall from section \ref{sec:sutured_surfaces} our notation for sets of sutures on $(D^2, F_{n+1})$, including numbering of $F_{n+1}$ and the notation $n = n_- + n_+$, $e = n_+ - n_-$, base and root, eastside and westside. We define annihilation operators $a_{\pm,i}: V(D^2, F_{n+1}) \To V(D^2, F_{n})$, for $0 \leq i \leq n_\pm + 1$.
\begin{enumerate}
\item \emph{$(-)$-annihilations}: For $0 \leq i \leq n_- + 1$, $a_{-,i}$ closes off the region between $(-2i, -2i+1)$.
\item \emph{$(+)$-annihilations}: For $0 \leq i \leq n_+ + 1$, $a_{+,i}$ closes off the region between $(2i-1,2i)$
\end{enumerate}
Note that every $(-)$-annihilation closes off a $+$ region. Every $(+)$-annihilation closes off a $-$ region. The \emph{initial} annihilations $a_{-,0}, a_{+,0}$ close off regions at the basepoint. The \emph{final} annihilations $a_{-, n_- + 1}, a_{+, n_+ + 1}$ close off regions at the root point. The \emph{internal} $(-)$-annihilations $a_{-,i}$, $1 \leq i \leq n_-$ close off regions on the westside; the internal $(+)$-annihilations $a_{+,i}$, $1 \leq i \leq n_+$ close off regions on the eastside. (Note this numbering of annihilations is different from \cite{Me09Paper}.) See figure \ref{fig:6}, which also shows how the basepoint behaves.

\begin{figure}
\centering
\includegraphics[scale=0.5]{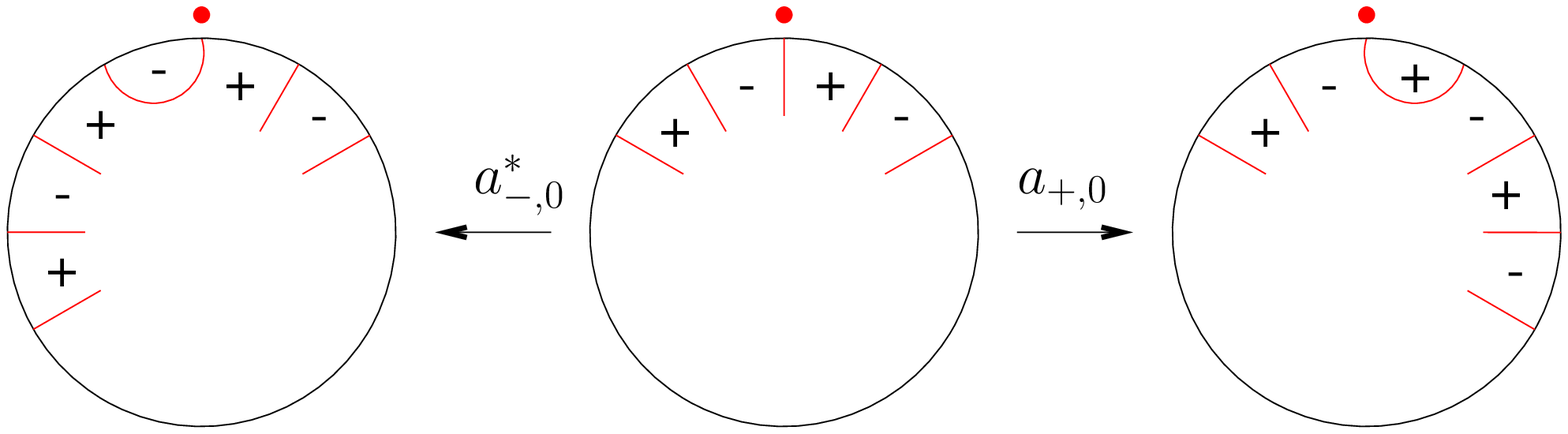}
\caption{Creation maps $a_{\pm,0}^*$.} \label{fig:5}
\end{figure}

Similarly we define creation operators $V(D^2, F_{n+1}) \to V(D^2, F_{n+2})$ as follows.
\begin{enumerate}
\item \emph{$(-)$-creations}: For $0 \leq i \leq n_- + 1$, $a_{-,i}^*$ creates a new chord joining $(-2i-1,-2i)$, between the points previously labelled $(-2i,-2i+1)$.
\item \emph{$(+)$-creations}: For $0 \leq i \leq n_+ + 1$, $a_{+,i}^*$ creates a new chord joining $(2i,2i+1)$, between the points previously labelled $(2i-1,2i)$.
\end{enumerate}
Every $(-)$-creation creates a new chord enclosing an outermost $-$ region; every $(+)$-creation creates a new chord enclosing an outermost $+$ region; \emph{initial} creations create new chords at the basepoint; \emph{final} creations create new chords at the root point; $(-)$-creations create new chords on the westside of the resulting diagram; $(+)$-creations create new chords on the eastside of the resulting diagram. All \emph{internal} $(-)$-creations $a_{-,i}^*$, $1 \leq i \leq n_-$, create new chords in locations on the westside of the original diagram; all internal $(+)$-creations $a_{+,i}^*$, $1 \leq i \leq n_+$, create new chords in locations on the eastside of the original diagram. See figure \ref{fig:5}.

It's easy to check that these creation and annihilation operators obey similar relations to those previously defined in section \ref{sec:algebraic_creation_annihilation}, substituting $(-,+)$ for $(x,y)$, including simplicial relations, as discussed in \cite{Me09Paper}. Any compositions of creations or annihilations which are asserted to be equal in section \ref{sec:algebraic_creation_annihilation} are true for these operators, but (so far!) only up to sign. In sections \ref{sec:choosing_coherent_basis}--\ref{sec:coherent_creation_and_annihilation} we will consider signs carefully and show the same relations hold.

\subsection{Basis, partial order}
\label{sec:basis_partial_order}

Let $\W$ denote the set of words on $\{-,+\}$. Given a word, let $n_\pm$ denote the number of $\pm$ signs, $n$ the total length, and let $e = n_+ - n_-$. Thus $(n,e,n_-.n_+)$ obey the same relations as $(n,e,n_x,n_y)$ in section \ref{sec:algebra}; we write $\W_n$ for words of length $n$, and $\W_n^e = \W_{n_-,n_+}$ for words with $n_\pm$ $\pm$ signs. Note $\W \cong \M$ (the free monoid) as graded monoids, identifying $-$ with $x$ and $+$ with $y$. In particular, the partial order $\leq$ carries naturally over to $\W$.

In \cite{Me09Paper}, we defined a distinguished subset of chord diagrams $\Gamma_w$ on $(D^2, F_{n+1})$, one for each word $w \in \W_n^e$. Starting from the vacuum $\Gamma_\emptyset$, we obtained $\Gamma_w$ by applying the sequence of initial creation operators $a_{\pm,0}^*$ corresponding to $w$; it has Euler class $e = n_+ - n_-$. We showed \cite[section 6.3.2]{Me09Paper} that creation and annihilation operators act on basis diagrams $\Gamma_w$ in perfect analogy to the corresponding creation and annihilation operators on words $w$ in $\{x,y\}$ in $\M_n^e$. That is, defining $a_{\pm,i}$ and $a_{\pm,i}$ to act on words in $\{-,+\}$ in $\W_n^e$, analogously as on words in $\{x,y\}$ in $\M_n^e$, then for a basis diagram $\Gamma_w$, we have $a_{\pm,i}^* \Gamma_w = \Gamma_{a_{\pm,i}^* w}$ and $a_{\pm,i} \Gamma_w = \Gamma_{a_{\pm,i} w}$ (provided $a_{\pm,i}w \neq 0$; if $a_{\pm,i}w = 0$ then $a_{\pm,i} \Gamma_w$ contains a closed loop).

We showed \cite[proposition 1.16]{Me09Paper} that over $\Z_2$ coefficients (where the $c(\Gamma_w)$ are single elements), the $c(\Gamma_w)$ form a basis for $SFH_{comb}(T,n+1)$: ``there is a basis of states given by applying creation operators to the vacuum''. These arguments immediately carry over to sutured TQFT, with either $\Z_2$ or $\Z$ coefficients, showing that representatives of the $c(\Gamma_w)$ form a basis for $V(D^2, F_{n+1})$. Over $\Z_2$, this proves proposition \ref{prop:mod_2_isomorphism} on the isomorphism $V(D^2, F_n) \cong SFH_{comb}(T,n)$ mod 2; we briefly rerun these arguments here.

For every word $w \in \W_n$, form $\Gamma_w$ (which has $n+1$ chords) and choose a representative $c_w \in c(\Gamma_w)$, arbitrarily for now. The argument that the $c_w$ are independent applies verbatim: a suture element of the vacuum is axiomatically nonzero; a suture element of any set of sutures on $V(D^2, F_n)$ containing a closed curve is $0$; if some linear combination of $c_{w_i}$ is zero, $\sum k_i c_{w_i} = 0$, then we can apply annihilation operators to reduce $\Gamma_{w_1}$ to $\Gamma_\emptyset$ but every other $\Gamma_{w_i}$ to a set of sutures containing a closed curve; we obtain $k_1 c(\Gamma_\emptyset) = 0$, a contradiction.

The argument that the $c_w$ span $V(D^2, F_{n+1})$ is based on the observation that a chord diagram $\Gamma$ either has an outermost chord at the basepoint, or is part of a bypass triple in which the other two chord diagrams have outermost chords at the basepoint. In either case, each element of $c(\Gamma)$ is a linear combination of suture elements for chord diagrams which have outermost chords at the basepoint. But if a chord diagram $\Gamma$ has an outermost chord at the basepoint, then $c(\Gamma) = a_{s,0}^* c(\Gamma')$, where $s \in \{-,+\}$, for some chord diagram $\Gamma'$ with fewer chords. Applying this observation repeatedly gives an algorithm (the \emph{base point decomposition algorithm}, \cite[section 4.2]{Me09Paper}) which expresses any suture element of a chord diagram as a linear combination of suture elements, each of which is obtained by applying creation operators to the vacuum, hence as a linear combination of the $c_w$. Axiomatically, suture elements span $V(D^2, F_{n+1})$; so the $c_w$ span as desired. 

We may refine this argument by Euler class; recall $e$ is preserved under bypass surgery. If $\Gamma$ has $n+1$ chords and Euler class $e$, then so do all the $\Gamma_w$ obtained in its decomposition; thus a suture element of $\Gamma$ is a linear combination of $c_w$, where $w \in \W_n^e$. Write $V(D^2)_n^e$ for the submodule of $V(D^2, F_{n+1})$ spanned by suture elements of chord diagrams with Euler class $e$, following the algebraic notation of section \ref{sec:algebra}; also write $V(D^2)_n = V(D^2, F_{n+1})$, $V(D^2)_n^e = V(D^2)_{n_-,n_+}$, and $V(D^2) = \oplus_n V(D^2)_n$. The argument above shows that $c_w$, over $w \in \W_n^e$, form a basis for $V(D^2)_n^e$; so as an abelian group it is free of dimension $\binom{n}{n_-} = \binom{n}{n_+}$; moreover $V(D^2)_n = \oplus_e V(D^2)_n^e$. Thus, as graded abelian groups, $V(D^2) \cong \F$.

\begin{prop}
The $c_w$, for $w \in \W_n^e$, form a basis for $V(D^2)_n^e$.
\qed
\end{prop}

In \cite{Me09Paper} we considered stacking basis chord diagrams $\Gamma_{w_0}, \Gamma_{w_1}$, i.e. including into a cylinder and rounding corners. We showed \cite[proposition 1.28, section 4.3]{Me09Paper} that the set of sutures so obtained on the sphere is connected iff $w_0 \leq w_1$. In particular, if $w_0, w_1$ have distinct Euler class then $w_0 \nleq w_1$ and the sutures are disconnected; so the summands $V(D^2)_n^e$ are orthogonal with respect to $\langle \cdot | \cdot \rangle$. Choosing representatives $c_{w_0} \in c(\Gamma_{w_0})$, $c_{w_1} \in c(\Gamma_{w_1})$, we have $\langle c_{w_0} | c_{w_1} \rangle = \pm 1$ if $w_0 \leq w_1$, and $0$ otherwise. Thus up to sign, the algebraic (in $\F$) and topological (in $V(D^2)$) versions of $\langle \cdot | \cdot \rangle$ agree; so theorem \ref{thm:mini_main_isomorphism} is true up to sign. In the next several sections we consider signs carefully; having done this, we will be able to prove a much more general isomorphism in theorem \ref{thm:main_isomorphism}.

\subsection{Previous results: suture elements mod 2}
\label{sec:results_mod_2}

In \cite{Me09Paper} we effectively considered in detail the structure of chord diagrams and suture elements in $V(D^2, F_n)$ over $\Z_2$ coefficients. All those arguments carry over here, except with sign ambiguities. We briefly recall these arguments and definitions as they are needed here.

We proved several properties of stacking; these give us information about $\langle \cdot | \cdot \rangle$. We have already discussed how, applied to suture elements, this map gives $\pm 1$ or $0$, accordingly as inclusion into a cylinder and rounding gives a sphere with connected or disconnected sutures. Using these facts, we obtain the following.
\begin{itemize}
\item
\cite[lemma 3.7]{Me09Paper} For any chord diagram $\Gamma$, including $\Gamma$ to both ends of the cylinder gives a sphere with connected sutures; so $\langle c(\Gamma) | c(\Gamma) \rangle = \pm 1$.
\item
\cite[lemma 3.9]{Me09Paper} If $\Gamma'$ is obtained from $\Gamma$ by an upwards bypass surgery (so that $\Gamma$ is obtained from $\Gamma'$ by downwards bypass surgery) then $\langle c(\Gamma) | c(\Gamma') \rangle = \pm 1$ and $\langle c(\Gamma') | c(\Gamma) \rangle = 0$.
\item
\cite[lemma 3.8]{Me09Paper} If two chord diagrams $\Gamma, \Gamma'$ share a common outermost chord $\gamma$, then by rounding and re-sharpening corners we have $\langle c(\Gamma) | c(\Gamma') \rangle = \pm \langle c(\Gamma - \gamma) | c(\Gamma' - \gamma) \rangle$.
\item
A similar argument of rounding and re-sharpening corners shows that annihilation and creation operators satisfy the adjoint properties described in section \ref{sec:adjoints_bilinear_form}, up to sign, on any given chord diagrams: $\langle a_{-,i} c(\Gamma_0) | c(\Gamma_1) \rangle = \pm \langle c(\Gamma_0) | a_{-,i}^* c(\Gamma_1) \rangle$ and $\langle c(\Gamma_0) | a_{+,i} c(\Gamma_1) \rangle = \pm \langle a_{+,i}^* c(\Gamma_0) | c(\Gamma_1) \rangle$.
\end{itemize}

We considered words and \cite[section 5.1]{Me09Paper} defined the notion of \emph{elementary move}, \emph{generalised elementary move} and \emph{nicely ordered sequence} of generalised elementary moves. All these moves come in \emph{forwards} and \emph{backwards} versions; forwards/backwards moves move forwards/backwards with respect to $\leq$. A forwards elementary move on a word takes a substring of the form $(-)^a (+)^b$ and replaces it with $(+)^b (-)^a$. A generalised forwards elementary move takes a substring of the form $(-)^{a_1} (+)^{b_1} \cdots (-)^{a_k} (+)^{b_k}$ and replaces it with $(+)^{b_1 + \cdots + b_k} (-)^{a_1 + \cdots + a_k}$; for some $i$ and $j$, it moves the $i$'th $-$ sign (from the left) in $w$ and moves it past the $j$'th (from the left) $+$ sign, together with all $-$ signs in between, and is denoted $FE(i,j)$. A sequence of forwards generalised elementary moves $FE(i_1, j_1), \ldots, FE(i_k, j_k)$ is \emph{nicely ordered} if $i_1 < i_2 < \cdots < i_k$ and $j_1 \leq j_2 \leq \cdots \leq j_k$. The backwards version of all these are obtained by reversing the roles of $-$ and $+$.

We showed \cite[lemmas 5.7--8]{Me09Paper} that two basis chord diagrams $\Gamma_{w_0}, \Gamma_{w_1}$ are bypass related iff $w_0, w_1$ are related by an elementary move; then the third diagram in their bypass triple, mod 2, has suture element $c(\Gamma_{w_0}) + c(\Gamma_{w_1})$; from lemma \ref{lem:bypass_relation_pm}, it has a suture element $\pm c(\Gamma_{w_0}) \pm c(\Gamma_{w_1})$ over $\Z$. We also showed that from a basis chord diagram $\Gamma_{w_0}$, a single bypass surgery either gives a chord diagram of the form $\Gamma_{w_1}$, where $w_0, w_1$ are related by an elementary move; or gives a diagram obtained by performing one of these bypass surgeries in the opposite direction, giving a suture element $c(\Gamma_{w_0}) + c(\Gamma_{w_1})$ mod $2$ or $\pm c(\Gamma_{w_0}) \pm c(\Gamma_{w_1})$ over $\Z$.

For any words $w_0 \leq w_1$, there is a nicely ordered sequence of forwards generalised elementary moves taking $w_0$ to $w_1$. Corresponding to these, there is a set $FBS(w_0,w_1)$ of upwards bypass surgeries on $\Gamma_{w_0}$ which effect the moves, eventually giving $w_1$. We showed, if we performed all these bypass surgeries \emph{downwards} instead, we obtained a chord diagram $\Gamma$ whose suture element, expressed in terms of the basis, was a sum (mod 2) $\sum c(\Gamma_{w_i})$, where $w_0, w_1$ occur in the sum, and for every other $w_i$ occurring in the sum, $w_0 \leq w_i \leq w_1$. In fact the $\Gamma_{w_i}$ in this sum can be obtained by performing upwards bypass surgeries along subsets of $FBS(w_0,w_1)$. Since \cite[proposition 1.19]{Me09Paper} the number of comparable pairs $w_0 \leq w_1$ in $\W_n^e$, and the number of chord diagrams with $n+1$ chords and Euler class $e$, are both the Narayana number $N_n^e$, all suture elements of chord diagrams are of this form \cite[theorem 1.20]{Me09Paper}.

The same argument applies immediately in sutured TQFT with $\Z$ coefficients: for $\Gamma$ a chord diagram and $c \in c(\Gamma)$ a suture element, we have $c = \sum a_i c_{w_i}$, where $c_{w_i} \in c(\Gamma_{w_i})$, the $w_i$ all satisfy $w_0 \leq w_i \leq w_1$. Following the notation of \cite{Me09Paper}, we denote $\Gamma = [\Gamma_{w_0}, \Gamma_{w_1}]$, and write $\pm[w_0,w_1]$ for the corresponding suture elements. All the coefficients $a_i \in \{-1,1\}$, as when we decompose $\Gamma$ in terms of the basis, using the decomposition algorithm, any basis diagram can occur at most once. We also showed \cite[proposition 1.23]{Me09Paper} that if $\Gamma$ is not a basis diagram, then an \emph{even} number of basis diagrams appear in this decomposition; the same diagrams appear over $\Z$ as over $\Z_2$, and so we also see an even number of basis diagrams in the $\Z$ case. There is still a bijection between comparable pairs $w_0 \leq w_1$ and chord diagrams $\Gamma = [\Gamma_{w_0}, \Gamma_{w_1}]$ with minimum and maximum words $w_0, w_1$ occurring in their suture elements.

\subsection{Choosing a coherent basis}
\label{sec:choosing_coherent_basis}

All suture elements suffer irreparably from sign ambiguity. We cannot hope to find canonical representatives for suture elements. But we will choose representatives for our \emph{basis} suture elements, so that we can write other suture elements with respect to them, up to sign. We do this using the creation operators $a_{s,i}^*$ and the bilinear form $\langle \cdot | \cdot \rangle \; : \; V(D^2)_n \otimes V(D^2)_n \To \Z$. As the creation operators and inner product themselves have ambiguities, being gluing maps, we first make some choices to fix signs.

Choose a vacuum representative $c_\emptyset \in c(\Gamma_\emptyset)$ arbitrarily. For each pair $(n_-,n_+)$, choose an arbitrary representative for the creation map $a_{+,0}^*  \; : \; V(D^2)_{n_-,n_+} \To V(D^2)_{n_-,n_+ + 1}$, which is defined up to sign. For each pair $(n_-,n_+) = (n_-, 0)$, choose an arbitrary representative for the creation map $a_{-,0}^* \; : \; V(D^2)_{n_-,0} \To V(D^2)_{n_- + 1,0}$. For each pair $(n_-,n_+)$ there is then a unique sequence of these chosen $a_{-,0}^*$ and $a_{+,0}^*$ operators which lead from $V(D^2)_{0,0}$ to $V(D^2)_{n_-,n_+}$; and this sequence of operators takes the suture element $c(\Gamma_{\emptyset})$ to $c(\Gamma_{w_{max}})$, where $w_{max}$ is the maximum word in $\W_{n_-,n_+}$, i.e. $w_{max} = (+)^{n_+} (-)^{n_-}$. For each $(n_-, n_+)$, we choose a representative $c_{w_{max}} = c_{(+)^{n_+} (-)^{n_-}} \in c(\Gamma_{w_{max}}) \subset V(D^2)_{n_-,n_+}$ to be the image of our chosen vacuum $c_\emptyset$ under this sequence of creations.

Recall that with respect to $\langle \cdot | \cdot \rangle_n : V(D^2)_n \otimes V(D^2)_n \To \Z$, the summands $V(D^2)_n^e$ are orthogonal:
\[
\langle \cdot | \cdot \rangle_n = \bigoplus_e \langle \cdot | \cdot \rangle_n^e
\]
where each $\langle \cdot | \cdot \rangle_n^e$ is a linear map $V(D^2)_n^e \otimes V(D^2)_n^e \To \Z$. We note that we may adjust each $\langle \cdot | \cdot \rangle_n^e$ by a sign, and the inner product map is still a gluing map of the desired type, taking suture elements to suture elements.

In each $\W_n^e$, we have chosen the maximum word $w_{max} = (+)^{n_+} (-)^{n_-}$. We have $\langle c_{w_{max}} | c_{w_{max}} \rangle = \pm 1$. Now adjust each $\langle \cdot | \cdot \rangle_n^e$ by a sign if necessary so that
\[
\langle c_{w_{max}} \; | \; c_{w_{max}} \rangle_n^e = 1.
\]
We now fix once and for all our ``inner product'' $\langle \cdot | \cdot \rangle$ to be the direct sum of these $\langle \cdot | \cdot \rangle_n^e$.

For each $w \in \W_n^e$ and $c_w \in c(\Gamma_w)$, we have $\langle c_w | c_{w_{max}} \rangle = \pm 1$; we choose, once and for all, the representative $c_w$ such that
\[
\langle c_w | c_{w_{max}} \rangle = 1.
\]
We next show this basis is coherent.

\subsection{Coherent creation and annihilation}
\label{sec:coherent_creation_and_annihilation}

We have already chosen signs on some initial creation operators (arbitrarily) in defining a coherent basis. We now choose representatives for all creation and annihilation operators, so that they are coherent with respect to this basis. Recall we defined $a_{\pm,i}$ and $a_{\pm,i}^*$ to act on words in $\W_n^e$, analogously as on words in $\M_n^e$; this corresponds to the action on basis diagrams $\Gamma_w$. So we have $a_{\pm,i}^* c_w = \pm c_{a_{\pm,i}^* w}$ and $a_{\pm,i} c_w = \pm c_{a_{\pm,i} w}$ (where we set $c_0 = 0$).

In $\W_n^e$ we have the maximal word $w_{max} = (+)^{n_+} (-)^{n_-} \in \W_n^e$ and minimal word $w_{min} = (-)^{n_-} (+)^{n_+}$. We choose each creation operator $a_{\pm,i}^*$ to take $c_{w_{max}} \mapsto c_{a_{\pm,i}^* w_{max}}$; and we choose each annihilation operator $a_{\pm,i}$ to take $c_{w_{max}} \mapsto c_{a_{\pm,i} w_{max}}$, unless this is $0$; else we define it to take $c_{w_{min}} \mapsto c_{a_{\pm,i} w_{min}}$. Clearly this agrees with our previous choices of some initial creation operators.

We now show, in several steps, that the creation and annihilation operators, and our choice of basis, are coherent.

\begin{lem}
\label{lem:w0-minus-w1}
If $w_0 \leq w_1$ are words related by an elementary move, then $c_{w_0} - c_{w_1}$ is a suture element.
\end{lem}

\begin{Proof}
There is a suture element of the form $c_{w_0} \pm c_{w_1}$, from the chord diagram obtained by performing downwards bypass surgery on $\Gamma_{w_0}$, where upwards bypass surgery would give $\Gamma_{w_1}$. Then mod $2$ we have $\langle c_{w_0} + c_{w_1} | c_{w_{max}} \rangle = 0$, so rounding corners on the cylinder with the corresponding chord diagrams gives disconnected sutures; hence $\langle c_{w_0} \pm c_{w_1} | c_{w_{max}} \rangle = 0$; but by our choice of basis, $\langle c_{w_0} | c_{w_{max}} \rangle = \langle c_{w_1} | c_{w_{max}} \rangle = 1$; thus the sign is minus, and $c_{w_0} - c_{w_1}$ is a suture element. 
\end{Proof}

\begin{lem}
\label{lem:mini-coherence}
For any word $w \in \W_n^e$, $\langle c_w | c_w \rangle = 1$. For any two words $w_0 \leq w_1$ in $\W_n^e$ related by an elementary move, $\langle c_{w_0} | c_{w_1} \rangle = 1$.
\end{lem}

\begin{Proof}
We use the following fact:
\begin{itemize}
\item
If $w_0 \leq w_1$ are words related by an elementary move, and $\langle c_{w_1} | c_{w_1} \rangle = 1$, then $c_{w_0} - c_{w_1}$ is a suture element and $\langle c_{w_0} | c_{w_0} \rangle = \langle c_{w_0} | c_{w_1} \rangle = 1$.
\end{itemize}
To see why the fact is true, note that from the previous lemma $c_{w_0} - c_{w_1}$ is a suture element. We have, mod $2$, and hence over $\Z$, $\langle c_{w_0} - c_{w_1} | c_{w_1} \rangle = 0$, so that $\langle c_{w_0} | c_{w_1} \rangle = \langle c_{w_1} | c_{w_1} \rangle = 1$. Then $\langle c_{w_0} | c_{w_0} - c_{w_1} \rangle = 0$ so $\langle c_{w_0} | c_{w_0} \rangle = \langle c_{w_0} | c_{w_1} \rangle = 1$.

By our choice of $\langle \cdot | \cdot \rangle$, $\langle c_{w_{max}} | c_{w_{max}} \rangle = 1$. Using the fact repeatedly, we obtain the desired result.
\end{Proof}

We now show that creation and annihilation on $V(D^2)$ behave entirely analogously to creation and annihilation on $\F$. In particular, creation operators are isometries, and creation and annihilation are partial adjoints.

\begin{prop}[Coherence of creation and annihilation]
\label{prop:coherence_creation_annihilation}
For all $w \in \W_n^e$, all $s \in \{-,+\}$ and all $0 \leq i \leq n_s + 1$:
\[
a_{s,i}^* c_w = c_{a_{s,i}^* w} \quad \text{and} \quad a_{s,i} c_w = c_{a_{s,i} w}.
\]
\end{prop}

\begin{Proof}
Consider an annihilation operator $a_{s,i}$. Let $A \subset \W_n^e$ consist of those words $w$ for which $a_{s,i} c_w = c_{a_{s,i} w}$. By definition of $a_{s,i}$, $A$ contains at least one word ($w_{max}$ or $w_{min}$). The result for annihilation operators now follows obviously from repeated application of the following fact.
\begin{itemize}
\item 
If $w \in A$ and $w'$ is obtained from $w$ by an elementary move, then $w' \in A$.
\end{itemize}
To see why it is true, note that from lemma \ref{lem:w0-minus-w1} above, $c_{w'} - c_w$ is a suture element, hence is taken by $a_{s,i}$ to a suture element, possibly $0$. We have $a_{s,i} ( c_{w'} - c_w ) = \pm c_{a_{s,i} w'} - c_{a_{s,i} w}$. But $a_{s,i} w'$ and $a_{s,i} w$ are related by an elementary move, or are identical; hence $c_{a_{s,i} w'} - c_{a_{s,i} w}$ is a suture element (possibly $0$) and the other alternative is not. Thus $w' \in A$.

The argument for creation operators is similar, and simpler, since the $0$ case does not arise.
\end{Proof}

\begin{prop}[Coherence of basis]\
\label{prop:coherent_basis}
\label{prop:suture_elt_structure}
\begin{enumerate}
\item
For words $w_0, w_1 \in \W_n^e$ with $w_0 \leq w_1$, $\langle c_{w_0} | c_{w_1} \rangle = 1$.
\item
For any chord diagram $\Gamma$ and suture element $c \in c(\Gamma)$, decomposing $c$ in terms of the basis,
\[
c = \sum_w a_w c_w \quad \text{where} \quad a_w \in \{-1,1\}.
\]
If $\Gamma$ is not a basis chord diagram then $\sum a_w = 0$.
\item
For each chord diagram $\Gamma$ and suture element $c \in c(\Gamma) \subset V(D^2)_n^e$, $\langle c | c \rangle = 1$.
\end{enumerate}
\end{prop}

\begin{Proof}
We first prove (ii). We already showed in section \ref{sec:results_mod_2} that all $a_w = \pm 1$, considering the decomposition algorithm; we only need show $\sum a_w = 0$ when $\Gamma$ is a non-basis diagram. So if $\Gamma$ is a basis diagram there is nothing to prove. If $\Gamma$ is a non-basis diagram obtained by a single bypass surgery on a basis diagram, then as discussed in \ref{sec:results_mod_2} we have $c = c_{w_0} + c_{w_1}$ mod $2$, where $w_0 \leq w_1$ are related by an elementary move; then lemma \ref{lem:w0-minus-w1} says $c = \pm(c_{w_0} - c_{w_1})$, so coefficients sum to $0$.

Proof by induction on the number of chords in $\Gamma$. All chord diagrams with $1$ chord are basis chord diagrams! If $\Gamma$ has an outermost chord at the basepoint, then $c(\Gamma) = a_{\pm,0}^* c(\Gamma')$ for some $\Gamma'$; since $a_{\pm,0}^*$ acts naturally on each basis element as $a_{\pm,0}^* c_w = c_{a_{\pm,0}^* w}$, applying $a_{\pm,0}^*$ preserves the property that the sum of coefficients is $0$, giving the desired result. Otherwise, we may assume $\Gamma$ has no outermost chord at the basepoint. In this case, apply one stage of the base point algorithm to $\Gamma$; then $c$ is expressed as $\pm c_- \pm c_+$, where each $c_\pm \in c(\Gamma_\pm)$, and $\Gamma_\pm$ has an outermost chord enclosing a $\pm$ region at the base point. As we found when discussing the decomposition algorithm, either $0$ or $2$ of $\Gamma_-, \Gamma_+$ are basis chord diagrams; as we assume $\Gamma$ is not obtained by a single bypass surgery on a basis chord diagram, both $\Gamma_\pm$ are non-basis chord diagrams; and $c_\pm = a_{\pm,0}^* c'_\pm$ for some suture elements $c'_\pm$ for non-basis chord diagrams with fewer chords. By induction, when decomposed, the coefficients of $c'_\pm$ both sum to $0$; hence the coefficients of $c$ also sum to $0$.

Next we prove (iii). We will repeatedly apply this fact:
\begin{itemize}
\item
Let $x \in c(\Gamma_x)$, $y \in c(\Gamma_y)$ be suture elements, where $\Gamma_x, \Gamma_y$ are chord diagrams related by bypass surgery. Suppose $\langle x | x \rangle = \langle y | y \rangle = 1$. Reorder $x,y$ if necessary so that $\langle x | y \rangle = \pm 1$ and $\langle y | x \rangle = 0$ (i.e. so $\Gamma_y$ is obtained from $\Gamma_x$ by an upwards bypass surgery), and replace $y$ with its negative if necessary, so that $\langle x | y \rangle = 1$. Then $x-y$ is a suture element, $\langle x | y \rangle = 1$, and $\langle x - y | x - y \rangle = 1$.
\end{itemize}
To see why, note that there is a suture element of the form $x \pm y$. We have $\langle x \pm y | x \pm y \rangle = 2 \pm \langle x | y \rangle$; hence $2 \pm 1 = \pm 1$; of course we must have $2 - 1 = 1$, and the conclusions follow.

Consider the set $A$ of suture elements $c$ such that $\langle c | c \rangle = 1$. We know from lemma \ref{lem:mini-coherence} above that $A$ contains all $c_w$, and (by the above fact) if it contains two suture elements for chord diagrams related by a bypass surgery, then it contains the suture element for the third chord diagram in their bypass triple. Repeated application of the above fact (e.g. considering a decomposition of a chord diagram into basis diagrams, using the base point decomposition algorithm) gives that $\langle c | c \rangle = 1$ for all suture elements of chord diagrams, proving (iii).

We now prove (i). Take $w_0 \in \W_n^e$; suppose there is some $w_1$ such that $w_0 \leq w_1$ but $\langle c_{w_0} | c_{w_1} \rangle \neq 1$, i.e. $\langle c_{w_0} | c_{w_1} \rangle = -1$. Among such $w_1$ we may take the one which is least with respect to the lexicographic ordering on $\W_n^e$; hence such $w_1$ is also minimal with respect to the partial order $\leq$. We have shown (lemma \ref{lem:mini-coherence}) that $\langle c_{w_0} | c_{w_0} \rangle = 1$; hence $w_0 \neq w_1$. Now take the chord diagram $\Gamma = [\Gamma_{w_0}, \Gamma_{w_1}]$ and take $c \in c(\Gamma)$. As $w_0 \neq w_1$, $\Gamma$ is not a basis diagram. In mod 2 we have $\langle c_{w_0} | c \rangle = 0$; the same is true over $\Z$. Writing $c = \sum_w a_w c_w$, where all $a_w \neq 0$, we now have that each $a_w = \pm 1$ and $\sum_w a_w = 0$. Thus
\[
\sum_w a_w \langle c_{w_0} | c_w \rangle = 0.
\]
Now every $a_w = \pm 1$ and for each $w$ occurring in the sum, $w_0 \leq w \leq w_1$. By minimality of $w_1$ then, we have $\langle c_{w_0} | c_w \rangle = 1$ for all $w$, except possibly $w = w_1$. Thus
\[
\sum_w a_w = 0 \quad \text{and} \quad a_{w_1} \langle c_{w_0} | c_{w_1} \rangle + \sum_{w \neq w_1} a_w = 0.
\]
These two sums are identical, except that $a_{w_1}$ in the first is replaced with $a_{w_1} \langle c_{w_0} | c_{w_1} \rangle$ in the second. Thus $\langle c_{w_0} | c_{w_1} \rangle = 1$, a contradiction, and we are done.
\end{Proof}

\subsection{Multiplication}

Given two chord diagrams $\Gamma_0, \Gamma_1$ on $(D^2, F_{n_0 + 1}), (D^2, F_{n_1 + 1})$ respectively, we now give the operation of \emph{multiplying} them, or rather, their suture elements. Consider a gluing $\tau$ which identifies a neighbourhood of the root point on the boundary of $(D^2, F_{n_0})$, with a neighbourhood of the basepoint on the boundary of $(D^2, F_{n_1})$. Choosing a gluing map $\Phi_\tau$, we have a linear operator
\[
\times_{n_0,n_1} \; : \; V(D^2)_{n_0} \otimes V(D^2)_{n_1} \To V(D^2)_{n_0 + n_1}.
\]
If $\Gamma_0, \Gamma_1$ have Euler classes $e_0, e_1$ then the gluing gives a chord diagram on $(D^2, F_{n_0 + n_1 + 1})$ with Euler class $e_0 + e_1$. Restricting to the $e_0, e_1$ summands, we obtain an operator
\[
\times_{(n_0,e_0),(n_1,e_1)} \; : \; V(D^2)_{n_0}^{e_0} \otimes V(D^2)_{n_1}^{e_1} \To V(D^2)_{n_0 + n_1}^{e_0 + e_1}
\]
where $\times_{n_0,n_1} = \oplus_{e_0,e_1} \times_{(n_0,e_0),(n_1,e_1)}$. Thus this operator respects all gradings. Moreover it is easily checked that it takes basis suture elements to basis suture elements: taking two words $w_0, w_1$, we have
\[
c(\Gamma_{w_0}) \otimes c(\Gamma_{w_1}) \mapsto c(\Gamma_{w_0 w_1}),
\]
where $w_0 w_1$ is the concatenation, or multiplication, of the two words $w_0 w_1$.

This operator can therefore be regarded as a $\Z$-bilinear multiplication, making $V(D^2)$ into a bigraded ring. There are sign ambiguities arising from the ambiguities in suture elements; we now resolve them. As $\times_{n_0, n_1} = \oplus_{e_0, e_1} \times_{(n_0, e_0),(n_1,e_1)}$, we may adjust each separate $\times_{(n_0,e_0),(n_1,e_1)}$ individually. We choose $\times_{(n_0,e_0),(n_1,e_1)}$ to send maximal basis elements, coherently oriented as chosen above, to basis elements, i.e.
\[
c_{(+)^{n_{+,0}} (-)^{n_{-,0}}} \otimes c_{(+)^{n_{+,1}} (-)^{n_{-,1}}} \mapsto c_{(+)^{n_{+,0}} (-)^{n_{-,0}} (+)^{n_{+,1}} (-)^{n_{-,1}}}
\]
where $n_{+, i}, n_{-,i}$ is the pair $n_+, n_-$ corresponding to $(n,e) = (n_i,e_i)$. Then we can show that all multiplication is coherently oriented with respect to our basis. The method of proof is by now familiar.
\begin{prop}
For any words $w_0 \in \W_{n_0}^{e_0}$ and $w_1 \in \W_{n_1}^{e_1}$, under $\times_{(n_0,e_0),(n_1,e_1)}$:
\[
c_{w_0} \otimes c_{w_1} \mapsto c_{w_0 w_1}.
\]
\end{prop}

\begin{Proof}
Consider the set $A$ of words for which this multiplication is coherent: $A = \{(w_0,w_1) \; : \; c_{w_0} \otimes c_{w_1} \mapsto c_{w_0 w_1} \}$. By our sign choice of $\times_{(n_0,e_0),(n_1,e_1)}$, the pair of maximum words
\[
\left( (+)^{n_{+,0}} (-)^{n_{-,0}}, (+)^{n_{+,1}} (-)^{n_{-,1}} \right) \in A.
\]
The result now follows obviously from the following two facts.
\begin{itemize}
\item
If $(w_0, w_1) \in A$ and $w'_0$ is related to $w_0$ by an elementary move, then $(w'_0, w_1) \in A$.
\item
If $(w_0, w_1) \in A$ and $w'_1$ is related to $w_1$ by an elementary move, then $(w_0, w'_1) \in A$.
\end{itemize}
To see why the first fact is true, note that $c_{w'_0} - c_{w_0}$ is a suture element, and multiplication must take $(c_{w'_0} - c_{w_0}) \otimes w_1$ to a suture element. The result of this multiplication must be $\pm c_{w'_0 w_1} - c_{w_0 w_1}$. But the concatenations $w'_0 w_1$ and $w_0 w_1$ are clearly related by an elementary move, hence $c_{w'_0 w_1} - c_{w_0 w_1}$ is a suture element and $-c_{w'_0 w_1} - c_{w_0 w_1}$ is not. So $c_{w'_0} \otimes c_{w_1} \mapsto c_{w'_0 w_1}$ and $(w'_0, w_1) \in A$. The second fact is identical but reversing first and second entries.
\end{Proof}

\subsection{Temperley--Lieb algebra}
\label{sec:Temperley-Lieb_STQFT}

\begin{figure}
\centering
\includegraphics[scale=0.25]{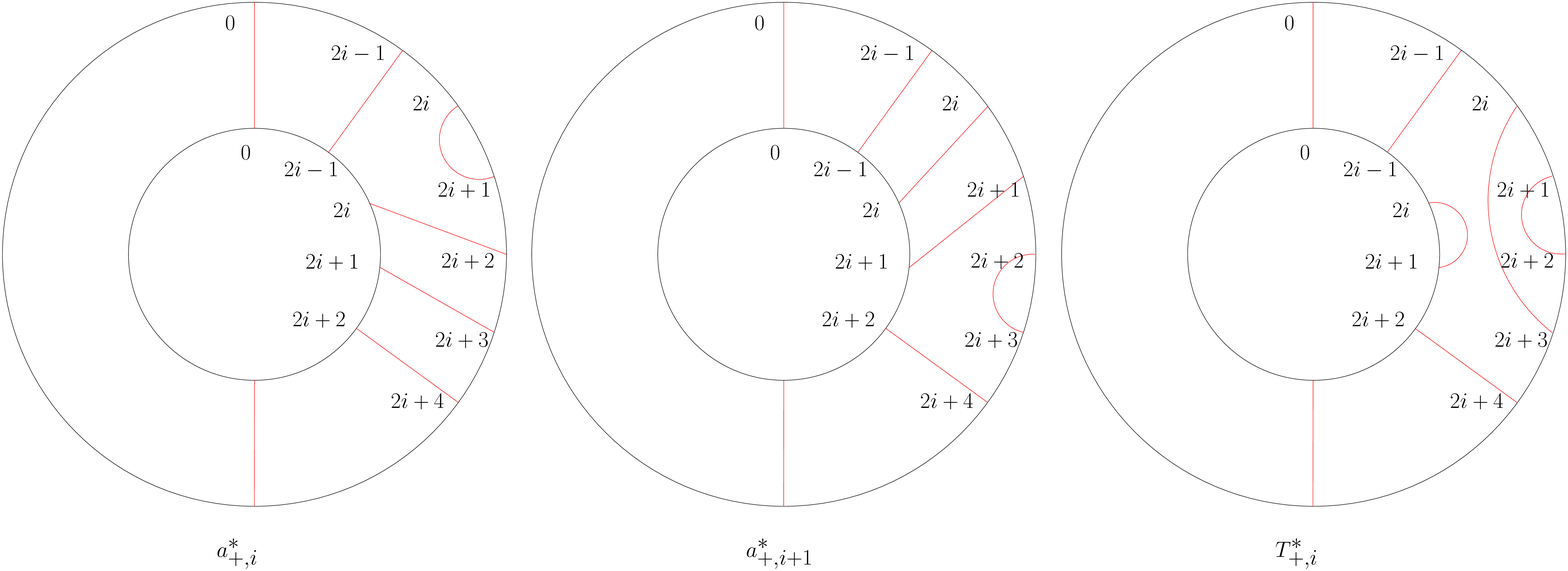}
\caption{Operators $a_{+,i}^*, a_{+,i+1}^*$ and $T_{+,i}^*$.} \label{fig:7}
\end{figure}

Consider the two creation operators $a_{+,i}^*$ and $a_{+,i+1}^*$, acting $V(D^2)_n^e \To V(D^2)_{n+1}^{e+1}$, for $0 \leq i \leq n_+$. Both are obtained by gluing an annulus with specific sutures to the exterior of a disc. We notice that the two sets of sutures are bypass related. Consider a third set of sutures on the annulus, forming a bypass triple. This can also be glued to the exterior of a disc, and we obtain a gluing map
\[
T_{+,i}^* \; : \; V(D^2)_n^e \To V(D^2)_{n+1}^{e+1}.
\]
See figure \ref{fig:7}. Similarly on the westside, we can consider $a_{-,i}^*$ and $a_{-,i+1}^*$, for $0 \leq i \leq n_-$, which give annuli with bypass-related sutures, and obtain $T_{-,i}^*: V(D^2)_n^e \To V(D^2)_{n+1}^{e-1}$.2

Similarly for annihilation operators, $a_{\pm,i}, a_{\pm,i+1}$, for $0 \leq i \leq n_\pm$, give annuli with bypass-related sutures, and taking a sutured annulus forming a bypass triple with them, we obtain a gluing map
\[
T_{\pm,i} \; : \; V(D^2)_n^e \To V(D^2)_{n-1}^{e \pm 1}.
\]

These gluing maps are ambiguous up to sign; we now choose coherent signs. We note that $T_{+,i}^*$ takes each basis element $c_w$ to $\pm c_{a_{+,i}^* w} \pm c_{a_{+,i+1}^* w}$. Note that $a_{+,i}^* w$ and $a_{+,i+1}^* w$ are words related by an elementary move, or are identical; hence $\pm \left( c_{a_{+,i}^* w} - c_{a_{+,i+1}^* w} \right)$ are suture elements, and the other two possibilities are not. Similar considerations apply to $T_{-,i}^*$. We can therefore choose a sign for $T_{\pm,i}^*$ by setting
\[
T_{\pm,i}^* = a_{\pm,i}^* - a_{\pm,i+1}^*
\]
Similarly, $T_{+,i}$ takes $c_w$ to $\pm c_{a_{+,i}w} \pm c_{a_{+,i+1}w}$, and the two words $a_{+,i} w, a_{+,i+1}w$ are either identical or related by an elementary move. Either way it is $\pm \left( c_{a_{+,i} w} - c_{a_{+,i+1} w} \right)$ that are suture elements; the same applies to $T_{-,i}$. We choose a sign for $T_{\pm,i}$:
\[
T_{\pm,i} = a_{\pm,i} - a_{\pm,i+1}.
\]

\begin{figure}
\centering
\includegraphics[scale=0.25]{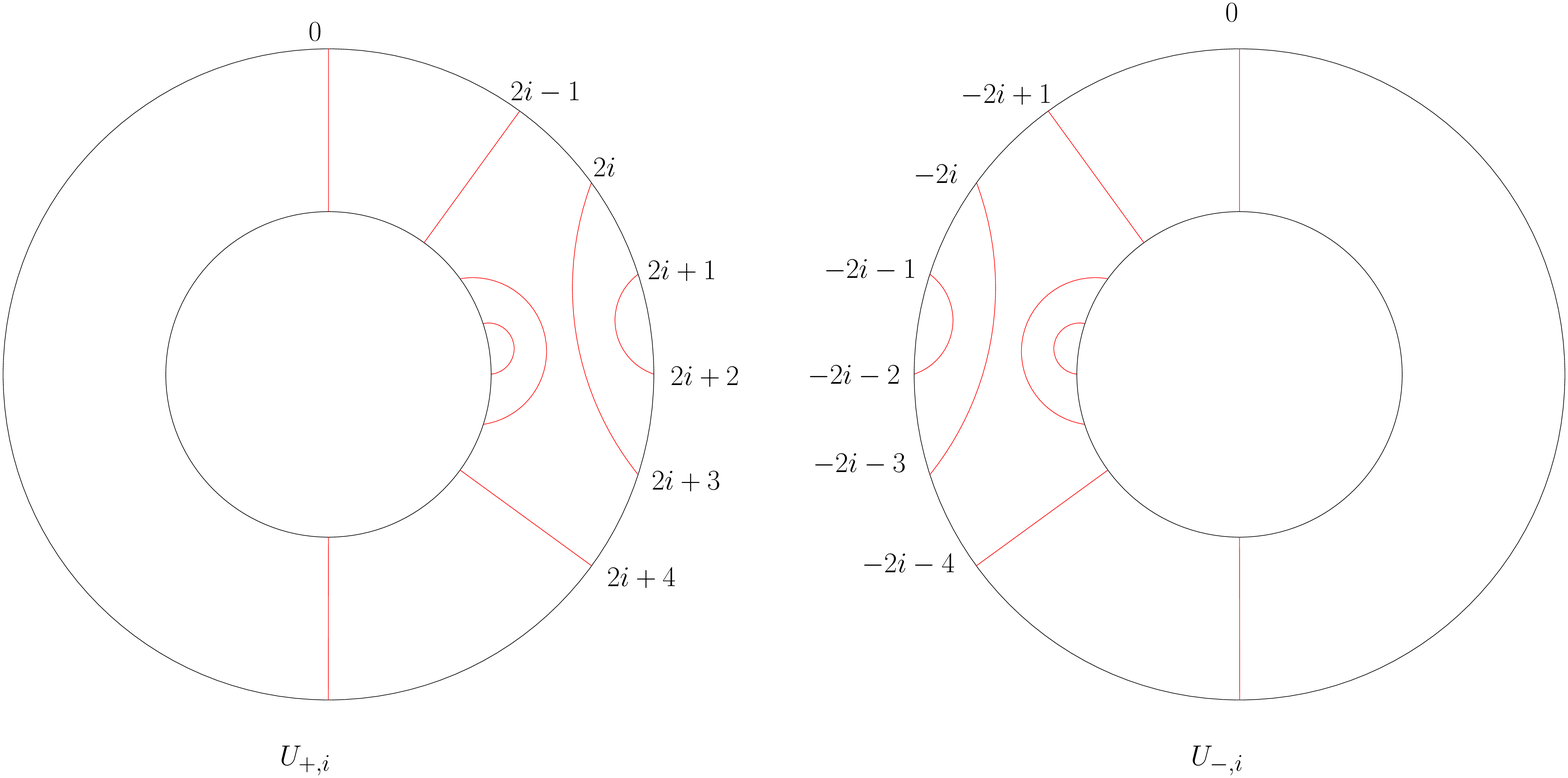}
\caption{Operators $U_{+,i}$, $U_{-,i}$.} \label{fig:8}
\end{figure}

\begin{figure}
\centering
\includegraphics[scale=0.25]{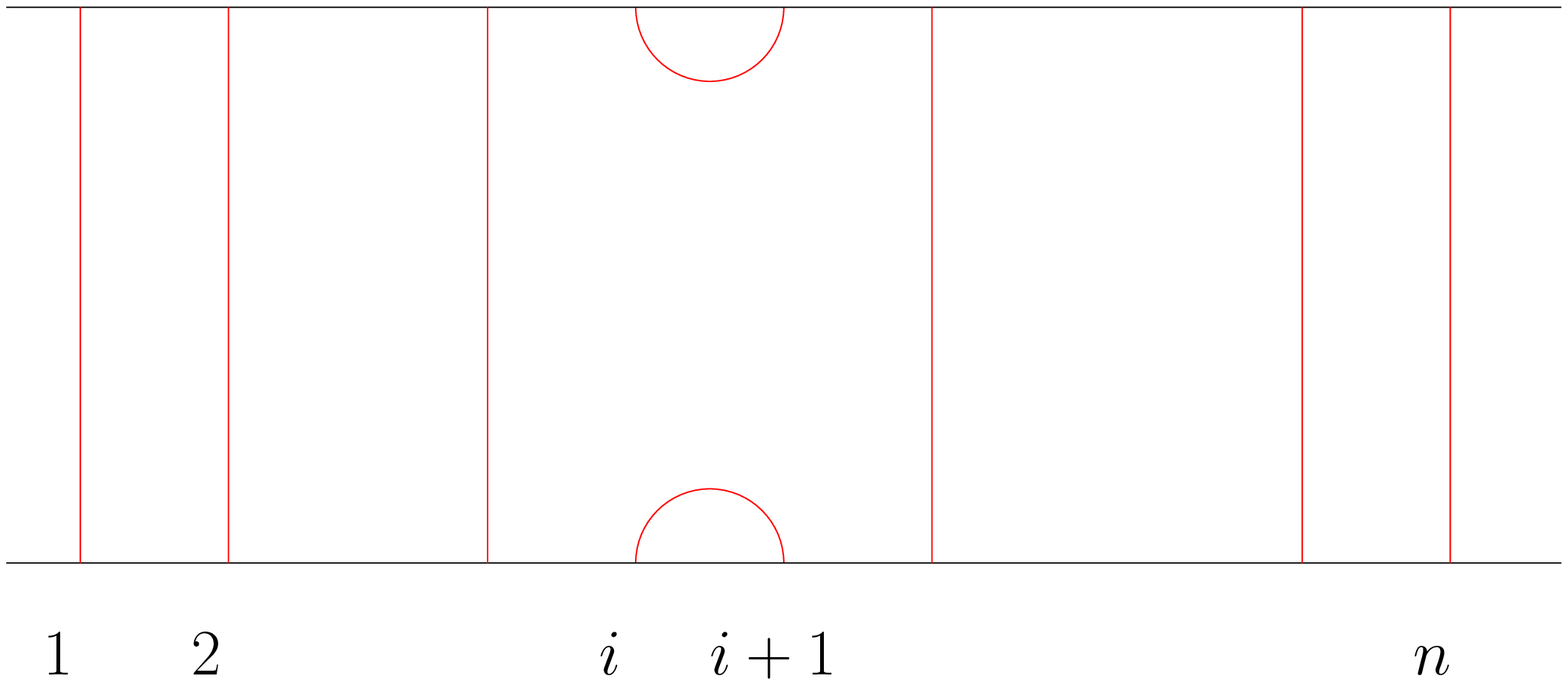}
\caption{Usual generators $U_i$ in Temperley--Lieb algebra.} \label{fig:9}
\end{figure}

Thus $T_{\pm,i}, T_{\pm,i}^*$ are algebraically identical to the operators of the same name in $\F$; and we may define
\[
U_{\pm,i} = T_{\pm,i}^* \; a_{\pm,i+1}.
\]
so that $U_{\pm,i}$ are operators $V(D^2)_n^e \To V(D^2)_n^e$ which glue annuli to the exterior of discs, having the effect shown in figure \ref{fig:8}. We see then (figure \ref{fig:9}) that the $U_{\pm,i}$ have sutures which are very similar to the usual generators given for the Temperley--Lieb algebra.

It follows immediately that the $U_{\pm,i}$ obey the relations of the Temperley--Lieb algebra (with $\delta = 0$), up to sign. In fact, since the $U_{\pm,i}$ are identical algebraically to the operators of the same name in section \ref{sec:Temperley-Lieb}, we have the ``twisted'' representation of the Temperley--Lieb algebra described there.

\subsection{Rotation}
\label{sec:rotation}

\begin{figure}
\centering
\includegraphics[scale=0.4]{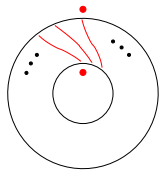}
\caption{The rotation operator.} \label{fig:10}
\end{figure}

The operation of rotating a chord diagram gives an operation in sutured TQFT. As discussed in \cite[section 7.1]{Me09Paper}, gluing a sutured annulus to the exterior of a disc, with sutures as in figure \ref{fig:10}, has the effect of rotating a sutured disc, or equivalently rotating the basepoint two places. (Moving the basepoint two places preserves the signs on either side of the basepoint.)

Consider computing an ``inner product'' $\langle \cdot | \cdot \rangle$, placing two chord diagrams into the ends of a cylinder. In particular, consider how sutures are rounded and which endpoints of the two sutured discs connect. We see that the set of sutures on the sphere obtained by inserting two chord diagrams $(\Gamma_0, \Gamma_1)$ is connected, if and only if the same is true for inserting the two chord diagrams $(\Gamma_1, R \Gamma_0)$; the rotation $R$ compensates precisely for the difference in rounding corners when the ends of the cylinder are swapped.

Thus we obtain an operator
\[
R \; : \; V(D^2)_n^e \To V(D^2)_n^e
\]
for each $n,e$, and for any chord diagrams $\Gamma_0, \Gamma_1$ of $n$ chords, $\langle c(\Gamma_0) | c(\Gamma_1) \rangle = \langle c(\Gamma_1) | R c(\Gamma_0) \rangle$. In particular, for any suture elements $c_0 \in c(\Gamma_0)$, $c_1 \in c(\Gamma_1)$ we obtain $\langle c_0 | c_1 \rangle = \pm \langle c_1 | R c_0 \rangle$. Recall the operator $H$ defined on $\F$; it follows that for suture elements $c_0, c_1$, $\langle c_0 | R c_1 \rangle = \pm \langle c_0 | H c_1 \rangle$.

We now choose a sign for $R$. For any $n,e$, we observe that rotation takes the basis chord diagram $\Gamma_{w_{min}}$ to $\Gamma_{w_{max}}$, and hence $R$ takes $c_{w_{min}}$ to $\pm c_{w_{max}}$; we choose a sign on $R$ so that $R c_{w_{min}} = c_{w_{max}}$.

\begin{prop}
With this choice of sign, $R=H$.
\end{prop}

\begin{Proof}
We first claim that for all $w \in \W_n^e$, $\langle c_w | R c_w \rangle = 1$. Let $A = \{ w \in \W_n^e \; : \; \langle c_w | R c_w \rangle = 1 \}$. Since $R c_{w_{min}} = c_{w_{max}}$, $w_{min} \in A$. The claim then follows obviously from the following fact:
\begin{itemize}
\item
If $w_0 \in A$, $w_0 \leq w_1$ and $w_0,w_1$ are related by an elementary move, then $w_1 \in A$.
\end{itemize}
To see why the fact is true, note that $c_{w_0} - c_{w_1}$ is a suture element and $\langle c_{w_0} | c_{w_0} - c_{w_1} \rangle = 0$. Thus $\langle c_{w_0} - c_{w_1} | R c_{w_0} \rangle = 0$ (maybe $\pm 0$, but nevertheless $0$!). Hence $\langle c_{w_0} | R c_{w_0} \rangle = \langle c_{w_1} | R c_{w_0} \rangle$ and as $w_0 \in A$ we have $\langle c_{w_1} | R c_{w_0} \rangle = 1$. Then we note that $\langle c_{w_0} - c_{w_1} | c_{w_1} \rangle = 0$, so that $\langle c_{w_1} | R ( c_{w_0} - c_{w_1} ) \rangle = 0$ also. This gives $\langle c_{w_1} | R c_{w_1} \rangle = \langle c_{w_1} | R c_{w_0} \rangle = 1$, hence $w_1 \in A$.

Next we claim that for any words $w_0, w_1$, $\langle c_{w_0} | c_{w_1} \rangle = \langle c_{w_1} | R c_{w_0} \rangle$. This is clearly true when $w_0 \nleq w_1$, so that both are $0$. So take a word $w \in \W_n^e$ given and let $A_w = \{ w' \in \W_n^e \; : \; \langle c_w | c_{w'} \rangle = \langle c_{w'} | R c_w \rangle \}$; we will show $A_w = \W_n^e$. From the previous claim, $w \in A_w$; and any $w'$ with $w \nleq w'$ lies in $A_w$. Thus we only need consider $w'$ with $w \leq w'$. Such a $w'$ can be reached from $w$ by a sequence of forwards elementary moves, and hence the claim follows from the following fact:
\begin{itemize}
\item
If $w \leq w_0 \leq w_1$, and $w_0 \in A_w$, and $w_0,w_1$ are related by an elementary move, then $w_1 \in A_w$.
\end{itemize}
To see that the fact is true, note that $c_{w_0} - c_{w_1}$ is a suture element, and $\langle c_w | c_{w_0} - c_{w_1} \rangle = 0$, hence $\langle c_{w_0} - c_{w_1} | R c_w \rangle = 0$. Thus $\langle c_{w_0} | R c_w \rangle = \langle c_{w_1} | R c_w \rangle$. But as $w_0 \in A_w$ then we have $\langle c_{w_0} | R c_w \rangle = \langle c_w | c_{w_0} \rangle = 1$. Hence $\langle c_{w_1} | R c_w \rangle = 1$, which is equal to $\langle c_w | c_{w_1} \rangle$.

Thus for any words, $\langle c_{w_0} | c_{w_1} \rangle = \langle c_{w_1} | R c_{w_0} \rangle$. As the $c_w$ form a basis for $V(D^2)_n^e$, it follows that for any $u,v \in V(D^2)_n^e$, $\langle u | v \rangle = \langle v | Ru \rangle$. Thus $R=H$.
\end{Proof}

Note it follows immediately that $H^{n+1} c_w = \pm c_w$ for every word $w \in \W_n^e$; and hence $H^{2n+2} = 1$. To obtain the precise sign $H^{n+1} = (-1)^{n_-  n_+}$, and precise period, proving theorem \ref{thm:H_periodic}, will take a little more work, which we delay to section \ref{sec:periodicity}.

It also follows that $R$ is an isometry: for any $u,v \in V(D^2)_n^e$, $\langle u | v \rangle = \langle R u | R v \rangle$. Of course this is obvious up to sign from the definition of $R$; the full equality is obvious now from the definition of $H$.

We have defined certain creation operators $a_{\pm,i}^*$, which have the effect on chord diagrams of inserting an outermost chord in a given position. These operators covered several locations for inserting chords, but not all. We showed in section \ref{sec:coherent_creation_and_annihilation} that creation operators are isometries; it's also clear they are injective. We now note that the operation of inserting an outermost chord \emph{anywhere} gives a linear operator which is an injective isometry: for it is a composition of some rotations and some $a_{\pm,0}^*$, which are all injective and isometries.

\subsection{Variations of nondegeneracy axioms}
\label{sec:variations_of_nondegeneracy}

We can now prove propositions \ref{prop:equivalent_axioms} and \ref{prop:inequivalent_axioms}, regarding variations of axiom \ref{ax:nondegeneracy}. Recall from section \ref{sec:results_mod_2} that given two words $w_0 \leq w_1$, $\pm[w_0,w_1]$ denotes the two suture elements of the chord diagram $[\Gamma_{w_0}, \Gamma_{w_1}]$; these have minimal and maximal basis elements $c_{w_0}, c_{w_1}$. And from section \ref{sec:partial_order}, recall the notion of difference $d(w_0, w_1)$, minimum and maximum of words. 

We will need the fact that any two words appear together in a suture element.
\begin{lem}
\label{lem:any_two_words_appear_together}
For any $w_0, w_1 \in \M_{n_-, n_+}$, both $c_{w_0}, c_{w_1}$ appear in the basis decomposition of the suture elements $\pm [\min(w_0,w_1), \max(w_0,w_1)]$. Moreover, in the decomposition of lemma \ref{lem:word_decomposition},
\[
\pm [\min(w_0,w_1), \max(w_0,w_1)] 
=
\pm [w_0^0, w_1^0] \; [w_1^1, w_0^1] \; [w_0^2, w_1^2] \; [w_1^3, w_0^3] \; \cdots [w_1^{2k-1}, w_0^{2k-1}].
\]
\end{lem}
(The multiplication of suture elements is multiplication in $V(D^2)$.)

\begin{Proof}
Suture elements are closed under multiplication, and there is a unique suture element (up to sign) with prescribed maximum and minimum. Now we know that in each $[w_0^i, w_1^i]$ or $[w_1^i, w_0^i]$, the basis elements $c_{w_0^i}, c_{w_1^i}$ each appear and are minimal/maximal. Hence, after multiplying, the basis elements $c_{w_0}, c_{w_1}, c_{\min(w_0,w_1)}, c_{\max(w_0,w_1)}$ all appear with $c_{\min(w_0,w_1)}$ minimal and $c_{\max(w_0,w_1)}$ maximal.
\end{Proof}

\begin{Proof}[of proposition \ref{prop:equivalent_axioms}]
The pairs of axioms (i)--(ii), (iii)--(iv), (v)--(vi), (vii)--(viii) are clearly equivalent for sets of sutures on $D^2$ with closed components, since axiom \ref{ax:6} then gives $c(\Gamma)=\{0\}$; for chord diagrams $\Gamma$, lemma \ref{lem:strong_sign_ambiguity} (which explicitly works for any of the 8 variations on axiom \ref{ax:nondegeneracy}) shows that $c(\Gamma)$ is of the form $\{x, -x\}$, immediately giving equivalence of the four pairs of axioms.

The implications $(ii) \Rightarrow (iv)$, and $(vi) \Rightarrow (viii)$ are obvious. Inclusion axiom (5 or 5') immediately gives $(viii) \Rightarrow (vi)$. Lemma \ref{lem:bypass_relation_pm} is precisely $(ii) \Rightarrow (v)$; in fact the argument shows $(iv) \Rightarrow (v)$. Thus we have the following, and it suffices to show $(v) \Rightarrow \{ (i), (ii) \}$.
\[
\{(i) \Leftrightarrow (ii)\} \Rightarrow \{(iii) \Leftrightarrow (iv)\} \Rightarrow \{(v) \Leftrightarrow (vi) \Leftrightarrow (vii) \Leftrightarrow (viii)\}.
\]

Assume $(v)$. Then we have a bypass relation. In fact, in developing the structure of $V(D^2)$ we have only used axiom \ref{ax:nondegeneracy} first to prove $(v)$ in lemma \ref{lem:bypass_relation_pm}, and then used $(v)$ repeatedly; we have not used $(i)$ or $(ii)$ independently. Thus, we obtain all the structure described in sections \ref{sec:creation_and_annihilation}--\ref{sec:rotation} above, including a coherent basis, creation and annihilation operators, and $\langle \cdot | \cdot \rangle$; all in isomorphism with the algebraic structure of $\F$. In particular $V(D^2) \cong \F$ as bigraded rings, and the bilinear forms $\langle \cdot | \cdot \rangle$ agree.

Now take two elements $\alpha, \beta \in V(D^2, F_n)$ with the property that for every chord diagram $\Gamma$ and one (hence any) $c \in c(\Gamma)$, $\langle \alpha | c \rangle = \pm \langle \beta | c \rangle$. We will show $\alpha = \pm \beta$, giving $(v) \Rightarrow \{ (i), (ii) \}$ and the desired equivalence of axioms. Using the isomorphism $V(D^2) \cong \F$, and the map $Q_-$ from section \ref{sec:duality}, we have
\[
\langle \alpha | c \rangle = \left( Q_-^{-1} \alpha \right) \cdot c, \quad \langle \beta | c \rangle = \left( Q_-^{-1} \beta \right) \cdot c.
\]
Hence, replacing $Q_-^{-1} \alpha, Q_-^{-1} \beta$ with $\alpha, \beta$, we have two elements $\alpha, \beta \in V(D^2, F_n)$ such that $\alpha \cdot c = \pm \beta \cdot c$ for all suture elements $c$; it is sufficient to show $\alpha = \pm \beta$. Consider now $c$ to run through basis elements $c_w$, for $w \in \W_n^e$. Writing the decompositions of $\alpha$ and $\beta$ as
\[
\alpha = \sum_{w \in \W_n^e} \alpha_w c_w, \quad \beta = \sum_{w \in \W_n^e} \beta_w c_w,
\]
we then have, for every $w \in \W_n^e$, $\alpha_w = \pm \beta_w$.

Partition $\W_n^e$ into three subsets:
\begin{align*}
S &= \{ w \in \W_n^e \; : \; \alpha_w = \beta_w = 0 \}, \\
T &= \{ w \in \W_n^e \; : \; 0 \neq \alpha_w = \beta_w \}, \\
U &= \{ w \in \W_n^e \; : \; 0 \neq \alpha_w = - \beta_w \}.
\end{align*}
We will show one of $T$ or $U$ is empty. Suppose not. Consider the set $A$ of suture elements $c \in V(D^2, F_n)$ whose basis decomposition contains words from both $T$ and $U$. By assumption, $T$ and $U$ are nonempty; then by lemma \ref{lem:any_two_words_appear_together} above, $A$ is nonempty. For $c \in A$, letting $\pm c = \pm [w_0, w_1]$, define $l(c) = d(w_0, w_1)$ (a difference as in definition \ref{defn:difference}). By definition $A$ contains no basis elements, so $l(c) > 0$.

Take $c \in A$ with $l(c)$ minimal; let $\pm c = \pm [w_0, w_1]$, so $w_0 \leq w_1$, $w_0 \neq w_1$. We claim that in the basis decomposition of $c$ there is only one word from $T$ and one from $U$. (There may be many from $S$.) Suppose not, so without loss of generality $c$ contains distinct words $t_0, t_1 \in T$ and $u \in U$. 

Let $w_- = \min(t_0, u)$ and let $w_+ = \max(t_0, u)$. By definition of $w_0$ we have $w_0 \leq t_0, u$ hence $w_0 \leq w_-$; indeed $w_0 \leq w_- \leq w_+ \leq w_1$. Thus $d(w_-, w_+) \leq d(w_0, w_1)$. By lemma \ref{lem:any_two_words_appear_together} above, $\pm [w_-, w_+]$ contains $t_0, u$ as monomials with nonzero coefficients; so $\pm [w_-, w_+] \in A$ and $l(\pm [w_-,w_+]) \leq l(\pm[w_0, w_1]) = l(c)$. By minimality of $c$ then equality must hold; so $w_0 = w_-$ and $w_+ = w_1$. In particular, $w_- \leq t_1 \leq w_+$. As $w_- = \min(t_0,u)$ and $w_+ = \max(t_0, u)$, then, for every $i \in \{1, \ldots, n_x \}$, the $i$'th $x$ in $t_1$ lies in a position between the respective positions of the $i$'th $x$ in $t_0$ and $u$:
\[
\text{either} \quad h_{t_0}^x (i) \leq h_{t_1}^x (i)  \leq h_u^x (i) \quad \text{or} \quad h_u^x (i) \leq h_{t_1}^x (i) \leq h_{t_0}^x (i).
\]

The same argument applies, reversing the roles of $t_0$ and $t_1$. We obtain that for every $i \in \{1, \ldots, n_x\}$,
\[
\text{either} \quad h_{t_1}^x (i) \leq h_{t_0}^x (i) \leq h_u^x (i) \quad \text{or} \quad h_u^x (i) \leq h_{t_0}^x (i) \leq h_{t_1}^x (i).
\]
Putting these together, we have that for all $i$, $h_{t_0}^x (i) = h_{t_1}^x (i)$, i.e. $t_0 = t_1$, a contradiction.

Thus the suture element $c \in V(D^2, F_n)$ contains precisely one word $t$ from $T$ and one $u$ from $U$, with coefficients $\epsilon_t, \epsilon_u \in \{-1,1\}$; the rest must be in $S$. Then $\alpha \cdot c = \alpha_{t} \epsilon_t + \alpha_{u} \epsilon_u$ and $\beta \cdot c = \beta_{t} \epsilon_t + \beta_{u} \epsilon_u$. By definition of $T$ and $U$, $\alpha_{t} = \beta_{t}$ and $\alpha_{u} = -\beta_{u}$, and recall $\alpha \cdot c = \pm \beta \cdot c$. Hence
\[
\alpha_{t} \epsilon_t + \alpha_{u} \epsilon_u = \left\{ \begin{array}{l} \alpha_{t} \epsilon_t - \alpha_{u} \epsilon_u \\ - \alpha_{t} \epsilon_t + \alpha_{u} \epsilon_u \end{array} \right.
\]
In the first case we have $2 \alpha_{u} \epsilon_u = 0$ so $\alpha_u = 0$ (these all lie in $\Z$, so there is no torsion); in the second $2 \alpha_{t} \epsilon_t = 0$ so $\alpha_t = 0$. These are both contradictions to the definitions of $T$ and $U$. Thus one of $T$ or $U$ is empty, and $\alpha = \pm \beta$.
\end{Proof}

Next, we prove proposition \ref{prop:inequivalent_axioms}, showing the equivalence of two variations on axiom \ref{ax:nondegeneracy}, and giving an example of a sutured TQFT satisfying this version of nondegeneracy. Recall from section \ref{sec:TQFT_axioms} that the ``free pre-sutured TQFT'' satisfies axioms 1--8, taking: $V(\Sigma, F) = \oplus_\Gamma \Z c_\Gamma$, a direct sum over isotopy classes of sutures $\Gamma$ on $(\Sigma,F)$ without contractible loops; $c(\Gamma) = c_\Gamma$ or $0$ as appropriate; and gluing maps defined in the most natural way. Our example will be a variation of this free pre-sutured TQFT, with a certain choice of bilinear form $\langle \cdot | \cdot \rangle$.

\begin{Proof}[of proposition \ref{prop:inequivalent_axioms}]
The equivalence of the two alternative axioms is immediate from axiom \ref{ax:8}, that each $V(\Sigma,F)$ is spanned by suture elements.

Define the sutured TQFT as follows: $V(\Sigma,F) = \oplus_\Gamma \Z c_\Gamma$ as for the free pre-sutured TQFT; so, for example, $V(D^2, F_n) \cong \Z^{C_n}$. For any set of sutures $\Gamma$ with contractible components,  $c(\Gamma) = \{0\}$; for $\Gamma$ without contractible components, $c(\Gamma) = \{c_\Gamma, - c_\Gamma\}$. For a gluing $\tau$ on a sutured background surface $(\Sigma,F)$, gluing maps $\Phi_\tau^i$ must satisfy $c_\Gamma \mapsto \pm c_{\#_\tau \Gamma}$ or $0$ as appropriate; we define \emph{many} gluing maps, covering all possibilities. More precisely, each $\Phi_\tau^i: V(\Sigma,F) \To V(\#_\tau(\Sigma,F))$ takes $c_\Gamma \mapsto 0$, if $\#_\tau \Gamma$ has a contractible component; else takes $c_\Gamma \mapsto \pm c_{\#_\tau \Gamma}$, where the sign can be chosen freely. As the $c_\Gamma$ form a free basis for $V(\Sigma,F)$, any choice of signs defines $\Phi_\tau^i$ completely and uniquely. We allow all choices of signs, and these form our collection $\Phi_\tau^i$. It is clear that axioms 1--8 are satisfied.

Now, as described in section \ref{sec:TQFT_axioms}, the bilinear map $\langle \cdot | \cdot \rangle: V(D^2, F_n) \otimes V(D^2, F_n) \To \Z$ is a gluing map obtained from stacking. But we have many sign choices when we define $\langle \cdot | \cdot \rangle$; we will show that these sign choices can be made to satisfy alternative axiom (i), proving the result. For two sets of sutures $\Gamma_0, \Gamma_1$ on $(D^2, F_n)$, $\langle c_{\Gamma_0} | c_{\Gamma_1} \rangle$ must be $0$, if $\Gamma_0, \Gamma_1$ are not stackable; otherwise $\langle c_{\Gamma_0} | c_{\Gamma_1} \rangle = \pm 1$, and moreover, as the $c_{\Gamma}$ form a free basis, either sign is available to define $\langle \cdot | \cdot \rangle$; signs may be chosen independently for different chord diagrams.

Ordering arbitrarily, write $\Gamma_1, \Gamma_2, \ldots \Gamma_{C_n}$ for the chord diagrams on $(D^2, F_n)$. Let $g_{ij} = \langle c_{\Gamma_i} | c_{\Gamma_j} \rangle$. For each pair $(i,j)$, either $g_{ij}$ is forced to be zero, or we are free to choose $g_{ij} = -1$ or $1$. Let $g$ be the matrix with entries $g_{ij}$. Clearly $g$ is the matrix for a bilinear map $\Z^{C_n} \otimes \Z^{C_n} \To \Z$; but equally it is the matrix for a bilinear map $\Q^{C_n} \otimes \Q^{C_n} \To \Q$. As $V(D^2, F_n)$ has a basis of suture elements, and these taken together with their negatives and $0$ form the complete set of suture elements in $V(D^2, F_n)$, (i) holds if the matrix $g$ over $\Q$ has full rank.

Thus we only need to choose signs for the $g_{ij}$ so that $g$ has full rank over $\Q$. This can be done as follows. For $i \neq j$, choose $g_{ij}$ arbitrarily (if there is a choice!) Now for any chord diagram $\Gamma$, $\langle c_\Gamma | c_\Gamma \rangle = \pm 1$, as we have seen above; thus each diagonal element $g_{ii} = \pm 1$. Consider row-reducing $g$. We successively choose the $g_{ii}$ so that, as row reduction proceeds to the $i$'th line, the $(i,i)$ element remains nonzero. For the first line this is obvious. Now suppose it is true up to the $i$'th row; the row-reducing ensures that at the $i$'th row, all $(i,j)$ elements with $j<i$ are zero. However, the $(i,i)$ element began as $g_{ij} = \pm 1$; in row reducing to ensure $(i,j)$ elements are zero for $j<i$, amounts are added to this $\pm 1$. For at least one of the two possible sign choices, the result will be nonzero; choose this sign. Then the $(i,i)$ element remains nonzero. After row reduction, we have an upper triangular matrix with nonzero entries on the diagonal; hence it has full rank, and (i) is satisfied.
\end{Proof}

\subsection{An additional axiom}
\label{sec:additional_axiom}

We consider a tenth axiom for sutured TQFT. Note that a gluing $\tau$ of a sutured background surface can only increase genus; thus we could simply set $V(\Sigma,F) = 0$ whenever $\Sigma$ has genus at least $1$, set all gluing maps to higher genus surfaces to be $0$, and we would obtain a consistent theory. Our additional axiom will require that certain gluing maps be isomorphisms. Recall that a sutured gluing map arises from an identification $\tau$ of arcs on the boundary of a sutured background surface $(\Sigma, F)$, which respects marked points $F$ and signs of complementary arcs.
\begin{ax}
\label{ax:10}
Let $\tau$ be a sutured gluing map on $(\Sigma,F)$, identifying two disjoint arcs $\gamma, \gamma'$ on $\partial \Sigma$. Suppose that $|\gamma \cap F| = |\gamma' \cap F| = 1$. Then any gluing map $\Phi_\tau$ associated to $\tau$ is an isomorphism.
\end{ax}

Any connected sutured background $(\Sigma, F)$ can be constructed from a sutured background disc $(D^2, F_n)$ by gluing maps of this type; in fact, any sutured background $(\Sigma, F)$ (possibly disconnected) can be constructed from a disjoint union of sutured background discs $\sqcup (D^2, F_2)$ by gluing maps of this type. So this final axiom strong enough to construct isomorphisms $V(\Sigma,F) \cong V(\sqcup_i (D^2, F_2)) \cong \otimes_i V(D^2, F_2)$ for any $(\Sigma, F)$. As we now see, it makes axiom 9 redundant, and almost does the same to axiom 8.
\begin{lem}
\label{lem:1-8_plus_10_implies_9}
In the presence of axioms \ref{ax:1}--\ref{ax:8}, axiom \ref{ax:10} implies axiom \ref{ax:nondegeneracy}.
\end{lem}

\begin{Proof}
On $(D^2, F_2)$ there are only two chord diagrams; as mentioned in section \ref{sec:basis_partial_order} inclusions $(D^2, F_2) \hookrightarrow (D^2, F_1)$ with intermediate sutures can easily be found to show that they have suture elements which are linearly independent; by axiom \ref{ax:8} two such suture elements span $V(D^2, F_1)$; thus they form a basis which we denote $x,y$ and in fact $V(D^2, F_2) = V(D^2)_1 \cong \F_1$ as graded abelian groups. 

Gluing two such discs, we obtain an isomorphism $V(D^2, F_3) = V(D^2)_2 \cong \F_2$ as graded abelian groups, and the summand $V(D^2)_2^0$ has basis $\{x \otimes y, y \otimes x \}$, which are suture elements (exactly like $xy, yx$ in the foregoing) for two bypass-related chord diagrams of $3$ chords. Including $(D^2, F_3) \hookrightarrow (D^2, F_1)$ in various ways, we obtain maps $\Z^2 \cong V(D^2)_2^0 \To V(D^2, F_1) \cong \Z$ which are coordinate projections; from these we see that the third chord diagram in their bypass triple is $\pm x \otimes y \pm y \otimes x$ (exactly like $\pm xy \pm yx$ in the foregoing). By proposition \ref{prop:equivalent_axioms}, this is equivalent to axiom \ref{ax:nondegeneracy}.
\end{Proof}

\begin{lem}
\label{lem:1-7_plus_10_implies_8}
Assume axioms \ref{ax:1}--\ref{ax:7}. Suppose on $(D^2, F_2)$ the two chord diagrams, one each in Euler class $-1$ and $1$, have suture elements $\pm x$, $\pm y$ respectively where $\{x,y\}$ is a basis for $V(D^2, F_2)$. Then axiom \ref{ax:10} implies axiom \ref{ax:8}. (And hence, by lemma \ref{lem:1-8_plus_10_implies_9}, also axiom \ref{ax:nondegeneracy}.)
\end{lem}

\begin{Proof}
Gluing several discs $(D^2, F_2)$ together, base-to-root, gives isomorphisms by axiom \ref{ax:10}. Since $x,y$ form a basis for $V(D^2, F_2)$, we obtain a basis of suture elements $\{x,y\}^{\otimes n}$ for each $V(D^2, F_{n+1})$ (exactly like words in $x,y$ in the foregoing). Thus each $V(D^2, F_n)$ is spanned by suture elements, and hence by gluing we obtain a basis of any $V(\Sigma,F)$ of suture elements; in particular a spanning set.
\end{Proof}

\section{Sutured Floer homology and sutured TQFT}
\label{sec:SFH}

We now show that the sutured Floer homology, defined by Juh\'{a}sz in \cite{Ju06}, of certain balanced sutured manifolds forms a sutured TQFT. As noted in \cite{HKM08}, the $SFH$ of manifolds of the type $(\Sigma \times S^1, F \times S^1)$, where $\Sigma$ is a surface with nonempty boundary and $F \subset \partial \Sigma$ is finite, has TQFT-like properties. In fact it is a sutured TQFT; for which sutured TQFT was designed to be an axiomatic model. Make the following assignments.
\begin{itemize}
\item
To a sutured background surface $(\Sigma,F)$, assign the abelian group $V(\Sigma,F) = SFH(- \Sigma \times S^1, - F \times S^1)$ with $\Z$ coefficients. It is known that $SFH$ splits as a direct sum over spin-c structures.
\item
A set of sutures $\Gamma$ on $(\Sigma, F)$ corresponds precisely to an isotopy class of contact structures $\xi$ on $\Sigma \times S^1$, such that the boundary $\partial \Sigma \times S^1$ is convex with dividing set $F \times S^1$ and positive/negative regions determined by the decomposition of $\partial \Sigma \backslash F$ into  positive and negative arcs $C_+ \cup C_-$ \cite{Gi00, GiBundles, Hon00II}. Let $c(\Gamma)$ be the \emph{contact invariant} $c(\xi) \subset V(\Sigma,F)$ \cite{OSContact, HKMContClass, HKM06ContClass}. This $c(\xi)$ is a subset of the form $\{\pm x\}$. The possible relative Euler classes of $\xi$ correspond to the spin-c structures on $(\Sigma,F)$; and $c(\xi)$ lies in the corresponding spin-c summand of $SFH$.
\item
For a gluing $\tau$ of the sutured background surface $(\Sigma,F)$, let $\Phi_\tau: SFH(- \Sigma \times S^1, - F \times S^1) \To SFH( - (\#_\tau \Sigma) \times S^1, -(\#_\tau F) \times S^1)$ be the map defined in \cite{HKM08} by the obvious inclusion of $\Sigma \times S^1 \hookrightarrow \#_\tau \Sigma \times S^1$, together with the canonical contact structure on $\#_\tau \Sigma \times S^1 - \Sigma \times S^1$ as convex neighbourhood of the boundary. In fact we can choose a sign on each $\Phi_\tau$ on each Euler class summand; let $\Phi_\tau^i$ be the collection of all maps obtained from $\Phi_\tau$ by all possible choices of signs.
\end{itemize}

\begin{prop}
\label{prop:SFH_is_STQFT}
These assignments satisfy all axioms of sutured TQFT (including axiom \ref{ax:10}).
\end{prop}

\begin{Proof}
The assignments above clearly satisfy axioms \ref{ax:1}--\ref{ax:3}. Axiom \ref{ax:4}, that $V(\sqcup_i (\Sigma_i, F_i)) = \otimes_i V(\Sigma_i, F_i)$ is clear from the definition of $SFH$. Axiom \ref{ax:5}, that a gluing map takes contact elements to contact elements in a natural way, is proved in \cite{HKM08}; as the contact invariant is of the form $\{\pm x\}$, each $\Phi_\tau^i$ is surjective on contact elements. If $\Gamma$ contains a contractible loop then the corresponding contact structure is overtwisted; by \cite{HKM06ContClass} then $c(\Gamma) = \{0\}$, so axiom \ref{ax:6} holds. That $V(D^2, F_1) \cong \Z$ is proved in \cite{HKM08} and follows from \cite{Ju08} or \cite{HKM06ContClass}; that $c(\Gamma_\emptyset) = \{\pm 1\}$ is proved in \cite{HKM08} and follows from \cite{HKM06ContClass}; so axiom \ref{ax:7} holds.

Axiom \ref{ax:10}, that a gluing map identifying precisely one pair of marked points is an isomorphism, is proved in \cite{HKM08}. Moreover, it's shown in \cite{HKM08}, following \cite{HKM06ContClass}, that $V(D^2, F_2) \cong \Z^2$, which splits as a direct sum $\Z \oplus \Z$ corresponding to chord diagrams of Euler class $-1$ and $1$; and contact elements for the two chord diagrams form a basis. By lemma \ref{lem:1-7_plus_10_implies_8}, this implies axiom \ref{ax:8}; by lemma \ref{lem:1-8_plus_10_implies_9}; this implies axiom \ref{ax:nondegeneracy}. 
\end{Proof}

It follows immediately that any result about sutured TQFT immediately gives a result about $SFH$. In particular, contact elements in the $SFH$ of sutured solid tori with longitudinal sutures have all the structure of the Fock space of section \ref{sec:algebra}, and all the structure of the sutured TQFT of discs.

\begin{cor}
\label{cor:STQFT_exists}
A sutured TQFT exists.
\qed
\end{cor}

Corollary \ref{cor:STQFT_exists} could of course also be proved by formally constructing a sutured TQFT from scratch. For instance, one could define $V(D^2, F_n)$ to be a direct sum of $\Z$ summands, one for each basis chord diagram; and then specify suture elements using the bypass relation and decomposition into a basis; and then go through the procedure of section \ref{sec:sutured_TQFT} to iron out signs. Defining the sutured TQFT on more complicated surfaces would however take more work.

\section{Non-commutative QFT = Sutured TQFT of discs}
\label{sec:sutured_TQFT_discs}

\subsection{Main isomorphism and suture elements}
\label{sec:main_isomorphism}

After section \ref{sec:sutured_TQFT}, it is clear that much of the algebraic non-commutative QFT structure developed in section \ref{sec:algebra} is equivalent to structure in the sutured TQFT of discs. To summarise, we make a detailed statement that ``the sutured TQFT of discs is the QFT of two non-commuting particles'', proving remaining details. This includes earlier theorems \ref{thm:main_thm_rough} and  \ref{thm:mini_main_isomorphism}.
\begin{thm}
\label{thm:main_isomorphism}
Every sutured TQFT obeying axioms 1--9 above satisfies
\[
V(D^2) \cong \F,
\]
an isomorphism of graded rings. In particular
\[
V(D^2)_{n+1} \cong \F_n \cong \Z^{2^n}, \quad V(D^2)_{n_-,n_+} \cong \F_{n_x, n_y} \cong \Z^{\binom{n}{n_x}} \cong \Z^{\binom{n}{n_y}},
\]
isomorphisms of abelian groups. Each $c_w \in V(D^2)_n^e$, for $w \in \W_n^e$, corresponds to $w \in \M_n^e \subset \F_n^e$, replacing $(-,+)$ with $(x,y)$. Moreover, under this isomorphism:
\begin{enumerate}
\item
The action of annihilation and creation operators $a_{x,i}, a_{y,i}, a_{x,i}^*, a_{y,i}^*$ on words $w \in \M_n^e \subset \F_n^e$ is identical to the action of annihilation and creation operators $a_{-,i}, a_{+,i}, a_{-,i}^*, a_{+,i}^*$ on basis suture elements $c_w$, for $w \in \W_n^e$.
\item
The $\langle \cdot | \cdot \rangle$ defined on $\F$ via the partial order $\leq$, and the $\langle \cdot | \cdot \rangle$ defined on $V(D^2)$ by stacking, agree.
\item
The operators $T_{x,i}$, $T_{y,i}$, $T_{x,i}^*$, $T_{y,i}^*$, $U_{x,i}$, $U_{y,i}$ act on each $\F_n^e$ identically to the operators $T_{-,i}$, $T_{+,i}$, $T_{-,i}^*$, $T_{+,i}^*$, $U_{-,i}$, $U_{+,i}$ on $V(D^2)_n^e$.
\item
The duality operator $H$ acts on $\F_n^e$ identically to the rotation operator $R$ on $V(D^2)_n^e$.
\item
The set of suture elements in $V(D^2)_n^e$ maps to the distinguished subset $\mathcal{C}_n^e$.
\end{enumerate}
\end{thm}

Most of this theorem has already been proved. The isomorphism of rings is clear; we chose signs on basis elements $c_w \in V(D^2)$, on $\langle \cdot | \cdot \rangle$, on multiplication, and on annihilation and creation and other operators, so that the isomorphism might apply to them; and then we proved that the isomorphism does indeed apply. We chose signs and proved $H=R$. It only remains to prove the last statement (v). But to do this we need to establish that $\mathcal{C}_n^e$ exists, in particular proposition \ref{prop:C1_equals_C2} that the three definitions $\mathcal{C}^1$, $\mathcal{C}^2$, $\mathcal{C}^3$ all agree. We will prove proposition \ref{prop:C1_equals_C2} and theorem \ref{thm:main_isomorphism} together by showing that all of $\mathcal{C}^1, \mathcal{C}^2, \mathcal{C}^3$ map to the set of sutures elements under this isomorphism. (Recall $\mathcal{C}^1 = \left\{ a_{s,i}^*, a_{s,i}, T_{s,i}^* \right\} \cdot 1$, $\mathcal{C}^2 = \left\{ a_{s,0}^*, a_{s,0}, H \right\} \cdot 1$ and $\mathcal{C}^3 = \left\{ a_{s, n_s + 1}^*, a_{s, n_s + 1}, H \right\} \cdot 1$.)

It's interesting to note that our proof of proposition \ref{prop:C1_equals_C2} relies upon an isomorphism to sutured TQFT, and hence upon the existence of a sutured TQFT; thus, at least in our presentation, there is a dependence upon the constructions of sutured Floer homology, holomorphic curves and all.

\begin{Proof}[Of proposition \ref{prop:C1_equals_C2} and theorem \ref{thm:main_isomorphism}(v)]
We first show that under our isomorphism, $\mathcal{C}^1$ corresponds to suture elements. Clearly the class of suture elements is preserved under the action of creations, annihilations, and the operators $T_{\pm,i}^*$, since we explained their action on suture elements in sections \ref{sec:creation_and_annihilation} and \ref{sec:Temperley-Lieb_STQFT} above. We must show that every suture element can be created from the vacuum by the action of these operators. Clearly the $0$ suture element can be obtained: annihilating the vacuum, for instance; we first show that we can obtain a suture element for any chord diagram.

Proof by induction on the number of chords. Clearly we can obtain the vacuum. Now suppose we have a chord diagram $\Gamma$ of $n+1$ chords and Euler class $e$. If $\Gamma$ contains an outermost chord at the base point, or an outermost chord at the root point, or has an outermost chord enclosing a positive region on the eastside, or has an outermost chord enclosing a negative region on the westside, then $\Gamma$ is obtained from a smaller chord diagram by applying a creation operator, and we reduce to a smaller diagram. So we may assume all outermost regions are negative and on the eastside, or positive and on the westside.

Suppose there is an outermost negative region on the eastside; the case of an outermost positive region on the westside is similar. Then there is one closest to the base point, so is enclosed by a chord running from the point $2i+1$ to $2i+2$; and this $i$ is minimal. (The points of $F_{n+1}$ on $(D^2, F_{n+1})$ are numbered as mentioned in section \ref{sec:sutured_surfaces} above.) Then, as shown in figure \ref{fig:11}, the chord diagram $\Gamma$ is obtained from a chord diagram $\Gamma'$ by applying the operator $a_{+,i+2} T_{+,i}^*$, where $\Gamma'$ is identical to $\Gamma$, except that the outermost chord between $(2i+1,2i+2)$ is moved to $(2i-1,2i)$. Thus the outermost negative region in $\Gamma'$ is closer to the base point; applying this procedure finitely many times, we obtain a chord diagram with an outermost chord at the base point, and can reduce to a smaller one.

\begin{figure}
\centering
\includegraphics[scale=0.3]{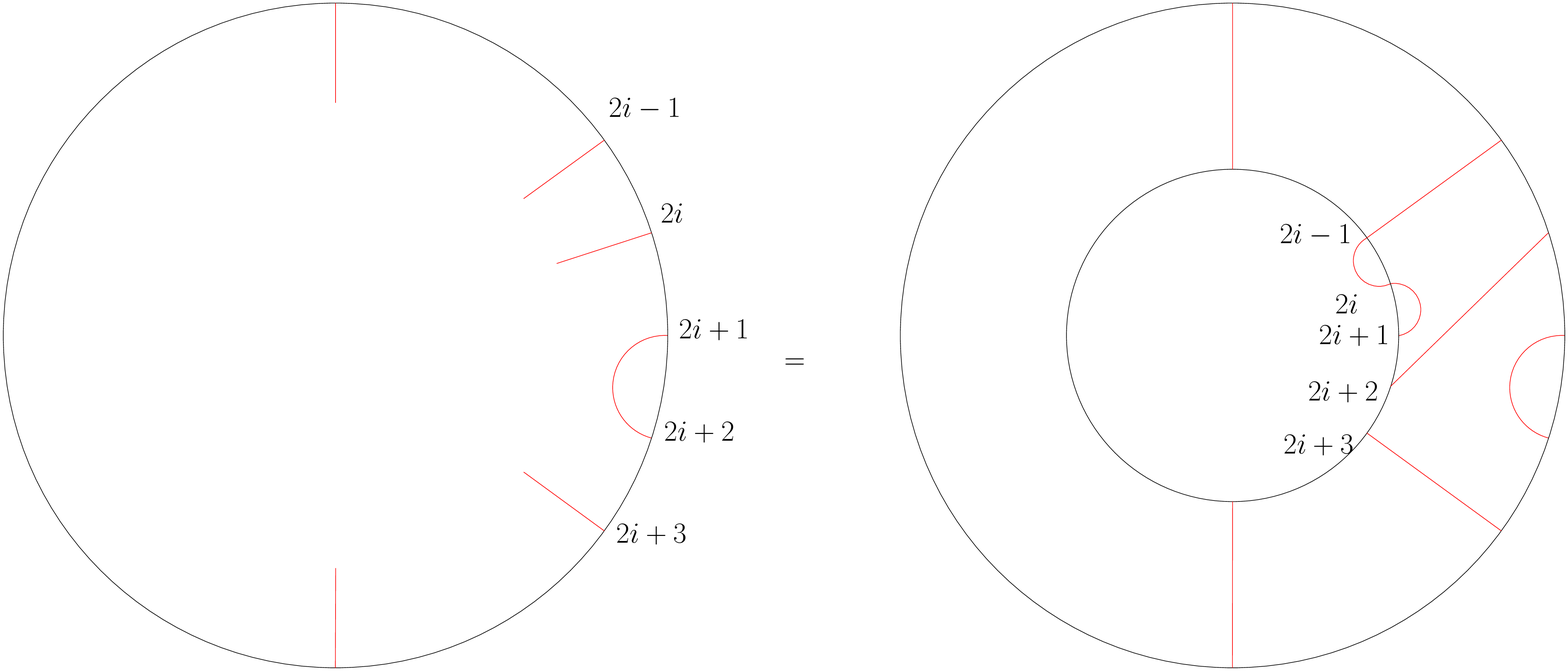}
\caption{A chord diagram with outermost region at $(2i+1,2i+2)$ is $a_{+,i+2} T_{+,i}^*$ of an otherwise identical chord diagram with outermost region at $(2i-1,2i)$.} \label{fig:11}
\end{figure}

Thus, for any chord diagram $\Gamma$, we can obtain a suture element in $c(\Gamma)$ by applying creation, annihilation and $T_{\pm,i}^*$ operators; but we have not yet shown we can obtain \emph{both} elements in $c(\Gamma)$. For this it suffices to show we can obtain $-1$. We can obtain $-1$ from $1$, for instance, by operating
\[
1 \quad \stackrel{a_{x,0}^*}{\mapsto} \quad x \quad \stackrel{T_{y,0}^*}{\mapsto} \quad yx - xy \quad \stackrel{a_{x,0}}{\mapsto} \quad -y \quad \stackrel{a_{y,0}}{\mapsto} \quad -1.
\]
Thus, $\mathcal{C}^1$ coincides with the set of suture elements. 

We now turn to $\mathcal{C}^2$. A similar but easier argument applies. Again, the class of suture elements is clearly preserved by creations and by $H$; we must show that we can obtain any suture element from the vacuum by applying initial creation and rotation operators. We can clearly obtain $0$. Again, proof by induction on the number of chords. One chord is clear. Let $\Gamma$ be a chord diagram. Then $\Gamma$ has an outermost chord; applying $H$ this may be rotated to the base point; so $\Gamma$ is obtained from a smaller chord diagram by an initial creation followed by some applications of $H$.

This shows that for any $\Gamma$ we can obtain an element in $c(\Gamma)$. We now show we can obtain $-1$. For this we need to use the fact that $H^{n+1} = -1$ sometimes; we have not proved this yet, but for our purposes it is sufficient to note that on $\F_{1,1}$, $H^3 = -1$; this can be computed by hand from the formula in corollary \ref{cor:H_formula}. Thus we can obtain $-1$, for instance, by
\[
1 \quad \stackrel{a_{y,0}^*}{\mapsto} \quad y \quad \stackrel{a_{x,0}^*}{\mapsto} \quad xy \quad \stackrel{H^3}{\mapsto} \quad -xy \quad \stackrel{a_{x,0}}{\mapsto} \quad -y \quad \stackrel{a_{y,0}}{\mapsto} \quad -1.
\]
Thus $\mathcal{C}^2$ coincides with the set of suture elements. The proof that $\mathcal{C}^3$ coincides with suture elements is identical, replacing base point with root point.
\end{Proof}

We can now also prove the algebraic statements about suture elements in theorem \ref{thm:properties_of_suture_elements}. (Again, this proof relies on sutured TQFT existence.)
\begin{Proof}[of theorem \ref{thm:properties_of_suture_elements}]
Part (i) is clear: we have defined $\times, a_{s,i}, a_{s,i}^*, T_{s,i}, T_{s,i}^*, U_{s,i}$ to preserve suture elements; and we have considered signs of suture elements at length. Part (ii) is a statement of the properties of suture elements, as discussed in sections \ref{sec:results_mod_2} and \ref{sec:coherent_creation_and_annihilation} above; including (ii) of proposition \ref{prop:suture_elt_structure}. Part (iii) is clear: it's easy to check $Q_\pm$ does not preserve $\mathcal{C}$; we have computed an inverse for $Q_\pm$, so that it is a bijection, and $|Q_\pm \mathcal{C}_n^e|  = |\mathcal{C}_n^e | = 2 N_n^e$; and since $Q_\pm \circ Q_\mp^{-1} = H^{\pm 1}$ preserves $\mathcal{C}$, we must have $Q_+ \mathcal{C} = Q_- \mathcal{C}$. Part (iv) follows from the above: $\langle v | v \rangle = 1$ is (iii) of proposition \ref{prop:suture_elt_structure}; $\langle v | Hv \rangle = 1$ then follows from the duality definition of $H$; and $\langle v_0 | v_1 \rangle \in \{-1,0,1\}$ follows from the sutured-TQFT definition of $\langle \cdot | \cdot \rangle$.

For part (v), note that if the sum/difference of two suture elements $u,v$ is also a suture element, then the same is true mod 2, hence by proposition 1.10 of \cite{Me09Paper}, they correspond to bypass-related chord diagrams $\Gamma_u,\Gamma_v$. Reorder $u,v$ if necessary so that $\langle u | v \rangle = \pm 1$ and $\langle v | u \rangle = 0$, and switch signs if necessary so $\langle u | v \rangle = 1$; then $u-v$ is a suture element for the third diagram. Since $u,v$ are bypass-related, $\Gamma_u, \Gamma_v$ may be isotoped to be identical, except in a disc $D' \subset D$; and we may choose a base point on $D'$ such that $\Gamma_u \cap D'$ is the basis chord diagram $\Gamma_{xy}$, and $\Gamma_v \cap D'$ is the basis chord diagram $\Gamma_{yx}$; these have suture elements $xy,yx$ respectively, and $\langle xy | yx \rangle = 1$. Now the chord diagrams $\Gamma_u, \Gamma_v$ can be obtained by applying initial creation operators and rotations to these two diagrams, so we obtain an operator $A^*$ taking $xy \mapsto \pm u$ and $yx \mapsto \pm v$. Since creations and rotations are isometries, we have $1 = \langle xy | yx \rangle = \langle A^*(xy) | A^*(yx) \rangle = \langle u | v \rangle$; thus under $A^*$, $(xy, yx) \mapsto \pm(u,v)$. If we get $-(u,v)$, then precompose $A^*$ with $H^3 = -1$ on $\F_{1,1}$, and we have the desired $A^*$.

Part (vi) is essentially lemma 3.1 of \cite{Me09Paper} without the contact geometry and with signs; we repeat the argument here with modifications necessary in this context. Given two chord diagrams $\Gamma_u, \Gamma_v$ and $u \in c(\Gamma_u)$, $v \in c(\Gamma_v)$ with $\langle u | v \rangle = 1$, we must find chord diagrams $\Gamma_u = \Gamma_0, \Gamma_1, \ldots, \Gamma_m = \Gamma_v$, and $u = c_0, c_1, \ldots, c_m = v$ where $c_i \in c(\Gamma_i)$, such that each $c_i - c_{i+1}$ is a suture element --- in particular, each pair $\Gamma_i, \Gamma_{i+1}$ is bypass-related --- and for all $i \leq j$, $\langle c_i | c_j \rangle = 1$. We prove this by induction on the number of chords in $\Gamma_u, \Gamma_v$. The assumption $\langle u | v \rangle = 1$ implies that $\Gamma_u, \Gamma_v$ have the same Euler class. With less than three chords there is nothing to prove; with three chords $\Gamma_u, \Gamma_v$ are either identical or bypass-related, and the result is clear.

Now consider general $\Gamma_u, \Gamma_v$ and $u,v$ with $\langle u | v\rangle = 1$. If $\Gamma_u, \Gamma_v$ share an outermost chord $\gamma$ then we simply consider $\Gamma_u - \gamma$ and $\Gamma_u - \gamma$; by induction we have a sequence of bypass-related chord diagrams with the desired properties; adding $\gamma$ to all these is an operation which gives an isometry, as noted at the end of section \ref{sec:rotation}. Hence we obtain chord diagrams and suture elements with the desired properties.

Thus we may assume $\Gamma_u, \Gamma_v$ have no outermost chords in common. Let $\gamma$ be an outermost chord of $\Gamma_v$, let its endpoints be $p,q$ where $q$ is clockwise of $p$. Let the next marked point of $\Gamma_v$ clockwise of $q$ be $r$. By assumption, there is no outermost chord connecting $p,q$ on $\Gamma_u$; there is also no outermost chord connecting $q,r$ on $\Gamma_u$, since then rounding would give $\langle u | v \rangle = 0$. Thus $\Gamma_u, \Gamma_v$ appear as shown in figure \ref{fig:12}, and we may perform upwards bypass surgery on $\Gamma_u = \Gamma_0$ near $p,q,r$, along the attaching arc $\delta$ shown, to obtain $\Gamma_1$.

\begin{figure}
\centering
\includegraphics[scale=0.4]{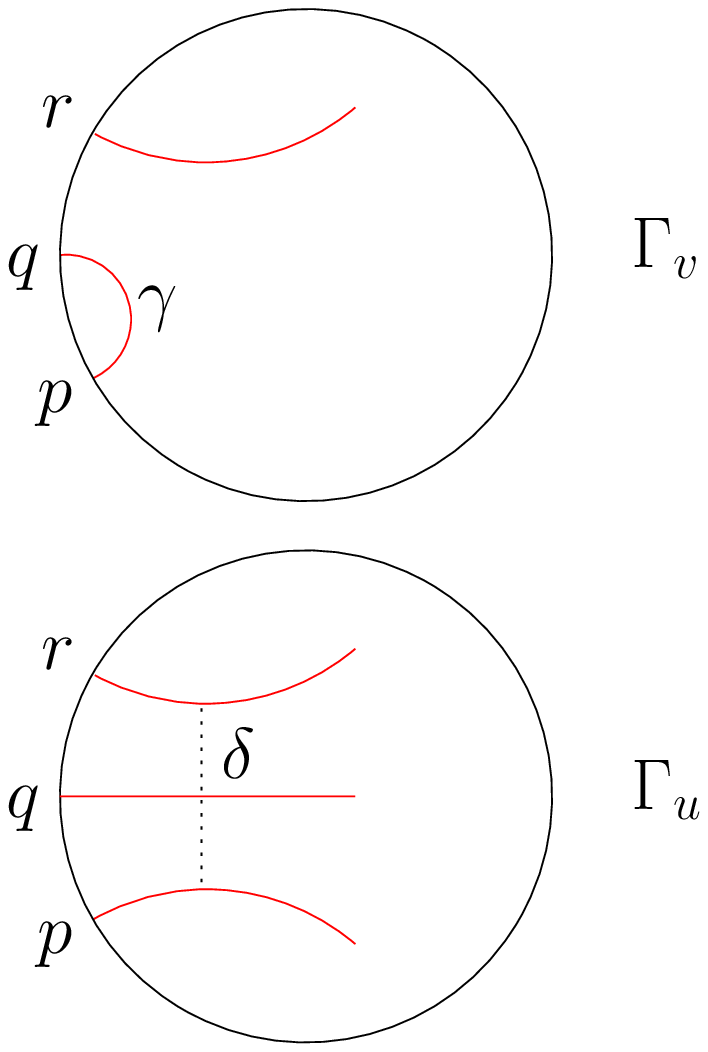}
\caption{Chord diagrams $\Gamma_u, \Gamma_v$.} \label{fig:12}
\end{figure}

Now, $\Gamma_1, \Gamma_v$ share the common outermost chord $\gamma$. Removing it, by induction we obtain a sequence of chord diagrams $\Gamma_1 - \gamma = \Gamma'_1, \Gamma'_2, \ldots, \Gamma'_m = \Gamma_v - \gamma$, and $c'_i \in c(\Gamma'_i)$, with each $c'_i - c'_{i+1}$ a suture element --- hence each pair $\Gamma'_i, \Gamma'_{i+1}$ bypass-related --- and $\langle c'_i | c'_j \rangle = 1$ for $1 \leq i \leq j \leq m$. Reinserting $\gamma$ we obtain a sequence of chord diagrams $\Gamma_1, \Gamma_1, \ldots, \Gamma_m$, where each $\Gamma_i, \Gamma_{i+1}$ are bypass-related. This inserting of $\Gamma$ gives a linear operator $L$ which is an isometry, and we may adjust $L$ by sign if necessary so that $L(c'_m) = v$. Then for $1 \leq i \leq m$, let $c_i = L(c'_i)$. As $L$ is an isometry, we have $\langle c_i | c_j \rangle = 1$ for $1 \leq i \leq j \leq m$. We set $c_0 = u$. As $L$ preserves suture elements, each $c_i - c_{i+1} = L(c'_i - c'_{i+1})$ for $1 \leq i \leq m-1$ is a suture element. 

We claim that the $\Gamma_i$ and $c_i$ are sequence of chord diagrams and suture elements with the desired properties. It only remains to verify that $c_0 - c_1$ is a suture element and that $\langle c_0 | c_i \rangle = 1$ for $0 \leq i \leq m$. Clearly $\langle c_0 | c_0 \rangle = 1$. To see $c_0 - c_1$ is  a suture element, note that since $\Gamma_0, \Gamma_1$ are bypass-related, either $c_0 - c_1$ or $c_0 + c_1$ is a suture element. Since $\langle c_0 | c_m \rangle = \langle u | v \rangle = 1$ and $\langle c_1 | c_m \rangle = 1$, we have $\langle c_0 \pm c_1 | c_m \rangle = 2$ or $0$; clearly $2$ is impossible, so it must be $c_0 - c_1$ that is a suture element.

For any $1 \leq i \leq m$, we next show that $\langle c_0 | c_i \rangle = \pm 1$. Since $\Gamma_1$ is obtained from $\Gamma_0$ by an upwards bypass surgery, $\langle c_0 | c_1 \rangle = \pm 1$. Now for $2 \leq i \leq m$, consider the sutured cylinders which are rounded in computing $\langle c_0 | c_i \rangle$ and $\langle c_1 | c_i \rangle$. Note $\Gamma_i$ contains the outermost chord $\gamma$, which can be pushed down by a ``finger move'' to $\Gamma_0$ or $\Gamma_1$ to give identical sets of sutures. Thus if one of these sutured cylinders has connected sutures, so does the other. Since $\langle c_1 | c_i \rangle = 1$ then $\langle c_0 | c_i \rangle = \pm 1$.

Finally we verify that for each $1 \leq i \leq m$, $\langle c_0 | c_i \rangle = 1$. We know $c_0 - c_1$ is a suture element and $\langle c_0 | c_i \rangle = \pm 1$, $\langle c_1 | c_i \rangle = 1$. Thus $\langle c_0 - c_1 | c_i \rangle = \pm 1 - 1 = 0$ or $-2$. Clearly $-2$ is impossible; hence $\langle c_0 | c_i \rangle = 1$.
\end{Proof}

In fact the proof of part (vi) of theorem \ref{thm:properties_of_suture_elements} above effectively gives another proof of lemma 3.1 of \cite{Me09Paper}. If $\Gamma_u, \Gamma_v$ give a tight contact structure on the solid cylinder then the proof constructs it as a thickened $\Gamma_u$ with a finite set of bypass attachments.

\subsection{Duality--Rotation explicitly}
\label{sec:duality_explicity}

We now prove some more detailed results about the $H=R$ duality/rotation operator, in addition to the results of section \ref{sec:duality}, including the formula of corollary \ref{cor:H_formula}. In \cite{Me09Paper} we obtained, mod $2$, a recursive formula for $H$, an explicit description for $H$. We reprove these here over $\Z$.

Recall the formula of corollary \ref{cor:H_formula}, and recall that $E_w^s$, for $s \in \{x,y\}$, is the set of those $i$ numbering the $s$'s which are followed (to the right) by a different symbol; and for $T \subseteq E_w^x$, we form $\psi_T^x w$ by taking each $x$ in $w$ corresponding to $T$, and the $y$ immediately following it, and replacing each $xy$ with $yx$; similarly for $\psi_T^y w$ where $T \subseteq E_w^y$.
\[
H w = Q_+ Q_-^{-1} w = \sum_{w_i \geq w} \sum_{T \subseteq E_{w_i}^y} (-1)^{|T|} \psi_T^y w_i
\]
Think of $x$'s as pawns on a 1-dimensional chessboard, and $y$'s as empty squares. Then an operation $\psi_T^x w$ moves takes a subset of the pawns and moves them each forward one square each; as they have all moved ahead we have $w \leq \psi_T^x w$. Call any $\psi_T^x w$ \emph{pawn-ahead of} $w$. Similarly, any $\psi_T^y w \leq w$; call any $\psi_T^y w$ \emph{pawn-behind} $w$. Further recall that $f_w^x(i)$ is the number of $y$'s (strictly) to the left of the $i$'th $x$ in $w$. 
\begin{prop}
\label{prop:H_term_by_term}
Let $w \in \M_n^e$. The words $v$ occurring in $Hw$ are precisely those such that:
\begin{enumerate}
\item If $i \in E_v^x$ then $f_v^x (i) = f_w^x (i) - 1$.
\item If $i \notin E_v^x$ then $f_v^x (i) \geq f_w^x$.
\end{enumerate}
The coefficient of $v$ in $Hw$ is $(-1)^{|E_v^x|}$.
\end{prop}

\begin{Proof}
From the formula for $H$, we have immediately that any $v$ occurring in $Hw$ is pawn-behind some $u$ with $w \leq u$. In fact, $v$ appears once for every $u \geq w$ pawn-ahead of $v$. So for $v$ occurring in $Hw$, $f_v^x (i) \geq f_w^x (i)$ for $i \notin E_v^x$, and $f_v^x (i) \geq f_w^x (i) - 1$ for $i \in E_v^x$. Partition the set $\{1, 2, \ldots, n_x\}$ numbering the $x$'s in the words into three sets $S_0, S_1, S_2$:
\begin{align*}
S_0 &= \{1, 2, \ldots, n_x\} - E_v^x  \\
S_1 &= \left\{ i \in E_v^x \; : \; f_v^x (i) = f_w^x (i) - 1 \right\} \\
S_2 &= \left\{ i \in E_v^x \; : \; f_v^x (i) \geq f_w^x (i) \right\}
\end{align*}
Thus, $v$ occurs in $Hw$ once for every word $u \geq w$ pawn-ahead of $v$; such $u$ exist if and only if for each $i \in S_0$ (i.e. those $i \notin E_v^x$), $f_v^x (i) \geq f_w^x (i)$. In this case, the $u$ occurring are precisely those obtained by moving up all the $x$-pawns in $S_1$, and any subset of the $x$-pawns in $S_2$. Thus, the complete set of all terms in the sum for $Hw$ involving $v$ is precisely
\[
\sum_{A \subseteq S_2} (-1)^{|A \cup S_1|} v = (-1)^{|S_1|} v \sum_{A \subseteq S_2} (-1)^{|A|}.
\]
This sum is clearly $0$ if $S_2$ is nonempty; and if $S_2$ is empty then $v$ appears in the final expression for $Hw$ with coefficient $(-1)^{|S_1|}$.

So, the $v$ occurring in $Hw$ are precisely those with $f_v^x (i) \geq f_w^x (i)$ for $i \notin E_v^x$, and such that $S_2 = \emptyset$, i.e. $f_v^x (i) = f_w^x (i) - 1$ for all $i \in E_v^x$. The coefficient of this $v$ is $(-1)^{|S_1|}$, and $S_1 = E_v^x$.
\end{Proof}

We immediately obtain, for instance,
\[
H(x^a) = x^a, \quad H(y^b) = y^b, \quad H( x^a y^b ) = y^b x^a, \quad H( y^b x^a ) = y^b x^a - y^{b-1} x y x^a.
\]
In general, take a word
\[
w = x^{a_1} y^{b_1} x^{a_2} y^{b_2} \cdots x^{a_k} y^{b_k}
\]
where $k \geq 2$, possibly $a_1 = 0$, possibly $b_k = 0$, but all other $a_i, b_i \neq 0$. Now consider a $v$ satisfying the conditions of the above proposition
\[
v = y^{\beta_1} x^{\alpha_1} \cdots y^{\beta_l} x^{\alpha_l}
\]
for some positive integer $l$, where possibly $\beta_1 = 0$, possibly $\alpha_l = 0$ but all other $\alpha_i, \beta_i \neq 0$. The values of $f_v^x$ are given by the cumulative sums $\beta_1 + \cdots + \beta_i$, and the values of $f_w^x$ by cumulative sums $b_1 + \cdots + b_j$. Considering each $x$ counted in $E_v^x$, by condition (i) of proposition \ref{prop:H_term_by_term} we must have
\[
\beta_1 + \cdots + \beta_i = b_1 + \cdots + b_{P(i)} - 1
\]
for $1 \leq i \leq l-1$, for some increasing sequence of integers $1 \leq P(1) < P(2) < \cdots < P(l-1) \leq k$. And of course $\beta_1 + \cdots + \beta_l = b_1 + \cdots + b_k$, so set $P(l) = k$. To maintain the two conditions of the proposition, we must also have the $x$'s ending each block in $v$ corresponding to the $x$ beginning the next block in $w$; thus
\[
\alpha_1 + \cdots + \alpha_i = a_1 + \cdots + a_{P(i)} + 1.
\]
And this term appears with sign $(-1)^{l-1}$. So $Hw$ is the sum of the terms
\[
\begin{array}{c}
(-1)^{l-1} y^{b_1 + \cdots + b_{P(1)} - 1} x^{a_1 + \cdots + a_{P(1)} + 1} y^{b_{P(1) + 1} + \cdots + b_{P(2)}} x^{a_{P(1) + 1} + \cdots + a_{P(2)}} \cdots \\
\cdots
y^{b_{P(l-2) + 1} + \cdots + b_{P(l-1)}} x^{a_{P(l-2) + 1} + \cdots + a_{P(l-1)}}
y^{b_{P(l-1)+1} + \cdots + b_{P(l)} + 1} x^{a_{P(l-1) + 1} + \cdots + a_{P(l)} - 1}
\end{array}
\]
over all $l \leq k$ and all increasing sequences $1 \leq P(1) < P(2) < \cdots < P(l) = k$. There are $2^{k-1}$ terms in the sum. This reproves proposition 7.2 of \cite{Me09Paper}, with signs now included.

\subsection{Duality--Rotation recursively}
\label{sec:duality_recursively}

We can also reprove the recursive properties of $R$ from \cite{Me09Paper}, now with signs. That is, we prove proposition 7.1 of that paper. Note the change of notation from that paper: $A_-, B_-, A_+, B_+$ are now respectively $a_{y,0}, a_{x,0}^*, a_{x,0}, a_{y,0}^*$.
\begin{lem} \
\begin{enumerate}
\item $a_{y,0} H a_{y,0}^* = H$.
\item $a_{y,0} H \left( a_{x,0}^* \right)^j a_{y,0}^* = H \left( a_{x,0}^* \right)^j$.
\item $a_{y,0} \left( a_{x,0} \right)^{j+1} H \left( a_{x,0}^* \right)^j a_{y,0}^* a_{x,0}^* = - H$.
\item
\begin{enumerate}
\item $a_{y,0} \left( a_{x,0} \right)^j H \left( a_{x,0} \right)^j = 0$ for $j = 1, \ldots, n_x$.
\item $\left( a_{x,0} \right)^{j+2} H \left( a_{x,0}^* \right)^j a_{y,0}^* = 0$.
\item $a_{x,0} H \left( a_{x,0}^* \right)^j \left( a_{y,0}^* \right)^2 = 0$.
\end{enumerate}
\end{enumerate}
\end{lem}

\begin{Proof}
Since all creations and annihilations here are initial, we drop the $0$ from the notation for the duration of this proof. For (i), we have
\[
\langle v | a_y H a_y^* w \rangle = \langle a_y^* v | H a_y^* w \rangle = \langle a_y^* w | a_y^* v \rangle = \langle w | v \rangle = \langle v | H w \rangle.
\]
In the first equality we use the adjoint property; then the definition of $H$; then creation operators are isometries; then definition of $H$ again. For (ii), we have
\[
\langle v | a_y H \left( a_x^* \right)^j a_y^* w \rangle = \langle a_y^* v | H \left( a_x^* \right)^j a_y^* w \rangle
= \langle \left( a_x^* \right)^j a_y^* w | a_y^* v \rangle 
= \langle \left( a_x^* \right)^j w | v \rangle
= \langle v | H \left( a_x^* \right)^j w \rangle.
\]
First we use the adjoint property; second the definition of $H$; and fourth the definition of $H$ again. The third equality follows from the observation that the inequality $x^j y w \leq y v$ is true iff $x^j w \leq v$.

Identity (iii) is the most difficult;  we begin by observing that in
\[
\langle v | a_y \left( a_x \right)^{j+1} H \left( a_x^* \right)^j a_y^* a_x^* w \rangle,
\]
any term of $H \left( a_x^* \right)^j a_y^* a_x^* w$ which does not begin with $x^{j+1} y$ is annihilated (to $0$!) by the operators $a_y \left( a_x \right)^{j+1}$; and then for each term $x^{j+1} y u$ occurring, the above inner product outputs $\langle v | u \rangle$. 

Now we note that if $v \leq u$ then $x^{j+1} y v \leq x^{j+1} y u$, and $x^j y x v \nleq x^{j+1} y u$. But, for any word $z$ whose first $j+2$ symbols are not $x^{j+1} y$, we have $x^{j+1} y v \leq z$ iff $x^j y x v \leq z$. Thus
\[
\langle x^{j+1} y v - x^j y x v | z \rangle = \left\{ \begin{array}{cl} \langle v | u \rangle & \text{if $z = x^{j+1} y u$}, \\
0 & \text{otherwise}. \end{array} \right\} 
= \langle v | a_y \left( a_x \right)^{j+1} z \rangle
\]
Hence
\begin{align*}
\langle v | a_y \left( a_x \right)^{j+1} H \left( a_x^* \right)^j a_y^* a_x^* w \rangle
&= \langle \left( a_x^* \right)^{j+1} a_y^* v | H \left( a_x^* \right)^j a_y^* a_x^* w \rangle
- \langle \left( a_x^* \right)^j a_y^* a_x^* v | H \left( a_x^* \right)^j a_y^* a_x^* w \rangle \\
&= \langle \left( a_x^* \right)^j a_y^*  a_x^* w | \left( a_x^* \right)^{j+1} a_y^* v \rangle
- \langle \left( a_x^* \right)^j a_y^* a_x^* w | \left( a_x^* \right)^j a_y^* a_x^* v \rangle \\
&= - \langle w | v \rangle = - \langle v | H w \rangle,
\end{align*}
where we first use the previous observation, then definition of $H$, then notice $x^jyx \nleq x^{j+1} y$ and that creations are isometries, and finally definition of $H$ again. This proves (iii).

We do not need to prove (iv) algebraically. As in \cite{Me09Paper}, it is sufficient to observe that the operators on chord diagrams produce closed loops, so are $0$.
\end{Proof}

As in \cite{Me09Paper}, this lemma describes the matrix for each $H_{n_x,n_y}$ recursively. Following notation there, for any two words $w_0, w_1$ (with $\leq n_x$ $x$'s and $\leq n_y$ $y$'s), we define the $w_0 \times w_1$ minor of this matrix to be the intersection of the rows corresponding to all words beginning with $w_0$ with the columns corresponding to all words beginning with $w_1$. The above lemma then gives the following description of $H_{n_x,n_y}$:
\begin{enumerate}
\item The $y \times y$ minor consists of $H_{n_x,n_y-1}$.
\item The $y \times xy$ minor contains the $x$-columns of $H_{n_x, n_y - 1}$. More generally, the $y \times x^j y$ minor contains the $x^j$-columns of $H_{n_x, n_y - 1}$, for any $j=1, \ldots, n_x$.
\item The $xy \times yx$ minor consists of $- H_{n_x - 1, n_y - 1}$. More generally, for any $j = 0, \ldots, n_x -1$, the $x^{j+1} y \times x^j y x$ minor consists of $- H_{n_x - j - 1, n_y - 1}$.
\item All other entries are zero.
\end{enumerate}

We may write this recursive structure as a formula, as in \cite{Me09Paper}: the ``fake'' commutator (commutators do not mean much mod 2!) there becomes a real commutator here.
\begin{thm}
\label{thm:H_recursive_formula}
\begin{align*}
H &= \sum_{i=0}^\infty 
a_{y,0}^* H \left( a_{x,0}^* \right)^i a_{y,0} \left( a_{x,0} \right)^i
-
\left( a_{x,0}^* \right)^{i+1} a_{y,0}^* H a_{x,0} a_{y,0} \left( a_{x,0} \right)^i \\
&= \sum_{i=0}^\infty 
 \left[ a_{y,0}^* \; H \; a_{x,0}, \; \left( a_{x,0}^* \right)^{i+1} \right] a_{y,0} \left( a_{x,0} \right)^i.
\end{align*}\
\qed
\end{thm}

\subsection{Periodicity}
\label{sec:periodicity}

We now prove theorem \ref{thm:H_periodic}, that $H$ is periodic, $H^{n+1} = (-1)^{n_x n_y}$. In fact we show that not only does $H$ rotate through chord diagrams, it also in a certain sense rotates through basis elements.

Consider a word $w \in \M_{n_x, n_y}$ as an $n \times 1$ chessboard, where $x$'s are pawns and $y$'s are empty squares. Form a sequence of words / chessboard configurations, which begins at $w_1$ with all pawns / $x$'s at the extreme right, i.e. $w_1 = y^{n_y} x^{n_x}$. Pawns can move one square left, if the square to their left is empty. Starting from $w_1$, we form a sequence $w_i$: on each move we move all possible pawns one square left, until we arrive at $w_n$ with all pawns at the extreme left, $w_n = x^{n_x} y^{n_y}$. Thus, for example, with $n_x = 2$ and $n_y = 3$ we have
\[
w_1 = yyyxx, w_2 = yyxyx, w_3 = yxyxy, w_4 = xyxyy, w_5 = xxyyy
\]
Set $w_0 = w_n = x^{n_x} y^{n_y}$ and number the $w_i$ mod $n+1$. (So $x^{n_x} y^{n_y}$ appears twice in a cycle.) We have obviously $w_1 > w_2 > \cdots > w_n = w_0$.

As above, denote by $\pm [u, v]$ denotes the two suture elements with first word $u$ and last word $v$ ($u \leq v$). For now adopt the notation that $[u,v]$ is the one of this pair of suture elements in which $u$ has coefficient $1$. Consider the suture elements $[w_i, w_{i-1}]$: it must take one of the following four forms
\[
[w_i, w_{i-1}] = x^\alpha (xy - yx)^\beta x^\gamma, \; x^\alpha (xy - yx)^\beta y^\gamma, \; y^\alpha (xy-yx)^\beta x^\gamma, \; \text{ or } y^\alpha (xy-yx)^\beta y^\gamma.
\]
To see why, note that all the above expressions are products of suture elements, hence suture elements; expanding them out we obtain a sum where the lexicographically first is some $w_i$, and the last is some $w_{i-1}$, for some $i$ and some $n_x, n_y$; for every $n_x, n_y$ and $i$, we can obtain such an expression with $w_i$ and $w_{i-1}$ occurring first and last; and we have proved there is a unique pair of suture elements $\pm [u,v]$ beginning and ending with each possible pair $u \leq v$.

Examining the explicit formula for $R=H$, we see that in any word $w$ beginning $w = x^\alpha (xy)^k (yx) \cdots$ or $y^\alpha (xy)^k (yx) \cdots$, where $k \geq 0$, the minimum word occurring in $H w$ begins:
\begin{align*}
\min H \left( x^\alpha (xy)^k (yx) \cdots \right) &= \left\{ \begin{array}{ll} x^{\alpha + 1} (xy)^k yx \cdots & k \geq 2 \\ y x^{\alpha + 2} \cdots & k = 1 \\ x^{\alpha + 1} y \cdots & k = 0 \end{array} \right. \\
\min H \left( y^\alpha (xy)^k (yx) \cdots \right) &= \left\{ \begin{array}{ll} y^{\alpha - 1} (xy)^k yx \cdots & k \geq 1 \\ y^\alpha x \cdots & k = 0 \end{array} \right.
\end{align*}
On the other hand, in $H w_i$, the (lexicographically) first word occurring is $w_{i+1}$; it occurs with coefficient $(-1)^{|E_{w_{i+1}}^x|}$. Thus, for any word $w \neq w_i$ occurring in $[w_i, w_{i-1}]$, the first word occurring in $H w$ is lexicographically after $w_{i+1}$. Thus $H ( [w_i, w_{i-1}] ) = (-1)^{|E_{w_{i+1}}^x|} [ w_{i+1}, v ]$ for some word $v$.

We can see what this $v$ is directly by examining chord diagrams. The basis chord diagrams $\Gamma_{w_i}$ are easily constructed, and so are forwards bypass systems taking each $\Gamma_{w_i}$ to $\Gamma_{w_{i-1}}$, as described in section 5 of \cite{Me09Paper}. Then perform downwards bypass surgeries along this bypass, following notation of \cite{Me09Paper}, we obtain a chord diagram $[ \Gamma_{w_i}, \Gamma_{w_{i-1}} ]$ such that, when its suture element is expanded in terms of basis elements, the first word occurring is $w_i$ and the last is $w_{i-1}$. We then note that $H [ \Gamma_{w_i}, \Gamma_{w_{i-1}} ] = [ \Gamma_{w_{i+1}}, \Gamma_{w_i} ]$, i.e. rotating the chord diagram by $H$ takes us from one to the next.  Hence
\[
H \left( \left[ w_i, w_{i-1} \right] \right) = (-1)^{|E_{w_{i+1}}^x|} \left[ w_{i+1}, w_i \right].
\]
It follows immediately that
\begin{align*}
  H^{n+1} \left( x^{n_x} y^{n_y} \right) &= H^{n+1} \left( \left[ x^{n_x} y^{n_y}, x^{n_x} y^{n_y} \right] \right) = H^{n+1} \left( \left[ w_n, w_0 \right] \right) \\
  &= \left( \prod_{i=0}^n (-1)^{|E_{w_i}^x|} \right) \left[ w_n, w_0 \right] = (-1)^{\sum_{i=0}^n |E_{w_i}^x|} x^{n_x} y^{n_y}.
\end{align*}
Recall that $E_w^x$ is the set of all block-ending $x$'s in $w$ (except possibly for the final $x$, i.e. all $x$'s which are followed by $y$'s. Now as we run through the set $w_i$, and pawns move from right to left, we see that for every $x$ and every $y$, they are adjacent to each other as $xy$ precisely once. Thus $\sum |E_{w_i}^x|$ counts every pair of an $x$ and a $y$ exactly once, and is equal to $n_x n_y$. So $H^{n+1} x^{n_x} y^{n_y} = (-1)^{n_x n_y} x^{n_x} y^{n_y}$.

\begin{Proof}[of theorem \ref{thm:H_periodic}]
We have shown above that for $w = w_{min} = x^{n_x} y^{n_y}$, $H^{n+1} w = (-1)^{n_x n_y} w$. Let now $A = \{ w \in \M_{n_x,n_y} : H^{n+1} w = (-1)^{n_x n_y} \}$, so $w_{min} \in A$. That $H^{n+1} = (-1)^{n_x n_y}$ follows immediately from the following fact:
\begin{itemize}
\item If $w \in A$ and $w, w'$ are related by an elementary move, then $w' \in A$.
\end{itemize}
To see this, let $H^{n+1} w' = \epsilon w'$ where $\epsilon = \pm 1$. As $w,w'$ are related by an elementary move, $w-w'$ is a suture element, and $H^{n+1} (w-w') = \pm (w- w')$. But we have $H^{n+1}(w-w') = (-1)^{n_x n_y} w - \epsilon w' = (-1)^{n_x n_y} \left( w - (-1)^{n_x n_y} \epsilon w' \right)$. Hence we must have $(-1)^{n_x n_y} \epsilon = 1$, so $\epsilon = (-1)^{n_x n_y}$ and $w' \in A$. 

Finally, it's clear from the discussion above that the least positive $j$ for which $H^j(x^{n_x} y^{n_y}) = \pm x^{n_x} y^{n_y}$ is $j=n+1$; and we now have $H^{n+1} = (-1)^{n_x n_y}$. Thus $H$ has period $n+1$ if $n_x n_y$ is even, and period $2n+2$ otherwise.
\end{Proof}

\addcontentsline{toc}{section}{References}

\small

\bibliography{danbib}
\bibliographystyle{amsplain}

\end{document}